\documentstyle[12pt]{article} 
\begin{document}
\newcommand{\singlespace}{
    \renewcommand{\baselinestretch}{1}
\large\normalsize}
\newcommand{\doublespace}{
   \renewcommand{\baselinestretch}{1.2}
   \large\normalsize}
\renewcommand{\theequation}{\thesection.\arabic{equation}}

\input amssym.def
\input amssym
\setcounter{equation}{0}
\def \ten#1{_{{}_{\scriptstyle#1}}}
\def \Z{\Bbb Z}
\def \C{\Bbb C}
\def \R{\Bbb R}
\def \Q{\Bbb Q}
\def \N{\Bbb N}
\def \l{\lambda}
\def \V{V^{\natural}}
\def \wt{{\rm wt}}
\def \tr{{\rm tr}}
\def \Res{{\rm Res}}
\def \End{{\rm End}}
\def \Aut{{\rm Aut}}
\def \mod{{\rm mod}}
\def \Hom{{\rm Hom}}
\def \im{{\rm im}}
\def \<{\langle} 
\def \>{\rangle} 
\def \w{\omega}
\def \o{\omega}
\def \t{\tau }
\def \char{{\rm char}}
\def \a{\alpha }
\def \b{\beta}
\def \e{\epsilon }
\def \la{\lambda }
\def \om{\omega }
\def \O{\Omega}
\def \qed{\mbox{ $\square$}}
\def \pf{\noindent {\bf Proof: \,}}
\def \voa{vertex operator algebra\ }
\def \voas{vertex operator algebras\ }
\def \p{\partial}

\singlespace
\newtheorem{thmm}{Theorem}
\newtheorem{th}{Theorem}[section]
\newtheorem{prop}[th]{Proposition}
\newtheorem{coro}[th]{Corollary}
\newtheorem{lem}[th]{Lemma}
\newtheorem{rem}[th]{Remark}
\newtheorem{de}[th]{Definition}

\begin{center}
{\Large {\bf Modular-invariance of trace functions in orbifold theory}} \\
\vspace{0.5cm}

Chongying Dong\footnote{Supported by NSF grant 
DMS-9303374 and faculty research funds granted by the University of 
California at 
Santa Cruz.}\\
Department of Mathematics, University of
California, Santa Cruz, CA 95064 \\
Haisheng Li\footnote{Supported by NSF grant DMS-9616630.}\\
Department of Mathematical Sciences,
Rutgers University, Camden, NJ 08102\\
Geoffrey Mason\footnote{Supported by NSF grant DMS-9401272 and faculty research funds granted by the University of 
California at 
Santa Cruz.}
\\
Department of Mathematics, University of
California, Santa Cruz, CA 95064
\end{center}
\hspace{1.5 cm}

\begin{abstract} The goal of the present paper is to provide a mathematically
rigorous foundation to certain aspects of the theory of rational orbifold
models in conformal field theory, in other words the theory of rational
vertex operator algebras and their automorphisms.

Under a certain finiteness condition on a rational vertex operator algebra
$V$ which holds in all known 
examples, we determine the precise number of $g$-twisted sectors for 
any automorphism $g$ of $V$ of finite order. We prove that
the trace functions and correlation functions associated with such 
twisted sectors are holomorphic functions in the upper half-plane and,
under suitable conditions, afford a representation of the modular
group of the type prescribed in string theory. We establish the rationality
of conformal weights and central charge.

In addition to conformal field theory itself, where our conclusions are 
required on physical grounds, there are applications to the generalized
Moonshine conjectures of Conway-Norton-Queen and to equivariant elliptic
cohomology.
\end{abstract}

Table of Contents

         1. Introduction

         2. Vertex Operator Algebras

         3. Twisted Modules

         4. $P$-functions and $Q$-functions

         5. The Space of 1-point Functions on the Torus

         6. The Differential Equations

         7. Formal 1-point Functions

         8. Correlation Functions

         9. An Existence Theorem for $g$-twisted Sectors

         10. The Main Theorems

         11. Rationality of Central Charge and Conformal Weights

         12. Condition $C_2$

         13. Applications to the Moonshine Module

\section{Introduction}

The goals of the present paper are to give a mathematically rigorous
study of rational orbifold models, more precisely we study the
questions of the {\em existence} and {\em modular-invariance of
twisted sectors of rational vertex operator algebras.} The idea of
orbifolding a vertex operator algebra with respect to an automorphism,
and in particular the introduction of twisted sectors, goes back to
some of the earliest papers in the physical literature on conformal
field theory [DHVW], while the question of modular-invariance
underlies the whole enterprise. Apart from a few exceptions such as
[DGM], the
physical literature tends to treat the existence and
modular-invariance of twisted sectors as axioms, while mathematical work
has been mainly limited to studying special cases such as
affine algebras [KP], lattice theories [Le], [FLM1] and
fermionic orbifolds [DM2].  Under some mild finiteness conditions on a 
rational vertex operator algebra V we will, among other things,
establish the following:

      (A) The precise number of inequivalent, simple g-twisted sectors that
            V possesses.

      (B) Modular-invariance (in a suitable sense) of the characters of
           twisted sectors.

   In order to facilitate the following discussion we assume that the
reader has a knowledge of the basic definitions concerning vertex operator
algebras as given, for example, in [FLM2], [FHL], [DM1] and below in
Sections 2 and 3.

   Let us suppose that $V$ is a vertex operator algebra. There are several
approaches to what it means for $V$ to be {\em rational}, each of them referring
to finiteness properties of $V$ of various kinds (c.f. [MS], [HMV], [AM] for
example). Our own approach is as follows (c.f. [DLM2], [DLM3]): an
{\em admissible}  $V$-module is a certain $\Z_{+}$-graded linear space
\begin{equation}\label{e1.1}
 M  =\bigoplus_{n=0}^{\infty}M(n)
\end{equation}
which admits an action of $V$ by vertex operators which satisfy certain
axioms, the most important of which is the Jacobi identity. However, the
homogeneous subspaces $M(n)$ are not assumed to be finite-dimensional. As
explained in [DLM2], this
definition includes as a special case the idea of an  ({\em ordinary})
V-module, which is a graded linear space of the shape
\begin{equation}\label{e1.2}
                M =\bigoplus_{n\in \C}M_{n}
\end{equation}
such that $M$ admits an appropriate action of $V$ by vertex operators and such
that each $M_{n}$ has finite dimension, $M_{\la+n}=0$ for fixed $\la\in \C$
and sufficiently small integer $n,$ and each $M_{n}$
 is the  $n$-eigenspace of the $L(0)$-operator. Although one
ultimately wants to establish the rationality of the grading of such
modules, it turns out to be convenient to allow the gradings to be 
{\em a  priori}
more general. One then has to {\em prove}  that the grading is rational after
all. By the way, the terms
``simple module'' and ``sector'' are synonymous in this context.

We call $V$ a {\em rational vertex operator algebra}  in case each admissible
$V$-module is completely reducible, i.e., a direct sum of simple admissible
modules. We have proved in [DLM3]
 that this assumption implies that $V$ has only finitely many
inequivalent simple admissible modules, moreover each such module is in
fact an ordinary simple $V$-module. These results are special cases of
results proved in (loc.cit.) in which one considers the same set-up, but
relative to an automorphism $g$ of $V$ of finite order. Thus one has the
notions of $g$-{\em twisted}  $V$-{\em module, admissible} $g$-{\em
twisted} $V$-{\em module} and
$g$-{\em rational vertex operator algebra} $V$, 
the latter being a vertex operator algebra all of
whose admissible g-twisted modules are completely reducible. We do not give
the precise definitions here (c.f. Section 3), noting only that if $g$ has
order $T$ then a simple $g$-twisted $V$-module has a grading of the shape
\begin{equation}\label{e1.3} 
M=\bigoplus_{n=0}^{\infty}M_{\lambda+ n/T}
\end{equation}
with $M_{\l}\ne 0$ 
for some complex number $\l$ which is called the {\em conformal weight}
 of $M.$ This
is an important invariant of $M$ which plays a r$\hat{\rm o}$le in the theory of the Verlinde algebra, for example.

    The basic result (loc.cit.) is that a $g$-rational vertex operator
algebra has only finitely many inequivalent, simple, admissible $g$-twisted
modules, and each of them is an ordinary simple $g$-twisted module. Although
the theory of twisted modules includes that of ordinary modules,   which
corresponds to the case $g = 1,$  it is nevertheless common and convenient to
refer to the {\em untwisted theory}  if $g = 1,$ and to the 
{\em twisted theory} otherwise.

   It seems likely that if $V$ is rational then in fact it is $g$-rational for
all automorphisms $g$ of finite order, but this is not known. Nevertheless,
our first main result shows that there is a close relation between the
untwisted and twisted theories. To explain this, we  first recall Zhu's
notion of ``Condition $C$" [Z], which we have renamed as Condition $C_2.$
Namely, for an element $v$ in the vertex operator algebra $V$ we use standard
notation for the vertex operator determined by $v$:
\begin{equation}\label{e1.4}
Y( v, z) =\sum_{n\in\Z}v(n) z^{-n-1}
\end{equation}
so that each $v(n)$ is a linear operator on $V.$ Define $C_2(V)$ to be
the subspace of $V$ spanned by the elements $u(-2)v$ for $u, v$ in
$V.$ We say that $V$ satisfies {\em Condition} $C_2$ if $C_2(V)$ has finite
codimension in $V.$ Zhu has conjectured in [Z] that every rational
vertex operator algebra satisfies Condition $C_2,$ and we verify the
conjecture for a number of the most familiar rational vertex operator
algebras in Section 12 of the present paper.  Next, any automorphism h
of $V$ induces in a natural way (c.f. Section 4 of [DM2] and Section 3 below) a bijection from the set of (simple,
admissible) $g$-twisted modules to the corresponding set of 
$hgh^{-1}$-twisted modules, so that if $g$ and $h$ 
commute then one may consider the
set of $h$-stable $g$-twisted modules. In particular we have the set of
$h$-stable ordinary (or untwisted) modules, which includes the vertex
operator algebra $V$ itself. Our first set of results may now be stated
as follows:
\begin{thmm}\label{thmm1}
Suppose that $V$ is a rational vertex operator algebra which
satisfies Condition $C_2.$ Then the following hold:

   (i) The central charge of $V$ and the conformal weight of each simple
$V$-module are rational numbers.

 (ii) If $g$ is an automorphism of $V$ of finite order then the number of
inequivalent, simple $g$-twisted $V$-modules is at most
equal to the
number of $g$-stable simple untwisted $V$-modules, and is at least 1 if $V$ is
simple.

(iii)  If $V$ is $g$-rational, the number of
inequivalent, simple $g$-twisted $V$-modules is precisely the
number of $g$-stable simple untwisted $V$-modules.
\end{thmm}

  We actually prove a generalization of this result in which $V$ is also
assumed to be $g^i$-rational for all integers $i$ ($g$ as in (ii)), where the
conclusion is that each simple $g^i$-twisted $V$-module has rational conformal
weight. There is a second variation of this theme involving the important
class of {\em holomorphic}  vertex operator algebras. This means that $V$ is
assumed to be both simple and rational, moreover $V$ is assumed to have a
{\em unique}  simple module - which is necessarily the {\em adjoint}
  module $V$ itself.
Familiar examples include the Moonshine Module [D2] and vertex operator
algebras associated to positive-definite, even, unimodular lattices [D1].
We establish
\begin{thmm}\label{thmm2}
Suppose that $V$ is a holomorphic vertex operator algebra which
satisfies Condition $C_2,$ and that $g$ is an automorphism of $V$ of finite
order. Then the following hold:

   (i) $V$ posesses a {\em unique}
  simple $g$-twisted $V$-module, call it $V(g).$

  (ii) The central charge of $V$ is a positive integer divisible by $8,$ and
the conformal weight of $V(g)$ is a rational number.
\end{thmm}

    The assertion of Theorem \ref{thmm2} (ii) concerning the central 
charge follows from results in [Z] as a consequence of 
the modular-invariance of the
character of a holomorphic vertex operator algebra, and we turn now to a
discussion of the general question of modular-invariance. One is concerned
with various trace functions, the most basic of which is the formal
character  of a (simple) $g$-twisted sector $M.$ If $M$ has grading 
(\ref{e1.3}) we define the formal character as
\begin{equation}\label{e1.5}
\char_qM=q^{\lambda-c/24}\sum_{n=0}^{\infty}\dim M_{\lambda + n/T} q^{n/T}
\end{equation}
where $c$ is the central charge and $q$ a formal variable. More generally, 
if $M$ is an $h$-stable $g$-twisted sector as before, then $h$ induces a 
linear map on $M$
which we denote by $\phi(h),$ and one may consider the corresponding graded
trace
\begin{equation}\label{e1.6}
Z_M(g, h) =  q^{\la-c/24}\sum_{n=0}^{\infty}\tr_{M_{\lambda+n/T}}\phi(h)q^{n/T}
\end{equation}
In the special case when $V$ is holomorphic and the twisted sector $V(g)$ is
unique, we set
\begin{equation}\label{e1.7}
   Z_{V(g)}(g,h) = Z(g, h)
\end{equation}

Note that in this situation, the uniqueness of $V(g)$ shows that $V(g)$ is
$h$-stable whenever $g$ and $h$ commute. As discussed in [DM1]-[DM2], this allows us to consider $\phi$  as a
{\em projective}  representation of the {\em centralizer}
  $C(g)$ on $V(g),$ in particular
(1.7) is defined for all commuting pairs $(g, h).$
   As is well-known, it is important to consider these trace functions as
special cases of so-called $(g,h)$ {\em correlation  functions.}  These may be
defined for homogeneous elements $v$ in $V_k$ and any pair of commuting
elements $(g, h)$ via
\begin{equation}\label{e1.8}
T_M(v, g, h)  =  q^{\lambda-c/24}\sum_{n=0}^{\infty}\tr_{M_{\lambda+n/T}}(v(k-1) \phi(h))q^{n/T}.
\end{equation}
Note that $v(k-1)$ induces a linear map on each homogeneous subspace of $M.$ If
we take $v$ to be the vacuum vector then (\ref{e1.8}) reduces to (\ref{e1.6}).
   One also considers these trace functions less formally by taking $q$ to be
the usual local parameter at infinity in the upper half-plane
\begin{equation}\label{e1.9}
{\frak h} =  \{\tau\in\C| \im\tau> 0 \}
\end{equation}
i.e., $q = q_{\tau}=e^{2\pi i \tau},$ in which case one writes 
$T_M(v, g, h, \tau)$ etc.
for the corresponding trace functions. By extending $T_M$ linearly to the
whole of $V$ one obtains a function
\begin{equation}\label{e1.10}
T_M: V \times P(G) \times {\frak h}\to P^1(\C)
\end{equation}
where $P(G)$ is the set of commuting pairs of elements in $G.$

  We take  $\Gamma=SL(2,\Z)$  to operate on ${\frak h}$ in the usual way via bilinear
transformations, that is
\begin{equation}\label{e1.11}
\gamma: \tau\mapsto\frac{ a\tau + b}{c\tau+d},\ \ \ \gamma=\left(\begin{array}{cc}a & b\\ c & d\end{array}\right)\in\Gamma
\end{equation}
and we let it act on the right of $P(G)$ via
\begin{equation}\label{e1.12}
         (g, h)\gamma =(g^a h^c, g^b h^d).
\end{equation}

   Zhu has introduced in [Z] a second vertex operator algebra associated in
a certain way to $V;$ it has the same underlying space, however the grading
is different. We are concerned with elements $v$ in $V$ which are homogeneous
of weight $k,$ say, with regard to the second  vertex operator algebra. We
write this as $\wt [v] = k.$ For such $v$ we define an action of the modular
group   $\Gamma$ on $T_M$ in a familiar way, namely
\begin{equation}\label{e1.13}
T_M|\gamma(v, g, h, \tau)  =  (c\tau+ d)^{-k} T_M (v, g, h,  \gamma\tau).
\end{equation}

   We can now state some of our main results concerning modular-invariance.

\begin{thmm}\label{thmm3}
Suppose that $V$ is a vertex operator algebra which satisfies
Condition $C_2,$ and let $G$ be a finite group of automorphisms of $V.$

   (i) For each triple $(v, g, h )\in V \times  P(G)$ and for each $h$-stable
$g$-twisted sector $M,$ the trace function  $T_M(v,g,h,\tau)$  is holomorphic 
in ${\frak h}.$

  (ii) Suppose in addition that $V$ is $g$-rational for each $g\in G.$ 
Let $v\in V$ satisfy $\wt[v] = k.$ Then the space of (holomorphic) functions 
in ${\frak h}$ spanned
by the trace functions $T_M(v, g, h, \tau)$ for all choices of $g, h,$ and $M$
 is  a (finite-dimensional)  $\Gamma$-module with respect to the action 
(\ref{e1.13}).

   More precisely, if $\gamma\in \Gamma$ then we have an equality
\begin{equation}\label{e1.14}
T_M|\gamma(v, g, h,\tau)=\sum_{W}\sigma_WT_W(v,(g,h)\gamma,\tau),
\end{equation}
where $(g, h)\gamma$ is as in (\ref{e1.12}) and  $W$ ranges over the 
$g^a h^c$-twisted sectors which are $g^b h^d$-stable. 
The constants $\sigma_W$ depend only on $g, h,$ $\gamma$ and $W.$
\end{thmm}
\begin{thmm}\label{thmm4}
Suppose that $V$ is a holomorphic vertex operator algebra
which satisfies Condition $C_2,$ and let $G$ be a cyclic group of automorphisms
of $V.$ If $(g, h)\in P(G),$ $v\in V$ satisfies $\wt[v] = k,$ and 
$\gamma\in \Gamma$ then $Z(v, g, h,\tau)$ is a holomorphic function in 
${\frak h}$  which satisfies
\begin{equation}\label{e1.15}
             Z|\gamma(v, g, h,\tau)  = 
\sigma(g, h, \gamma) Z(v, (g, h)\gamma, \tau)
\end{equation}
for some constant $\sigma(g, h, \gamma).$
\end{thmm}

   We can summarize some of the results above by saying that 
the functions $T_M(v, g, h, \tau)$  and  $Z(v, g, h, \tau)$  are 
{\em generalized modular forms}
 of weight $k$ in the sense of [KM]. This means essentially that each
of these functions and each of their transforms under the modular group
have $q$-expansions in {\em rational}  powers of $q$ with bounded denominators, and
that up to scalar multiples there are only a finite number of such
transforms. One would like to show that each of these functions is, in
fact, a modular form of weight $k$ and some level $N$  in the usual sense 
of being invariant under the principal congruence subgroup  
$\Gamma(N).$ This will require further argument, as it is shown in (loc.cit.)
 (and is
certainly known to the experts) that there are generalized modular forms
which are {\em not} modular forms in the usual sense. All we can say in general
at the moment is that  $T_M$  and  $Z$  have finite level in the sense of
Wohlfahrt [Wo]. This is true of all generalized modular forms, and means that
$T_M$ and $Z$ and each of their  $\Gamma$ -transforms are invariant under the 
group $\Delta(N)$ for some $N,$ where we define  $\Delta(N)$ to be the normal 
closure of $T^N=\left(\begin{array}{cc} 1 & N\\ 0 & 1\end{array}\right)$
in $\Gamma.$

   Let us emphasize the differences between Theorem \ref{thmm3} (ii) and 
Theorem \ref{thmm4}.
In the former we assume $g$-rationality, which in practice is hard to verify,
even for known vertex operator algebras. In Theorem \ref{thmm4} there is no 
such
assumption, however we have to limit ourselves to {\em cyclic}
  pairs $(g, h)$ in
$P(G).$ One expects Theorem \ref{thmm4} to hold for all commuting pairs 
$(g, h).$ In [N],
Simon Norton conjectured that (\ref{e1.15}) holds in the special case that 
$v$ is
the vacuum vector and $V$ is the Moonshine Module [FLM2] whose automorphism
group is the Monster (loc.cit.). His argument was based on extensive
numerical evidence in Conway-Norton's famous paper Monstrous Moonshine
[CN] which was significantly expanded in the thesis of Larrissa Queen [Q].
A little later it was given a string-theoretic interpretation in [DGH]. We
will see in Section 12 that the Moonshine Module satisfies Condition $C_2,$
so that Theorem \ref{thmm4} applies. Norton also conjectured that each 
$Z(g, h, \tau)$ is
either constant or a {\em hauptmodul,} the latter being a modular function
(weight zero) on some congruence subgroup  $\Gamma_0$ of $\Gamma$ of genus zero such
that the modular function in question generates the full field of
meromorphic modular functions on  $\Gamma_0.$  By utilizing the results of Borcherds
[B2] which establish the original Moonshine conjectures in [CN] we obtain

\begin{thmm}\label{thmm5} (Generalized Moonshine): Let $V^{\natural}$
 be the Moonshine Module, and
let $g$ be an element of the Monster simple group $M.$ The following hold:

   (i) $V^{\natural}$ has a unique g-twisted sector $V^{\natural}(g).$

  (ii) The formal character $\char V^{\natural}(g)$ is a hauptmodul.

 (iii) More generally, if $g$ and $h$ in $M$ generate a cyclic subgroup then the
graded trace $Z(g, h, \tau)$ of  $\phi(h)$  on $V^{\natural}(g)$ is a 
hauptmodul.
\end{thmm}

   These results essentially establish the Conway-Norton-Queen
conjectures about the Monster for cyclic pairs $(g, h).$ As we have
said, the spaces $V^{\natural}(g)$ support faithful projective
representations of the corresponding centralizer $C_M(g)$ of $g$ in
$M.$ It was these representations that were conjectured to exist in
[N] and [Q], and they are of considerable interest in their own
right. See [DLM1] for more information in some special cases.  There
also appears to be some remarkable connections with the work of
Borcherds and Ryba [Ry], [BR], [B3] on {\em modular moonshine} which we hope
to consider elsewhere.  Finally, we mention that Theorem \ref{thmm4}
can be considerably strengthened, indeed the best possible results can
be established, if we assume that $g$ has small order. For
simplicity we state only a special case:
\begin{thmm}\label{thmm6} Let $V$ be as in Theorem \ref{thmm4}, 
and assume that the central
charge of $V$ is divisible by 24 (cf. Theorem \ref{thmm2} (ii)). Suppose 
that $g$ has
order $p = 2$ or 3. Then the following hold:

   (i) The conformal weight  $\lambda$  of the (unique) $g$-twisted sector 
$V(g)$ is in $\frac{1}{p^2}\Z.$

  (ii) The graded trace $Z(v, 1, g, \tau)$ of $v(k-1)g$ on $V$ is a modular form of
weight $k$ and level $p^2.$

 (iii) We have $\lambda\in\frac{1}{p}\Z$ if, and only if, $Z(v, g, h, \tau)$
 is a modular form on
the congruence subgroup  $\Gamma_0(p).$
\end{thmm}

   These results follow from the previous theorems, and will be 
proved elsewhere.  They can be used to make rigorous some of the assumptions
commonly made by physicists (cf. [M], [S], [T1], [T2], [T3], [V] for
example) in the theory of $\Z_p$-orbifolds, and we hope that Theorem 
\ref{thmm6} may
provide the basis for a more complete theory of such orbifolds.

   We have already mentioned the work of Zhu [Z] on several occasions, and
it is indeed this paper to which we are mainly indebted intellectually. In
essence, we are going to prove an equivariant version of the theory laid
down by Zhu (loc.cit.), though even in the special case that he was
studying our work yields improvements on his results. The equivariant
refinement of Zhu's theory began with our paper [DLM3]
 which also plays a basic r$\hat{\rm o}$le in the present paper. In
this paper we constructed so-called twisted Zhu algebras $A_g(V)$ which are
associative algebras associated to a vertex operator algebra $V$ and
automorphism $g$ of finite order. They have the property that, at least for
suitable classes of vertex operator algebras, the module category for $A_g(V)$
and the category of $g$-twisted modules for $V$ are {\em equivalent.}
 This reduces
the construction of $g$-twisted sectors to the corresponding problem for
$A_g(V)$ (not a priori known to be non-zero!) As in [Z], the r$\hat{\rm o}$le of the
finiteness condition $C_2$ is to show that the $(g, h)$ correlation functions
that we have considered above satisfy certain differential equations of
regular - singular type. These differential equations have coefficients
which are essentially modular forms on a congruence subgroup, a fact which
is ultimately attributable to the Jacobi identity satisfied by the vertex
operators. One attempts to characterize the space of correlation functions
as those solutions of the differential equation which possess other
technical features related to properties of $V$ and $A_g(V).$ This space is
essentially what is sometimes referred to as the 
$(g, h)$-{\em conformal block},
and our results follow from the technical assertion that, under suitable
circumstances, the $(g, h)$-conformal block is indeed spanned by the 
$(g, h)$-correlation functions.
   
We point out that the holomorphy of trace functions follows from the
fact that they are solutions of suitable differential equations, moreover
the Frobenius - Fuchs theory of differential equations with regular -
singular points leads to the representation of elements of the conformal
block as $q$-expansions. One attempts to identify coefficients of such
$q$-expansions with the trace function defined by some $A_g(V)$-module, and at
the same time show that elements of the conformal blocks are free of
logarithmic singularities. The point is that the Frobenius - Fuchs theory
plays a critical r$\hat{\rm o}$le, as it does also in [Z].

   The paper is organized as follows: after some preliminaries in
Sections 2 and 3, we take up in Section 4 the study of certain modular
and Jacobi-type forms, the main goal being to write down the
transformation laws which they satisfy. Our methods (and probably the
results too) will be well-known to experts in elliptic functions, but
it is fascinating to see how such classical topics such as Eisenstein
series and Bernoulli distributions play a r$\hat{\rm o}$le in an
abstract theory of vertex operator algebras. In Section 5 we introduce
the space of abstract $(g, h)$ 1-point functions associated to a
vertex operator algebra and establish that it affords an action by the
modular group (Theorem 5.4). In Sections 6 and 7 we continue the study
of such functions and in particular write down the differential
equation that they satisfy and the general shape of the solutions in
terms of q-expansions and logarithmic singularities. Next we prove in
Section 8 (Theorems 8.1 and 8.7) that if $g$ and $h$ commute then
distinct $h$-stable $g$-twisted sectors give rise to linearly
independent trace functions which lie in the $(g, h)$-conformal block,
and in Section 9 we give as an application of the ideas developed so
far a general existence theorem for twisted sectors. Section 10
contains the main theorems which give sufficient conditions under
which the $(g, h)$ conformal block is spanned by trace
functions. Having reached Section 11, we have enough information to be
able to apply the methods and results of Anderson and Moore [AM], and
this leads to the rationality results stated above in Theorems
\ref{thmm1} and \ref{thmm2} as well as the applications to
modular-invariance and to the generalized Moonshine Conjectures, which
are discussed in Section 13.  We also point out how one can use Theorem
\ref{thmm4} to describe  not only other correlation functions but also 
``Monstrous Moonshine of weight k.''
Section 12 establishes
condition $C_2$ for a number of well-known rational vertex operator
algebras, so that our theory applies to all of these examples.

\section{Vertex operator algebras}
\setcounter{equation}{0}
The definition of vertex operator algebra entails
a $\Z$-graded complex vector space:
\begin{equation}\label{g2.1}
V=\coprod_{n\in{\Bbb Z}}V_n
\end{equation}
satisfying $\dim V_{n}< \infty$ for all $n$ and  $V_{n}=0$ for $n << 0.$
If $v\in V_n$ we write $\wt v=n$ and say that $v$ is {\em homogeneous} and has 
(conformal) {\em weight} $n.$ For each $v\in V$ there are linear operators
$v_n\in\End V,$ $n\in\Z$ which are assembled into a {\em vertex operator}
$$Y(v,z)=\displaystyle{\sum_{n\in\Z}v(n)z^{-n-1}} \in  (\mbox{End}\,V)[[z,z^{-1}]].$$ 
Various axioms are imposed:

(i) For $u,v\in V,$ 
\begin{equation}\label{g2.2}
u(n)v=0\ \ \ \ \ \mbox{for}\ \  n\ \ \mbox{sufficiently large}.
\end{equation}

(ii) There is a distinguished {\em vacuum element} ${\bf 1}\in V_0$ satisfying
\begin{equation}\label{g2.3}
Y({\bf 1},z)=1
\end{equation}
and
\begin{equation}\label{g2.4}
Y(v,z){\bf 1}=v+\sum_{n\geq 2}v(-n){\bf 1}z^{n-1}.
\end{equation}

(iii) There is a distinguished {\em conformal vector} $\o\in V_2$ with
generating function 
$$Y(\w,z)=\sum_{n\in\Z}L(n)z^{-n-2}$$
such that the
component operators generate a copy of the Virasoro algebra represented on $V$
with {\em central charge} $c.$ That is
\begin{equation}\label{g2.5}
[L(m),L(n)]=(m-n)L(m+n)+\frac{1}{12}(m^3-m)\delta_{m+n,0}c.
\end{equation}
Moreover we have
\begin{eqnarray}
V_n=\{v\in V|L(0)v=nv\}\label{g2.6}\\
\frac{d}{dz}Y(v,z)=Y(L(-1)v,z).\label{g2.7}
\end{eqnarray}

(iv) The {\em Jacobi identity} holds, that is 
\begin{equation}\label{2.2}
\begin{array}{c}
\displaystyle{z^{-1}_0\delta\left(\frac{z_1-z_2}{z_0}\right)
Y(u,z_1)Y(v,z_2)-z^{-1}_0\delta\left(\frac{z_2-z_1}{-z_0}\right)
Y(v,z_2)Y(u,z_1)}\\
\displaystyle{=z_2^{-1}\delta
\left(\frac{z_1-z_0}{z_2}\right)
Y(Y(u,z_0)v,z_2)}
\end{array}
\end{equation}
where 
$\delta(z)=\sum_{n\in {\Bbb Z}}z^n$ and  $\delta(\frac{z_1-z_2}{z_0})$ is 
expanded as a formal power series in $z_2.$

Such a vertex operator algebra may be denoted  by the 4-tuple 
$(V,Y,{\bf 1},\omega)$ or, more usually, by $V$. 

Zhu has introduced a second vertex operator algebra
$(V,Y[\ ],{\bf 1},\tilde\o)$ associated to $V$ in Theorem 4.2.1 of [Z].
 It plays a crucial r$\hat{\rm o}$le
in Zhu's theory and also in the present paper, so we give some details
of the construction\footnote{Our formulae differ from that of Zhu by a factor
of $2\pi i$ in the exponent $e^{z \wt v}.$}.

The conformal vector $\tilde \o$ is defined to be     
$\omega-\frac{c}{24}.$ The vertex operators $Y[v,z]$  are defined
for homogeneous $v$ via the equality  
\begin{equation}\label{g2.8}
Y[v,z]=Y(v,e^z-1)e^{z{\rm wt}v}=\sum_{n\in\Z}v[n]z^{-n-1}.
\end{equation}

For integers $i,m,p$ with $i,m\geq 0$ we may define $c(p,i,m)$ in either
of two equivalent ways:
\begin{eqnarray}
{p-1+z\choose i}=\sum_{m=0}^ic(p,i,m)z^m\\ \label{g2.9}
m!\sum_{i=m}^{\infty}c(p,i,m)z^i=(\log(1+z))^m(1+z)^{p-1}\label{g2.10}
\end{eqnarray}
where, as usual
\begin{equation}\label{g2.11}
\log(1+z)=\sum_{n=1}^{\infty}\frac{(-1)^{n-1}}{n}z^n.
\end{equation}
Using a change of variable we calculate that 
\begin{eqnarray*}
& &v[m]=\Res_zY[v,z]z^m\\
& &\,\ \ \ \ \ \ =\Res_zY[v,\log(1+z)](\log(1+z))^m(1+z)^{-1}
\end{eqnarray*}
i.e.,
\begin{equation}\label{g2.12}
v[m]={\rm Res}_zY(v,z)(\log(1+z))^m(1+z)^{{\rm wt}v-1}.
\end{equation}
So if $m\geq 0$ then
\begin{equation}\label{g2.13}
v[m]=m!\sum_{i=m}^\infty c(\wt v,i,m)v(i).
\end{equation}
For example we have
\begin{equation}\label{g2.14}
v[0]=\sum_{i=0}^{\infty}{\wt v-1\choose i}v(i).
\end{equation}

We also write 
\begin{equation}\label{g2.15}
Y[\tilde \o,z]=\sum_{n\in\Z}L[n]z^{-n-2}.
\end{equation}
For example, one has
\begin{eqnarray}
& &L[-2]=\o[-1]-\frac{c}{24}\label{g2.16}\\
& & L[-1]=L(-1)+L(0)\label{g2.17}\\
& &\ \ \,L[0]=L(0)+\sum_{n=1}^{\infty}\frac{(-1)^{n-1}}{n(n+1)}L(n).\label{g2.18}
\end{eqnarray}

Care must be taken to distinguish between the notion of conformal weight
in the {\em original} vertex operator algebra and in 
$(V,Y[\ ],{\bf 1},\tilde\o).$ If $v\in V$ is homogeneous in the {\em latter} 
vertex operator algebra we denote its conformal weight by $\wt [v],$ and
set 
\begin{equation}\label{g2.19}
V_{[n]}=\{v\in V|L[0]v=nv\}.
\end{equation}

In general, $V_n$ and $V_{[n]}$ are distinct, though it follows from
(\ref{g2.18}) that for each $N$ we have
\begin{equation}\label{g2.20}
\bigoplus_{n\leq N}V_n=\bigoplus_{n\leq N}V_{[n]}.
\end{equation}

\section{Twisted modules}
\setcounter{equation}{0}

Let $V$ be a vertex operator algebra. 
An {\em automorphism} $g$ of the vertex operator algebra $V$ is a
linear automorphism of $V$ preserving ${\bf 1}$ and $\omega$ such that
the actions of $g$ and $Y(v,z)$ on $V$ are compatible in the sense
that $$gY(v,z)g^{-1}=Y(gv,z)$$ for $v\in V.$  Let $\Aut V$ be the group
of automorphisms of $V.$ If $g\in \Aut V$ has finite order $T,$ say, there
are various classes of $g$-twisted $V$-modules of concern to us
(cf. [DLM2], [DLM3]). A {\em weak} $g$-twisted $V$-module 
is a $\C$-linear space $M$ equipped with a linear map
$V\to (\End M)[[z^{1/T},z^{-1/T}]],$ $v\mapsto Y_M(v,z)=\sum_{n\in\Q}v(n)z^{-n-1}$ satisfying the following:

For $v\in V, w\in M,$ 
\begin{equation}\label{g3.10}
v(m)w=0
\end{equation}
if $m$ is sufficiently large;
\begin{equation}\label{g3.11}
Y_M({\bf 1},z)=1;
\end{equation}
set 
\begin{equation}\label{d1}
V^r=\{v\in V|gv=e^{-2\pi ir/T}v\}
\end{equation}
for  $0\leq r\leq T-1.$ Then
\begin{eqnarray}
Y_M(v,z)=\sum_{n\in r/T+\Z}v(n)z^{-n-1}\ \ \ \ {\rm for}\ \ 
v\in V^r;\label{g3.12}
\end{eqnarray}
and the twisted Jacobi identity holds
\begin{equation}\label{g3.13}
\begin{array}{c}
\displaystyle{z^{-1}_0\delta\left(\frac{z_1-z_2}{z_0}\right)
Y_M(u,z_1)Y_M(v,z_2)-z^{-1}_0\delta\left(\frac{z_2-z_1}{-z_0}\right)
Y_M(v,z_2)Y_M(u,z_1)}\\
\displaystyle{=z_2^{-1}\left(\frac{z_1-z_0}{z_2}\right)^{-r/T}
\delta\left(\frac{z_1-z_0}{z_2}\right)
Y_M(Y(u,z_0)v,z_2)}
\end{array}
\end{equation}
if $u\in V^r.$ We often write $(M,Y_M)$ for this module. It can be shown
(cf. Lemma 2.2 of [DLM2], [DLM3]) that $Y_M(\o,z)$ has component operators
which still satisfy (\ref{g2.5}) and (\ref{g2.7}). If we take $g=1,$ we get a weak 
$V$-module.

For $m\in \C$ we use the notation $\iota_{z_1,z_2}(z_1-z_2)^m$ to denote
the  expansion of $(z_1-z_2)^m$
in terms the nonnegative integral powers of $z_2:$
\begin{equation}\label{ad1}
\iota_{z_1,z_2}((z_1-z_2)^m)=z_1^m\sum_{i\geq 0}{m\choose i}(-1)^iz_2^iz_1^{-i}
\end{equation}
where ${m\choose i}=\frac{m(m-1)\cdots (m-i+1)}{i!}.$ 	It follows from
the twisted Jacobi identity
(\ref{g3.13}) that
\begin{eqnarray}
& &\Res_{z_1}Y_M(u,z_1)Y_M(v,z_2)z_1^r\iota_{z_1,z_2}((z_1-z_2)^m)z_2^n
\nonumber\\
& &\ \ \ -\Res_{z_1}Y_M(v,z_2)Y_M(u,z_1)z_1^r\iota_{z_2,z_1}((z_1-z_2)^m)
z_2^n\nonumber\\
& &=\Res_{z_1-z_2}Y_M(Y(u,z_1-z_2)v,z_2)\iota_{z_2,z_1-z_2}z_1^r(z_1-z_2)^mz_2^n\label{ad2}
\end{eqnarray} 
for $m\in\Z, n\in\C.$

An {\em admissible} $g$-twisted $V$-module 
is a  weak $g$-twisted $V$-module $M$ which carries a 
${1\over T}{\Z}_{+}$-grading 
\begin{equation}\label{m3.12}
M=\oplus_{n\in\frac{1}{T}\Z_{+}}M(n)
\end{equation}
which satisfies the following
\begin{eqnarray}\label{m3.13}
v(m)M(n)\subseteq M(n+\wt v-m-1)
\end{eqnarray}
for homogeneous $v\in V.$ We may assume that $M(0)\ne 0$ if $M\ne 0.$
Again if $g=1$ we have an admissible $V$-module.

An (ordinary) $g$-twisted $V$-module is a weak $g$-twisted $V$-module $M$ which carries a 
$\C$-grading induced by the spectrum of $L(0).$ That is, we have
\begin{equation}\label{g3.14}
M=\coprod_{\lambda \in{\C}}M_{\lambda} 
\end{equation}
where $M_{\l}=\{w\in M|L(0)w=\l w\}.$ Moreover we require that 
$\dim M_{\l}$ is finite and for fixed $\l,$ $M_{{n\over T}+\l}=0$
for all small enough integers $n.$ If $g=1$ we have an ordinary $V$-module.
If $M$ is a simple $g$-twisted $V$-module, then
\begin{equation}\label{3.10'}
M=\bigoplus_{n=0}^{\infty}M_{\la+n/T}
\end{equation}
for some $\l\in\C$ such that $M_{\l}\ne 0.$ 
 We define $\l$ as the {\em conformal weight} of  $M.$

The \voa $V$ is called $g$-{\em rational} in case every admissible $g$-twisted
$V$-module is completely reducible, i.e., a direct sum of simple admissible
$g$-twisted modules.

\begin{rem} There is a subtlety in the definition of these twisted
modules. Namely, the definition of $V^r$ given in (\ref{d1}) is not
quite the same as that which we have used elsewhere ([DLM3], [DM1]
etc.).  Previously we set $V^r=\{v\in V|gv=e^{2\pi ir/T}v\}$, so in
effect we have interchanged the notions of $g$- and $g^{-1}$-twisted modules.
The reason for doing this is that it makes the theorem of modular-invariance 
(Theorem \ref{t5.4} below) have the expected form.
\end{rem}

\begin{th}\label{dt3.3} \ [DLM3] Suppose that $V$ is a $g$-rational vertex 
operator algebra, where $g\in\Aut V$
has finite order. Then the following hold:

(a) Every simple admissible $g$-twisted $V$-module is an ordinary $g$-twisted
$V$-module.

(b) $V$ has only finitely many isomorphism classes of simple admissible
 $g$-twisted modules.
\end{th}

Let $V$ be a \voa with $g$ an automorphism of $V$ of finite order $T.$
In [DLM3] a certain associative algebra $A_g(V)$ was introduced which
plays a basic r$\hat{\rm o}$le in the present work. In case $g=1,$ $A_1(V)=A(V)$
was introduced and used extensively in [Z]. We need to review
certain aspects of these constructions here.

Let $V^r$ be as in (\ref{d1}). For $u\in V^r$ and $v\in V$ we define 
\begin{eqnarray}\label{g3.15}
u\circ_g v=\Res_{z}\frac{(1+z)^{{\wt u}-1+\delta_{r}+{r\over T}}}{z^{1+\delta_{r}}}Y(u,z)v
\end{eqnarray}
\begin{equation}\label{g3.16}
u*_gv=\left\{
\begin{array}{ll}
\Res_z(Y(u,z)\frac{(1+z)^{{\wt}\,u}}{z}v)
 & {\rm if}\ r=0\\
0  & {\rm if}\ r>0.
\end{array}\right.
\end{equation} 
where $\delta_{r}=1$ if $r=0$ and 
$\delta_{r}=0$ if $r\ne 0$.
Extend $\circ_g$ and $*_g$ to bilinear products  on $V.$ We let $O_g(V)$ be the linear span of all $u\circ_g v.$
\begin{th}\label{t3.4} ([DLM3],[Z]) The quotient $A_g(V)=V/O_g(V)$
is an associative algebra with respect to $*_g.$
\end{th}

Note that if $g=1$ then $A(V)$ is nonzero, whereas if $g\ne 1$ the
analogous assertion may not be true. But if $A_g(V)\ne 0$ then the
vacuum element maps to the identity of $A_g(V),$ and the conformal
vector maps into the center of $A_g(V).$

Let $M$ be a weak $g$-twisted $V$-module. 
Define $\O(M)=\{w\in V|u(\wt u-1+n)w=0,
u\in V, n>0\}.$ 
For homogeneous $u\in V$ define
\begin{equation}\label{g3.18}
o(u)=u(\wt u-1)
\end{equation}
(sometimes called the {\em zero mode} of u).
\begin{th}\label{t3.5} [DLM3] Let $M$ be a weak $g$-twisted $V$-module. The following
hold:

(a) The map $v\mapsto o(v)$ induces a representation of the associative algebra
$A_g(V)$ on $\O(M).$

(b) If $M$ is a simple admissible $g$-twisted $V$-module then $\Omega(M)=M(0)$
is a simple $A_g(V)$-module. Moreover, $M\mapsto M(0)$
induces a bijection between (isomorphism classes of) simple
admissible $g$-twisted $V$-modules and simple $A_g(V)$-modules.
\end{th}

When combined with Theorem \ref{dt3.3} one finds
\begin{th}\label{t3.6} [DLM3] Suppose that $V$ is a \voa with an 
automorphism $g$ of finite order, and that $V$ is $g$-rational 
(possibly $g=1$). Then the following
hold:

(a) $A_g(V)$ is a finite-dimensional, semi-simple associative algebra
(possibly 0).

(b) The map $M\mapsto \O(M)$ induces an equivalence between the
category of ordinary $g$-twisted $V$-modules and the category of
finite-dimensional $A_g(V)$-modules.
\end{th}

There are various group actions that we need to explain. 
Let $g,h$ be automorphisms of $V$ with $g$ of finite order.
If $(M,Y_g)$ is a weak $g$-twisted module for $V$ there is a weak 
$hgh^{-1}$-twisted $V$-module  
$(M, Y_{hgh^{-1}})$ where for $v\in V$ we define
\begin{equation}\label{g3.19}
Y_{hgh^{-1}}(v,z)=Y_g(h^{-1}v,z).
\end{equation}
This defines a left action of $\Aut(V)$ on the family of all weak
twisted modules. Symbolically, we write
\begin{equation}\label{g3.20}
h\circ (M,Y_g)=(M,Y_{hgh^{-1}})=h\circ M.
\end{equation}
The action (\ref{g3.20}) induces an action
\begin{equation}\label{g3.21}
h\circ \O(M)=\O(h\circ M).
\end{equation}

Next, it follows easily from the definitions (\ref{g3.15}) and (\ref{g3.16}) that
the action of $h$ on $V$ induces an isomorphism of associative algebras
\begin{equation}\label{g3.22}
\begin{array}{cccc}
h: & A_g(V) & \to & A_{hgh^{-1}}(V)\\
   & v& \mapsto & hv
\end{array}
\end{equation}
which then induces a  map
\begin{equation}\label{g3.23}
h: A_{hgh^{-1}}(V)-\mod\to A_{g}(V)-\mod.
\end{equation}
To describe (\ref{g3.23}), let $(N,*_{hgh^{-1}})$ be a left
$A_{hgh^{-1}}(V)$-module (extending the notation of
(\ref{g3.16})). Then $h\circ(N,*_{hgh^{-1}})=(N,*_{g})$ where, for $n\in
N, v\in V,$
\begin{equation}\label{g3.24}
v*_{hgh^{-1}}n=h^{-1}v*_gn.
\end{equation}

Now if $(M,Y_g)$ is as before then (\ref{g3.21}) and (\ref{g3.23}) both define
actions of  $h$ on $\O(M);$ they are the same. For if $v\in V$ and we
consider the image of $v$ in $A_{hgh^{-1}}(V),$ it acts on $\O(h\circ M)$ via
the zero mode $o_{hgh^{-1}}(v)$ of $v$ in the vertex operator 
$Y_{hgh^{-1}}(v,z)=Y_g(h^{-1}v,z).$ In other words, if $m\in\O(h\circ M)=\O(M)$ then
$o_{hgh^{-1}}(v)m=o_g(h^{-1}v)m,$ which is precisely what (\ref{g3.24}) says.

We can summarize this by saying that the functor $\O$ is 
$\Aut(V)$-equivariant, so that in the previous notation there is a commuting
diagram (up to isomorphism)
\begin{equation}\label{3.25}
\begin{array}{ccc}
\{{\rm simple}\ g{\rm -twisted}\ V{\rm -modules}\} & \stackrel{h}{\rightarrow} &\{{\rm simple}\ hgh^{-1}{\rm -twisted}\  V{\rm -modules}\}\\
\O\downarrow & & \downarrow\O\\
\{{\rm simple}\ A_g(V){\rm -modules}\} & \stackrel{h}{\rightarrow} &\{{\rm simple}\ A_{hgh^{-1}}(V){\rm -modules}\}
\end{array}
\end{equation}

We say that the $g$-twisted $V$-module $M$ is $h$-invariant if $h\circ M\cong
M.$ The set of all such automorphisms, the {\em stablizer} of $M,$ is a subgroup $C$ of $\Aut V$ which contains $g.$ This is because one knows [DLM3] that
$g$ acts trivially on $A_g(V),$ whence the assertion follows from 
(\ref{3.25}). Moreover $g\in Z(C)$ (the {\em center} of $C$). There is
a {\em projective} representation of $C$ on $M$ induced by the action 
(\ref{g3.20}). See [DM1] or [DM2] for more information on this point.

Next we discuss the so-called $C_2$-{\em condition} introduced by Zhu [Z].
Let $V$ be a \voa and $M$ a $V$-module. We define
\begin{equation}\label{3.26}
C_2(M)=\<v(-2)m|v\in V, m\in M\>,
\end{equation}
that is, $C_2(M)$ is spanned by the indicated elements of $M.$ We say that
$M$ {\em satisfies condition} $C_2$ in case $C_2(M)$ has finite codimension
in $M.$ The most important case is that in which $M$ is taken to
be $V$ itself.

\begin{prop}\label{p2.2} Suppose that $V$ satisfies condition $C_2,$
and let $g$ be an automorphism of $V$ of finite order. Then the 
algebra $A_g(V)$ has finite dimension. 
\end{prop}

\noindent{\bf Proof:\ } Note from (\ref{3.26}) that $C_2(V)$ is a
$\Z$-graded subspace of $V.$ Since the codimension of $C_2(V)$ is finite,
there is an integer $k$ such that
$V=C_2(V)+W$ where $W$ is the sum of the first $k$ homogeneous subspaces of $V.$

We will show that $V_m\subset W+O_g(V)$ for each $m\in\Z,$ in which case
$V=W+O_g(V)$ and $\dim A_g(V)\leq \dim W.$ We proceed by induction on
$m.$

Recall that $V^r$ is the eigenspace of $g$ with eigenvalue $e^{-2\pi ir/T}$
for $0\leq r\leq T-1$ where $g$ has order $T.$ Since $C_2(V)$ is a
homogeneous and $g$-invariant subspace of $V$ then 
we may write any $c\in C_2(V)\cap V_m$ in the form
\begin{equation}\label{3.27}
c=\sum_{i=1}^nu_i(-2)v_i
\end{equation}
for homogeneous elements $u_i,v_i\in V$ satisfying $u_i\in V^r$
for some $r=r(i)$ and $\wt u_i+\wt v_i+1=m.$ 

Suppose first that $u_i\in V^r$ with $r=0.$ According to (\ref{g3.15}),
$O_g(V)$ contains
\begin{equation}\label{3.28}
\Res_{z}Y(u_i,z)\frac{(1+z)^{\wt u_i}}{z^{2}}v_i=\sum_{j\geq 0}{\wt u_i\choose
j}u_i(j-2)v_i.
\end{equation}
Now $\wt u_i(j-2)v_i=\wt u_i+\wt v_i-j+1=m-j,$ so if $j\geq 1$ then
$u_i(j-2)v_i\in W+O_g(V)$ by induction. But then $W+O_g(V)$ also
contains the remaining summand $u_i(-2)v_i$ of (\ref{3.28}).

Now suppose that $u_i\in V^r$ with $r\geq 1.$ By Lemma 2.2 (i) of
[DLM3] we have that $O_g(V)$ contains the element
$\Res_{z}Y(u_i,z)\frac{(1+z)^{\wt u_i-1+r/T}}{z^{2}}v_i,$ and we conclude
once again that $u_i(-2)v_i$ lies in $W+O_g(V).$ So we have shown that
all summands of (\ref{3.27}) lie in $W+O_g(V).$
The proposition follows.\qed

\begin{prop}\label{p3.8} Suppose that $V$ satisfies condition $C_2,$ and
 let $g$ be an automorphism of $V$ of finite order. If $A_g(V)\ne 0$ then
$V$ has a simple $g$-twisted $V$-module.
\end{prop}

\pf $A_g(V)$ has finite dimension by Proposition \ref{p2.2}. Now the
result follows from Theorem 9.1 of [DLM3]. \qed

The following lemma will be used in Section 12.
\begin{lem}\label{l3.7a} Let $M$ be a $V$-module. Then
$C_2(M)$ is invariant under the operators
$v(0)$ and $v(-1)$ for any $v\in V.$
\end{lem}

\pf Consider $u(-2)w\in C_2(M)$ for $u\in V$ and $w\in M.$ Then for
$k=0,-1$
$v(k)u(-2)w=u(-2)v(k)w+\sum_{i=0}^{\infty}{k\choose i}(v(i)u)(-2+k-i)w\in C_2(M)$
as required. \qed

\section{$P$-functions and $Q$-functions}
\setcounter{equation}{0}

We study certain functions, which we denote by $P$ and $Q,$ which play a ${\rm r}\hat{\rm o}{\rm le}$ 
in later sections. The $P$-functions are essentially Jacobi forms [EZ] and
the $Q$-functions are certain modular forms. The main goal is to write
down the transformation laws of these functions under the
action of the modular group $\Gamma=SL(2,\Z).$

Let $\frak h$ denote the upper half plane $\frak h=\{z\in\C|{\rm im}\,z>0\}$
with the usual left action of $\Gamma$ via  M\"obius 
transformations
\begin{equation}\label{g4.1}
\left(\begin{array}{cc}
a & b\\
c & d
\end{array}
\right):\tau\mapsto \frac{a\tau+b}{c\tau+d}.
\end{equation}
$\Gamma$ also acts on the right of $T^2=S^1\times S^1$ via
\begin{equation}\label{g4.2}
(\mu,\lambda)
\left(
\begin{array}{cc}
a & b\\
c & d
\end{array}
\right)
=(\mu^a\lambda^c,\mu^b\lambda^d).
\end{equation}

Let $t$ be a torsion point of $T^2.$ Thus $t=(\mu,\l)$ with
  $\mu=e^{2\pi ij/M}$ and $\lambda=e^{2\pi il/N}$ for 
integers $j,l,M,N$ with $M,N>0.$ For each integer $k=1,2,\cdots$ and
each $t$ we define a function $P_k$ on $\C\times \frak h$ as follows
\begin{equation}\label{g4.3}
P_k(\mu,\lambda,z,q_{\tau})=P_k(\mu,\lambda,z,\tau)=
\frac{1}{(k-1)!}\sum^{\ \ \ \ \ \prime}_{n\in \frac{j}{M}+\Z}\frac{n^{k-1}q_z^n}{1-\l q_{\tau}^n}
\end{equation}
where the sign $\sum'$ means omit the term $n=0$ if $(\mu,\l)=(1,1).$
Here and below we write $q_x=e^{2\pi ix}.$

We will prove
\begin{th}\label{t4.1} Suppose that $(\mu,\l)\ne (1,1).$ Then
$$P_k(\mu,\lambda,\frac{z}{c\tau+d},\gamma\tau)=
(c\tau+d)^kP_k((\mu,\lambda)\gamma,z,\tau).$$
for all $\gamma=
\left(
\begin{array}{cc}
a & b\\
c & d
\end{array}
\right)\in \Gamma.$
\end{th}

\begin{rem} (i) (\ref{g4.3}) converges uniformly and absolutely on compact subsets of the 
region $|q_{\tau}|<|q_z|<1.$

(ii) Theorem \ref{t4.1} holds also for $(\mu,\la)=(1,1)$ in case $k\geq 3$
but {\em not} if $k=1,2$ (cf. [Z]).
\end{rem}

We can reformulate Theorem \ref{t4.1} as follows:
for suitable functions $F(\mu,\lambda,z,\tau)$ on $(\Q/\Z)^2\times \C\times
\frak h,$ and for an integer $k,$ we set 
\begin{equation}\label{g4.4}
F|_k\gamma(\mu,\lambda,z,\tau)=(c\tau+d)^{-k}F((\mu,\lambda)\gamma^{-1},
\frac{z}{c\tau+d},\gamma\tau).
\end{equation}
As is well-known, this defines a right action of $\Gamma$ on such functions $F.$ Theorem \ref{t4.1} says precisely that $P_k$ 
is an {\em invariant} of this action. So it is enough to prove the
theorem for the two standard generators 
$S=\left(
\begin{array}{cc}
0 & -1\\
1 & 0
\end{array}
\right)$ and
$T=\left(
\begin{array}{cc}
1 & 1\\
0 & 1
\end{array}
\right)$
of $\Gamma.$ If $\gamma=T$ then Theorem \ref{t4.1} reduces to the assertion 
$P_k(\mu,\lambda,z,\tau+1)=P_k(\mu,\mu\lambda,z,\tau),$ which
follows 
immediately from the definition (\ref{g4.3}).

We also note the equality
\begin{equation}\label{g4.5}
\frac{d}{dz}P_k(\mu,\lambda,z,\tau)=2\pi ikP_{k+1}(\mu,\lambda,z,\tau).
\end{equation}
So if Theorem \ref{t4.1} holds for $k$ then it holds for $k+1$ by (\ref{g4.5})
and the chain rule. These reductions reduce us to proving Theorem \ref{t4.1}
in the case that $k=1$ and $\gamma=S,$ when it can be restated in the
form
\begin{th}\label{t4.2}
If $(\mu,\lambda)\ne (1,1)$ then 
$$P_1(\mu,\lambda,\frac{z}{\tau},\frac{-1}{\tau})=\tau P_1(\lambda,\mu^{-1},z,\tau).$$
\end{th}

We will need to make use of several other functions in the proof of Theorem
\ref{t4.2}.

First there is the usual Eisenstein series $G_2(\tau)$ with
$q$-expansion
\begin{equation}\label{g4.6}
G_2(\tau)=\frac{\pi^2}{3}+2(2\pi i)^2\sum_{n=1}^{\infty}\frac{nq_{\tau}^n}{1-q_\tau^n}
\end{equation}
and well-known transformation law
\begin{equation}\label{g4.7}
G_2(\gamma\tau)=(c\tau+d)^2G_2(\tau)-2\pi ic(c\tau+d),\ \ 
\gamma=\left(
\begin{array}{cc}
a & b\\
c & d
\end{array}
\right)\in \Gamma.
\end{equation}
(For these and other facts about elliptic functions, the reader may consult 
[La].) Let 
\begin{equation}\label{g4.8}
\wp_1(z,\tau)=G_2(\tau)z+\pi i\frac{q_z+1}{q_z-1}
+2\pi i\sum_{n=1}^{\infty}\left(\frac{q_\tau^nq_z^{-1}}{1-q_{\tau}^nq_z^{-1}}
- -\frac{q_zq_{\tau}^n}{1-q_zq_{\tau}^n}\right).
\end{equation}
The function $\wp_1(z,\tau)$ is not elliptic, but satisfies 
\begin{equation}\label{g4.9}
\wp_1(z+1,\tau)=\wp_1(z,\tau)+G_2(\tau)
\end{equation}
\begin{equation}\label{g4.10}
\wp_1(z+\tau,\tau)=\wp_1(z,\tau)+G_2(\tau)\tau-2\pi i
\end{equation}
\begin{equation}\label{g4.11}
\wp_1(\frac{z}{c\tau+d},\gamma\tau)=(c\tau+d)\wp_1(z,\tau),
\gamma\in\left(\begin{array}{cc}
a & b\\
c& d
\end{array}
\right)\in\Gamma.
\end{equation}

Now introduce a further $P$-type function
\begin{equation}\label{g4.12}
P_{\lambda}(z,\tau)=2\pi i\sum^{\ \ \ \ \ \prime}_{n\in \Z}\frac{q_z^n}{1-\l q_{\tau}^n}
\end{equation}
where $\lambda$ is a root of unity.

\begin{lem}\label{l4.3}
Suppose that $|q_\tau|<|q_z|<1$ and that $\la^N=1.$  Then
$$P_{\la}(z,\tau)=\sum_{k=0}^{N-1}\la^k\left(
G_2(N\tau)(z+k\tau)-\wp_1(z+k\tau,N\tau)-\pi i\right).$$
\end{lem}

\noindent{\bf Proof \ } As $|\lambda q_\tau^n|<1$ for $n\geq 1$ then
\begin{eqnarray*}
& & \sum_{n=1}^{\infty}\frac{q_z^n}{1-\la q_{\tau}}
=\sum_{n=1}^{\infty}\sum_{m=0}^{\infty}q_z^n\la^mq_{\t}^{mn}\\
& &\ \ \ \ \ \ =\sum_{m=0}^{\infty}\sum_{n=0}^{\infty}\la^mq_zq_{\t}^m(q_zq_{\t}^m)^n\\
& &\ \ \ \ \ \ =\sum_{m=0}^{\infty}\frac{\la^mq_zq_{\t}^m}{1-q_zq_{\t}^m}
\ \ \ \ \ \ \ \ ({\rm as}\ |q_zq_{\t}^m|<1\ {\rm for}\ m\geq 0)\\
& &\ \ \ \ \ \ =\frac{q_z}{1-q_z}+\sum_{m=1}^{\infty}\frac{\la^mq_zq_{\t}^m}{1-q_zq_{\t}^m}.
\end{eqnarray*}
Using $|q_z^{-1}q_{\t}^m|<1$ for $m\geq 1,$ a similar calculation
yields 
$$\sum_{n=1}^{\infty}\frac{\la^{-1}q_z^{-n}q_{\t}^n}{1-\la^{-1}q_{\t}^n}
=\sum_{m=1}^{\infty}\frac{\la^{-m}q_z^{-1}q_{\t}^m}{1-q_z^{-1}q_{\t}^m}.$$
{}From this and (\ref{g4.12}) we get
\begin{equation}\label{g4.13}
(2\pi i)^{-1}P_{\la}(z,\t)=
\frac{q_z}{1-q_z}+\sum_{m=1}^{\infty}\left(\frac{\la^mq_zq_{\t}^m}{1-q_zq_{\t}^m}-\frac{\la^{-m}q_z^{-1}q_{\t}^m}{1-q_z^{-1}q_{\t}^m}\right).
\end{equation}
Next, use the expansion
$$\sum_{m=0}^{\infty}\frac{\la^mq_zq_{\t}^m}{1-q_zq_{\t}^m}=
\sum_{n=0}^{\infty}\sum_{k=0}^{N-1}\frac{\la^kq_zq_{\t}^{nN+k}}{1-q_zq_{\t}^{nN+k}}
 =\sum_{k=0}^{N-1}\la^k\sum_{n=0}^{\infty}\frac{q_{z+k\t}q_{N\t}^{n}}{1-q_{z+k\t}q_{\t}^{n}}$$
and a similar expression for the second term under the summation sign in (\ref{g4.13}) to see that
\begin{equation}\label{g4.14}
P_{\la}(z,\t)
=2\pi i\sum_{k=0}^{N-1}\la^k\left(\sum_{n=0}^{\infty}\frac{q_{z+k\t}q_{N\t}^{n}}{1-q_{z+k\t}q_{N\t}^{n}}-\sum_{n=1}^{\infty}\frac{q_{z+k\t}^{-1}q_{N\t}^{n}}{1-q_{z+k\t}^{-1}q_{N\t}^{n}}\right).
\end{equation}
Using the formula (\ref{g4.8}) for $\wp_1(z+k\t,N\t),$ the lemma
follows readily from (\ref{g4.14}). \qed

\bigskip

For a root of unity $\la,$ set
\begin{equation}\label{g4.15}
\epsilon(\la)=\left\{
\begin{array}{lr}
1-\la, & \la \ne 1\\
0, & \la=1.
\end{array}\right.
\end{equation}

Now we are ready for the proof of Theorem \ref{t4.2}. Let $\nu=e^{2\pi i/M}$
with $\mu$ and $\la$ as before. For 
$t\in\Z$ we then have
\begin{eqnarray*}
\sum_{j=1}^M\nu^{jt}P_1(\nu^j,\la,z,\t)
=\sum_{j=1}^M\nu^{jt}\sum^{\ \ \ \ \ \prime}_{n\in \frac{j}{M}+\Z}\frac{q^n_z}{1-\l q_{\tau}^n}
=\sum^{\ \ \ \ \ \prime}_{n\in \Z}\frac{q^n_{\frac{z+t}{M}}}{1-\l q_{\tau\over M}^n}.
\end{eqnarray*}
{}From (\ref{g4.12}) and (\ref{g4.15}) we conclude that
\begin{equation}\label{g4.16}
\sum_{j=1}^M\nu^{jt}P_1(\nu^j,\la,z,\t)={1\over 2\pi i}P_{\l}({z+t\over M},
{\tau\over M})+\epsilon(\l).
\end{equation}

Regarding this as a system of linear equations in 
$P_1(\nu^j,\la,z,\t)$ for $t=0,1,..., M-1,$ we may invert to find
that (with $\mu=\nu^j$)
$$P_1(\mu,\la,z,\t)
=\frac{1}{2\pi iM}\sum_{t=0}^{M-1}\mu^{-t}\left(P_1(\frac{z+t}{M},\frac{\t}{M})
+2\pi i\epsilon(\la)\right).$$
Now use Lemma \ref{l4.3} to obtain
\begin{eqnarray}
& &\ \ \ \ P_1(\mu,\lambda,z,\t)\nonumber\\
& &=\frac{1}{2\pi iM}\sum_{t=0}^{M-1}\sum_{k=0}^{N-1}\mu^{-t}\la^k
\left(G_2(\frac{N\t}{M})\frac{(t+k\t)}{M}-\wp_1(\frac{z+t+k\t}{M},\frac{N\t}{M})\right)\nonumber\\
& &\ \ \ \ \ \ \ \ +\frac{\epsilon(\la)}{M}\sum_{t=0}^{M-1}\mu^{-t}.\label{g4.17}
\end{eqnarray}
(We used $(\mu,\l)\ne(1,1)$ to eliminate some terms in (\ref{g4.17}).)

So now we get, using (\ref{g4.7}) and (\ref{g4.17}):
\begin{eqnarray}
& &\ \ \ \ P_1(\mu,\la,\frac{z}{\t},\frac{-1}{\t})\nonumber\\
& &=\frac{1}{2\pi iM}\sum_{t=0}^{M-1}\sum_{k=0}^{N-1}\mu^{-t}\la^k
\left(G_2(\frac{-N}{M\t})\frac{(-k+t\t)}{M\t}-\wp_1(\frac{z-k+t\t}{M\t},\frac{-N}{M\t})\right)\nonumber\\
& &\ \ \ \ \ \ \ \ +\frac{\epsilon(\la)}{M}\sum_{t=0}^{M-1}\mu^{-t}\nonumber\\
& &=\!\frac{1}{2\pi iM}\sum_{t=0}^{M-1}\sum_{k=0}^{N-1}\mu^{-t}\la^k
\!\left(\!\left(\!\left(\!\frac{M\t}{N}\right)^2\!G_2(\frac{M\t}{N})\!-\!2\pi i\frac{M\t}{N}\right)\!\frac{(-k\!+\!t\t)}{M\tau}\!\right.\nonumber\\
& & \ \ \ \ \ \ \left.-\!\frac{M\t}{N}\wp_1(\frac{z\!-\!k\!+\!t\t}{N},\frac{M\t}{N})\!\right)+\frac{\epsilon(\la)}{M}\sum_{t=0}^{M-1}\mu^{-t}.\label{g4.18}
\end{eqnarray}

\noindent{\bf Case 1:} $\mu\ne1\ne \la.$ Here (\ref{g4.18}) simplifies to
\begin{eqnarray*}
& &P_1(\mu,\la,\frac{z}{\t},\frac{-1}{\t})
=-\frac{\t}{2\pi iN}\sum_{t=0}^{M-1}\sum_{k=0}^{N-1}\mu^{-t}\la^k\wp_1(\frac{z-k+t\t}{N},\frac{M\t}{N})\\
& &\ \ \ \ \ =-\frac{\t}{2\pi iN}\sum_{t=0}^{M-1}\sum_{k=0}^{N-1}\mu^{-t}\la^{k}\wp_1(\frac{z+N-k+t\t}{N},\frac{M\t}{N})
\end{eqnarray*}
using (\ref{g4.9}). From (\ref{g4.17}) this is indeed equal to $\tau P_1(\la,\mu^{-1},z,\t),$ as required.

\noindent{\bf Case 2: } $\la\ne 1=\mu.$ This time (\ref{g4.18}) reads
\begin{eqnarray}
& &\ \ \ \ P_1(\mu,\la,\frac{z}{\t},\frac{-1}{\t})\nonumber\\
& &=\frac{1}{2\pi iM}\sum_{t=0}^{M-1}\sum_{k=0}^{N-1}\la^k
\left(\left(\left(\frac{M\t}{N}\right)^2G_2(\frac{M\t}{N})-\frac{2\pi iM\t}{N}\right)\left(\frac{-k}{M\t}\right)\right.\nonumber\\
& &\ \ \ \ \ \ \  \ \left.-\frac{M\t}{N}\wp_1(\frac{z-k+t\t}{N},\frac{M\t}{N})\right)+\epsilon(\la).\label{m4.19}
\end{eqnarray}
It is easy to check that $\frac{1}{N}\sum_{k=0}^{N-1}k\l^k=-\epsilon(\l),$
so that (\ref{m4.19}) simplifies to read
\begin{eqnarray*}
& &\ \ \ \ P_1(\mu,\la,\frac{z}{\t},\frac{-1}{\t})\nonumber\\
& &=-\frac{\t}{2\pi iN}\sum_{t=0}^{M-1}\sum_{k=0}^{N-1}\la^k\left(
G_2(\frac{M\t}{N})\frac{k}{N}+\wp_1(\frac{z-k+t\t}{N},\frac{M\t}{N})\right)\\
& &=\frac{\t}{2\pi iN}\sum_{l=0}^{N-1}\sum_{t=0}^{M-1}\la^{-l}
\left(G_2(\frac{M\t}{N})\frac{l}{N}-\wp_1(\frac{z+l+t\t}{N},\frac{M\t}{N})\right)\\
& &=\t P_1(\la,\mu^{-1},z,\t).
\end{eqnarray*}
This completes the discussion of case 2. The final case $\l=1\ne \mu$ is completely analogous, and we accordingly omit details.This completes the proof of 
Theorem \ref{t4.2}, hence also that of Theorem \ref{t4.1}.
\bigbreak

We discuss some aspects of {\em Bernoulli polynomials.} Recall  [La],
[Ra] that
these polynomials $B_k(x)$ are defined by the generating function
\begin{equation}\label{m4.20}
\frac{te^{tx}}{e^t-1}=\sum_{k=0}^{\infty}B_k(x)\frac{t^k}{k!}.
\end{equation}
For example
\begin{equation}\label{m4.21}
B_0(x)=1, B_1(x)=x-\frac{1}{2}, B_2(x)=x^2-x+\frac{1}{6}.
\end{equation}

We will   need the following identities (loc.cit.)
\begin{equation}\label{m4.22}
\sum_{a=0}^{N-1}(a+x)^{k-1}=\frac{1}{k}(B_k(x+N)-B_k(x))
\end{equation}
\begin{equation}\label{m4.22A}
B_k(1-x)=(-1)^kB_k(x).
\end{equation}.
\begin{prop}\label{p4.4} If $\mu=e^{2\pi ij/M}$ with $1\leq j\leq M$ and
$k\geq 2$ then 
$$\frac{1}{(2\pi i)^{k}}\sum_{m\ne 0}{\mu^{m}\over m^k}=\frac{-B_k(j/M)}{k!}.$$
\end{prop}

\noindent{\bf Proof:\ } This is a typical sort of calculation which we give
using   results from [La]. Now
\begin{eqnarray*}
 & &\sum_{0\ne m\in\Z}\mu^{m}/m^k=\sum_{m=1}^{\infty}\left\{
\frac{\mu^{m}}{m^k}+(-1)^k\frac{\mu^{-m}}{m^k}\right\}\\
& &\ \ \ =\sum_{t=1}^M\sum_{n=0}^{\infty}\frac{\mu^t+(-1)^k\mu^{-t}}{(Mn+t)^k}\\
& &\ \ \ =M^{-k}\sum_{t=1}^M\zeta(k,t/M)(\mu^t+(-1)^k\mu^{-t})
\end{eqnarray*}
where $\zeta(k,x)$ is the Hurwitz zeta-function
$$\zeta(k,x)=\sum_{n=0}^{\infty}\frac{1}{(n+x)^k}.$$
For an $M$-th root of unity $\alpha$ define $f_{\alpha}: \Z/M\Z\to\C$ by
$f_\alpha(n)=\alpha^n.$ Set
$$\xi(k,f_\alpha)=M^{-k}\sum_{n=1}^Mf_{\alpha}(n)\zeta(k,n/M).$$
In fact, one can use this definition for any function
$f:\Z/M\Z\to\C.$ Thus
$$\sum_{0\ne m\in\Z}{\mu^{m}\over m^k}=\xi(k,f_{\mu})+(-1)^k\xi(k,f_{\mu^{-1}}).$$
Define 
$$f_j(n)=\left\{
\begin{array}{lr}
1 & n\equiv j\ \mod M\\
0 & n\not\equiv j \ \mod M
\end{array}
\right.$$
and
$$\hat f_j(n)=\sum_{a=1}^Mf_j(a)e^{-2\pi ian/M}.$$
Then we get $\hat f_j(n)=\mu^{-n}.$

Use (loc.cit. Theorem 2.1 P245) to get (remember $k\geq 2$)
\begin{eqnarray*}
& &\xi(1-k,f_j)=(2\pi i)^{-1}\left(\frac{2\pi}{M}\right)^{1-k}\Gamma(k)
\left(\xi(k,f_{\mu})e^{\pi i(1-k)/2}
- -\xi(k,f_{\mu^{-1}})e^{-\pi i(1-k)/2}\right)\\
& &\ \ \ =(2\pi i)^{-k}M^{k-1}(k-1)!\left(
\xi(k,f_{\mu})+(-1)^k\xi(k,f_{\mu^{-1}})\right).
\end{eqnarray*} 
($\Gamma$ is the gamma-function here!) So now we have
$$\sum_{0\ne m\in\Z}\frac{\mu^{m}}{m^k}=\frac{(2\pi i)^kM^{1-k}}{(k-1)!}
\xi(1-k,f_j).$$

On the other hand by definition, 
$$\xi(1-k,f_j)=M^{k-1}\sum_{n=1}^Mf_j(n)\zeta(1-k,n/M)=M^{k-1}\zeta(1-k,j/M)$$
so 
$$\sum_{0\ne m\in\Z}\frac{\mu^{m}}{m^k}=\frac{(2\pi i)^k}{(k-1)!}
\zeta(1-k,j/M).$$
Moreover (loc.cit. Corollary on P243)
$$\zeta(1-k,j/M)=-\Gamma(k)\Res_{z}\frac{z^{-k}e^{zj/M}}{e^z-1}=-(k-1)!B_k(j/M)/k!=-B_k(j/M)/k.$$
So finally
$$\sum_{0\ne m\in\Z}\frac{\mu^{m}}{m^k}=-\frac{(2\pi i)^k}{k!}B_k(j/M).\ \ \ \Box$$
\bigbreak

We next introduce the $Q$-functions. With $\mu=e^{2\pi ij/M}$ and 
$\la=e^{2\pi i l/N},$ we define for $k=1,2,\cdots$ and $(\mu,\la)\ne (1,1),$ 
\begin{eqnarray}
& &Q_{k}(\mu,\la,q_{\t})=Q_k(\mu,\l,\tau)\nonumber\\
& &\ \ \ \ \ \ =\frac{1}{(k-1)!}\sum_{n\geq 0}
\frac{\la(n+j/M)^{k-1}q_{\tau}^{n+j/M}}{1-\lambda q_{\t}^{n-j/M}}\nonumber\\
& &\ \ \ \ \ \ \ \ +\frac{(-1)^k}{(k-1)!}\sum_{n\geq 1}\frac{\la^{-1}(n-j/M)^{k-1}
q^{n-j/M}_{\tau}}{1-\lambda^{-1}q^{n-j/M}_{\tau}}-\frac{B_k(j/M)}{k!}\label{m4.23}.
\end{eqnarray}
Here $(n+j/M)^{k-1}=1$ if $n=0, j=0$ and $k=1.$ Similarly,
$(n-j/M)^{k-1}=1$ if $n=1, j=M$ and $k=1.$  
For good measure we also set 
\begin{equation}\label{m4.24}
Q_0(\mu,\l,\tau)=-1.
\end{equation}

We need to justify the notation, which suggests that $Q_k(\mu,\l,\t)$ depends 
only on $\t$ and the residue classes of $j$ and $l$ modulo $M$ and $N$ respectively. To see this, note that if we provisionally denote by $Q_k'(\mu,\l,\tau)$ the value of (\ref{m4.23}) in which $j$ is replaced by $j+M,$ then we
find that 
$$Q_k'(\mu,\l,\t)-Q_k(\mu,\la,\tau)=\frac{1}{(k-1)!}(j/M)^{k-1}-\frac{B_k(j/M+1)}{k!}+\frac{B_k(j/M)}{k!}=0,$$
the last equality following from (\ref{m4.22}). 

We are going to prove
\begin{th}\label{t4.5} If $k\geq 0$ then $Q_k(\mu,\l,\t)$  is a holomorphic 
modular form of weight $k.$ If $\gamma=\left(\begin{array}{cc} a & b\\ c & d
\end{array}\right)\in \Gamma$ it satisfies 
$$Q_k(\mu,\la,\gamma\tau)=(c\tau+d)^kQ_k((\mu,\la)\gamma,\tau).$$
\end{th}

As usual one needs to deal with the cases $k=1,2$ of Theorem \ref{t4.5} separately. To this end, for each element $\vec{a}=(a_1,a_2)\in \Q^2/\Z^2,$  we recall
the Klein and Hecke forms (loc.cit.) defined as follows:
$$g_{\vec{a}}(\tau)=-q^{B_2(a_1)/2}_{\tau}e^{2\pi ia_2(a_1-1)/2}
(1-q_{a_1\tau+a_2})\prod_{n=1}^{\infty}(1-q^n_{\tau}q_{a_1\tau+a_2})(1-q_{\tau}^nq_{a_1\tau+a_2}^{-1})$$
\begin{equation}\label{m4.24A}
h_{\vec{a}}(\tau)=2\pi i\left\{
a_1-\frac{1}{2}-\frac{q_{a_1\tau+a_2}}{1-q_{a_1\tau+a_2}}
- -\sum_{m=1}^{\infty}\left\{\frac{q_{\tau}^mq_{a_1\tau+a_2}}{1-q_{\tau}^mq_{a_1\tau+a_2}}-\frac{q_{\tau}^mq_{a_1\tau+a_2}^{-1}}{1-q_{\tau}^mq_{a_1\tau+a_2}^{-1}}\right\}\right\}.
\end{equation}. 

Using (\ref{m4.21}) we easily find 
\begin{prop}\label{p4.6} Let $a=(j/M,l/N)\not\in\Z.$ Then 

(i) $h_{\vec{a}}(\tau)=-2\pi iQ_1(\mu,\la,\tau)$

(ii) $\frac{d}{d\tau}(log g_{\vec{a}}(\tau))=-2\pi iQ_2(\mu,\la,\tau).$
\end{prop}

Now we can prove Theorem \ref{t4.5} in the case $k=1,2.$ For $k=1$ we
use (i) above together with Theorem 2 (i) and H3 of (loc.cit. P248). Similarly if $k=2$ we use the calculation of (loc.cit.,P251 et.
seq.)

We now consider the case $k\geq 3.$ In this case the result is a consequence of
the following:
\begin{th}\label{t4.7} If $k\geq 3$ then 
$$Q_k(\mu,\la,\tau)=\frac{1}{(2\pi i)^{k}}\sum^{\ \ \ \ \ \prime}_{m_1,m_2\in\Z}\frac{\la^{-m_1}\mu^{m_2}}{(m_1\tau+m_2)^k}$$
where $\sum'$ indicates that $(m_1,m_2)\ne (0,0).$
\end{th}

\pf The non-constant part of
$\sum'\frac{\la^{-m_1}\mu^{m_2}}{(m_1\tau+m_2)^k}$ is equal to
\begin{eqnarray*}
& &\ \ \  \sum_{m_1=1}^{\infty}\frac{\la^{-m_1}\mu^{m_2}}{(m_1\tau+m_2)^k}
+(-1)^k\sum_{m_1=1}^{\infty}\frac{\la^{m_1}\mu^{m_2}}{(m_1\tau-m_2)^k}\\
& &=\sum_{m_1=1}^{\infty}\sum_{t=1}^M\sum_{n\in\Z}\mu^t\!\left(\frac{\la^{-m_1}}{(m_1\tau+Mn+t)^k}+(-1)^k\frac{\la^{m_1}}{(m_1\tau-Mn-t)^k}\right)\\
& &=M^{-k}\sum_{m_1=1}^{\infty}\sum_{t=1}^M\sum_{n\in\Z}\mu^t\left(\frac{\la^{-m_1}}{(m_1\tau/M+n+t/M)^k}+(-1)^k\frac{\la^{m_1}}{(m_1\tau/M-n-t/M)^k}\right).
\end{eqnarray*}
Use (loc.cit. P155) to get this equal to
\begin{eqnarray*}
& &M^{-k}\!\sum_{m_1=1}^{\infty}\!\sum_{t=1}^M\!\frac{(-1)^k\mu^t(2\pi i)^k}{(k-1)!}\!
\sum_{n=1}^{\infty}\!\left(n^{k-1}\la^{-m_1}q^n_{m_1\tau/M+n+t/M}\!+\!(-1)^kn^{k-1}\la^{m_1}q^n_{m_1\tau/M-t/M}\right)\\
& &=\!(\!-\!1)^kM^{-k}\!\frac{(2\pi i)^k}{(k-1)!}\!
\sum_{m_1=1}^{\infty}\!\sum_{t=1}^M\!\mu^t\sum_{n=1}^{\infty}\!
\left(n^{k-1}\la^{-m_1}\nu^{tn}q^{m_1n}_{\tau/M}+\!(\!-1)^kn^{k-1}\la^{m_1}\nu^{-tn}q^{m_1n}_{\tau/M}\right)\\
& &=(-1)^kM^{-k}\frac{(2\pi i)^k}{(k-1)!}\sum_{t=1}^M\mu^t
\sum_{n=1}^{\infty}
\left(\sum_{d|n}d^{k-1}(\la^{-n/d}\nu^{td}+(-1)^k\la^{n/d}\nu^{-td}\right)q^n_{\tau/M}
\end{eqnarray*}
(where $\nu=e^{2\pi i/M}$). Using orthogonality relations for roots of unity, this is equal to
\begin{eqnarray*}
& &\frac{M^{1-k}(-2\pi i)^k}{(k-1)!}\sum_{n=1}^{\infty}
\left(\sum_{\stackrel{d|n}{d\equiv -j\,(\mod M)}}d^{k-1}\la^{-n/d}+(-1)^k
\sum_{\stackrel{d|n}{d\equiv j\,(\mod M)}}d^{k-1}\la^{n/d}\right)q^n_{\tau/M}\\
& &=\frac{M^{1-k}(2\pi i)^k}{(k-1)!}\sum_{n=0}^{\infty}\sum_{m=1}^{\infty}
\left(\la^{m}(Mn+j)^{k-1}q_{\tau/M}^{m(Mn+j)}\right.\\
& &\ \ \ \ \ \left.+(-1)^k\la^{-m}(Mn+M-j)^{k-1}q_{\tau/M}^{Mn+M-j)}\right)\\
& &=\frac{M^{1-k}(2\pi i)^k}{(k-1)!}\sum_{n=0}^{\infty}
\left(\frac{\la q_{\tau/M}^{Mn+j}}{1-\la q_{\tau/M}^{Mn+j}}
+(-1)^k\frac{\la^{-1}q_{\tau/M}^{Mn+M-j}}{1-\la^{-1}q_{\tau/M}^{Mn+M-j}}\right),
\end{eqnarray*}
where $1\leq j\leq M.$ 
Together with Proposition \ref{p4.4} we now get
\begin{eqnarray*}
& &\ \ \ \ \frac{1}{(2\pi i)^{k}}\sum^{\ \ \ \ \ \prime}\frac{\la^{-m_1}\mu^{m_2}}{(m_1\tau+m_2)^k}
=-\frac{B_k(j/M)}{k!}\\
& &+\frac{M^{1-k}}{(k-1)!}\sum_{n=0}^{\infty}
\left(\frac{\la q_{\tau/M}^{Mn+j}}{1-\la q_{\tau/M}^{Mn+j}}+
(-1)^k\frac{\la^{-1}q_{\tau/M}^{Mn+M-j}}{1-\la^{-1}q_{\tau/M}^{Mn+M-j}}\right)
\\
& &=Q_k(\mu,\la,\tau).
\end{eqnarray*}
This completes the proof of Theorem \ref{t4.7} and hence also Theorem 
\ref{t4.5}. \qed

We remark the well-known fact that if we take $(\mu,\la)=(1,1)$ in Theorem 
\ref{t4.7} we obtain the Eisenstein series of weight $k$ as long as $k\geq 3$
is even. Thus for $k\geq 4$ even,
\begin{eqnarray}\label{m4.25}
& &G_k(\tau)=\sum^{\ \ \ \ \ \prime}_{m_1,m_2\in\Z}\frac{1}{(m_1\tau+m_2)^k}\nonumber\\
& &\ \ \ \ =(2\pi i)^k\left(\frac{-B_k(0)}{k!}+\frac{2}{(k-1)!}\sum_{n=1}^{\infty}\sigma_{k-1}(n)q^n\right)
\end{eqnarray}
where $\sigma_{k-1}(n)=\sum_{d|n}d^{k-1}.$ For this range of values
of $k$, $G_k(\t)$ is a modular form on $SL(2,\Z)$ of weight $k.$

We make use of the normalized Eisenstein series 
\begin{equation}\label{m4.26}
E_k(\tau)=\frac{1}{(2\pi i)^k}G_k(\tau)=\frac{-B_k(0)}{k!}+\frac{2}{(k-1)!}\sum_{n=1}^{\infty}\sigma_{k-1}(n)q^n
\end{equation}
for even $k\geq 2.$ (Warning: this is not the same notation as used in (loc.cit.),
for example.) 

We also utilize the differential operator $\partial_k$ which acts on holomorphic functions on $\frak h$ via
\begin{equation}\label{m4.27}
\p_kf(\t)=\frac{1}{2\pi i}\frac{df(\tau)}{d\t}+kE_2(\t)f(\t).
\end{equation}
One checks using (\ref{g4.7}) that
\begin{equation}\label{m4.28}
(\p_kf)|_{k+2}\gamma=\p_k(f|_k\gamma).
\end{equation}

We complete this section with a discussion of some further functions related
to $P$ and $Q$ which occur later. Again we take $\mu=e^{2\pi ij/M}$ and
$\l=e^{2\pi il/N}.$ Set for $k\geq 1,$
\begin{equation}\label{m4.29}
\bar P_{k}(\mu,\la,z,q_{\t})=\bar P_{k}(\mu,\la,z,\tau)=\frac{1}{(k-1)!}\sum_{n\in j/M+\Z}^{\ \ \ \ \ \prime}\frac{n^{k-1}z^n}{1-\la q^n_{\tau}}
\end{equation}
with
\begin{equation}\label{m4.30}
\bar P_0=0.
\end{equation}

Recall (\ref{ad1}). We shall need the following result in Section 8.
\begin{prop}\label{p4.8}
If $m\in\Z,$ $\mu=e^{2\pi ij/M}, k\geq 0$ and $(\mu,\la)\ne(1,1),$ then
\begin{eqnarray*} 
& &\ \ Q_{k}(\mu,\la,\t)+\frac{1}{k}B_k(1-m+j/M)\\
& &=\Res_{z}\left(\iota_{z,z_1}((z-z_1)^{-1})z_1^{m-j/M}z^{-m+j/M}\bar P_{k}(\mu,\la,\frac{z_1}{z},\t)\right)\nonumber\\
& &\ \ -\Res_{z}\left(\la\iota_{z_1,z}((z-z_1)^{-1})z_1^{m-j/M}z^{-m+j/M}\bar P_{k}(\mu,\la,\frac{z_1q_{\tau}}{z},\tau)\right).
\end{eqnarray*}
\end{prop}

\pf The result is clear if $k=0,$ so take $k\geq 1.$ 
The first of the two residues we must evaluate is equal to
\begin{eqnarray*}
& &\Res_{z}\frac{1}{(k-1)!z}\sum_{n=0}^{\infty}\sum_{r\in\Z}
\left(\frac{(r+j/M)^{k-1}\left(\frac{z_1}{z}\right)^{r+j/M}}{1-\lambda q^{r+j/M}_{\t}}\right)\\
& &=\frac{1}{(k-1)!}\sum_{n=0}^{\infty}\frac{(-n-m+j/M)^{k-1}}{1-\la q_{\t}^{-n-m+j/M}}\\
& &=\frac{(-1)^k}{(k-1)!}\sum_{n=m}^{\infty}\frac{\la^{-1}(n-j/M)^{k-1}q_{\t}^{n-j/M}}{1-\la^{-1} q_{\t}^{n-j/M}}.
\end{eqnarray*}

Similarly the second residue is equal to 
$$ \frac{-1}{(k-1)!}\sum_{n=1-m}^{\infty}\frac{\la(n+j/M)^{k-1}q_{\t}^{n+j/M}}{1-\la q_{\t}^{n+j/M}}.$$
Comparing with the definition of $Q_k(\mu,\la,\tau),$ we see that
we are reduced to establishing that 
\begin{eqnarray}\label{m4.31}
& &\ \ \frac{1}{k!}(B_k(1-m+j/M)-B_k(j/M))\nonumber\\
& &=\left\{
\begin{array}{lr}
\displaystyle{\frac{-1}{(k-1)!}\sum_{n=0}^{-m}\frac{\la(n+j/M)^{k-1}q_{\t}^{n+j/M}}{1-\la q_{\t}^{n+j/M}}} & \\
\ \ \ \displaystyle{+\frac{(-1)^k}{(k-1)!}\sum_{n=m}^{0}\frac{\la^{-1}(n-j/M)^{k-1}q_{\t}^{n-j/M}}{1-\la^{-1} q_{\t}^{n-j/M}}},& m\leq 0\\
\displaystyle{\frac{1}{(k-1)!}\sum_{n=1-m}^{-1}\frac{\la(n+j/M)^{k-1}q_{\t}^{n+j/M}}{1-\la q_{\t}^{n+j/M}}} &\\
\ \ \ \displaystyle{-\frac{(-1)^k}{(k-1)!}\sum_{n=1}^{m-1}\frac{\la^{-1}(n-j/M)^{k-1}q_{\t}^{n-j/M}}{1-\la^{-1} q_{\t}^{n-j/M}}},& m\geq 2\\
0,& m=1
\end{array}
\right.
\end{eqnarray}

The case $m=1$ is trivial. Assume next that $m\leq 0.$ Then the two summations
in (\ref{m4.31}) are equal to 
\begin{eqnarray*}
& &\frac{-1}{(k-1)!}\sum_{n=0}^{-m}\frac{\la(n+j/M)^{k-1}q_{\t}^{n+j/M}}{1-\la q_{\t}^{n+j/M}}-\frac{1}{(k-1)!}\sum_{n=m}^{0}\frac{(-n+j/M)^{k-1}}{\l q_{\t}^{-n+j/M}-1}\\
& &=\frac{1}{(k-1)!}\sum_{n=0}^{-m}(n+j/M)^{k-1}.
\end{eqnarray*}
Now the desired result follows from (\ref{m4.22}). Similarly, if $m\geq 2$ the 
two summands in (\ref{m4.31}) sum to
$$\frac{(-1)^k}{(k-1)!}\sum_{n=1}^{m-1}(n-j/M)^{k-1}.$$
 Two applications of (\ref{m4.22}) then yield the desired result. \qed
\section{The space of 1-point functions on the torus}
\setcounter{equation}{0}

The following notation will be in force for some time.

(a) $V$ is a \voa.

(b) $g,h\in \Aut(V)$ have finite order and satisfy $gh=hg.$

(c) $A=\<g,h\>.$ 

(d) $g$ has order $T,$ $h$ has order $T_1$ and $A$ has exponent $N=lcm(T,T_1).$

(e) $\Gamma(T,T_1)$ is the subgroup of matrices $\left(\begin{array}{cc}
a & b\\
c & d
\end{array}\right)$ in $SL(2,\Z)$ satisfying $a\equiv d\equiv 1$ $(\mod N),$
$b\equiv 0$ $(\mod T),$ $c\equiv 0$ $(\mod T_1).$

(f) $M(T,T_1)$ is the ring of holomorphic modular forms on $\Gamma(T,T_1);$
it is naturally graded
$M(T,T_1)=\oplus_{k\geq 0}M_k(T,T_1)$ where $M_k(T,T_1)$ is the space of forms
of weight $k.$ We also set $M(1)=M(1,1).$

(g) $V(T,T_1)=M(T,T_1)\otimes_{\C}V.$

(h) $O(g,h)$ is the $M(T,T_1)$-submodule of $V(T,T_1)$ generated by the following
elements, where $v\in V$ satisfies $gv=\mu^{-1}v, hv=\la^{-1}v:$
\begin{eqnarray}
& & v[0]w, w\in V, (\mu,\la)=(1,1)\label{m5.1}\\
& & v[-2]w+\sum_{k=2}^{\infty}(2k-1)E_{2k}(\t)\otimes v[2k-2]w, (\mu,\la)=(1,1)\label{m5.2}\\
& & v, (\mu,\la)\ne (1,1)\label{m5.3}\\
& & \sum_{k=0}^{\infty} Q_k(\mu,\la,\t)\otimes v[k-1]w, (\mu,\la)\ne (1,1).
\label{m5.4}
\end{eqnarray}
Here, notation for modular forms is as in Section 4. These definitions
are sensible because of the following:

\begin{lem}\label{l5.1} $M(T,T_1)$ is a Noetherian ring which contains each 
$E_{2k}(\t),$ $k\geq 2,$ and each $Q_k(\mu,\la,\t),$ $k\geq 0,$ for
$\mu,\la$ a $T$-th., resp. $T_1$-th. root of unity.
\end{lem}

\pf It is well-known that the ring of holomorphic modular forms 
on any congruence subgroup
of $SL(2,\Z)$ is Noetherian. So the first statement holds.

Each $E_{2k}$ is a modular form on $SL(2,\Z),$ whereas the containment
$Q_k(\mu,\la,\t)\in M_k(T,T_1)$ follows from Theorem \ref{t4.5}. \qed

\begin{lem}\label{l5.2}  Suppose that $V$ satisfies condition $C_2.$
Then $ V(T,T_1)/O(g,h)$ is a finitely-generated $M(T,T_1)$-module. 
\end{lem}
   
\pf Since $C_2(V)$ is a graded subspace of $V$ of finite 
codimension, there is an integer $m$ such that $V_n\subset C_2(V)$
whenever $n>m.$ Let $M$ be the $M(T,T_1)$-submodule of $V(T,T_1)$ generated 
by $W=\oplus_{n\leq m}V_n.$ The lemma will follow from the assertion that
$V(T,T_1)=M+O(g,h).$ This will be established by proving that if $v\in V_{[k]}$
(cf. (\ref{g2.19})) then $v\in M+O(g,h).$

If $k\leq m$ then $v\in W$ by (\ref{g2.20}) and we are done, so we may
take $k>m.$ Since $V_{[k]}\subset C_2(V)+W$ then we may write $v$ in
the form $$ v=w+\sum_{i=1}^pa_i(-2)b_i$$ with $a_i,b_i\in V$
homogeneous in the \voa (V, [\ ]) such that $\wt [a_i]+\wt[b_i]+1=k.$
Clearly, it suffices to show that $a_i(-2)b_i\in M+O(g,h).$ We may
also assume that $ga=\mu^{-1}a_i, ha_i=\la^{-1}a_i$ for
suitable scalars $\mu, \la.$

Suppose first that $(\mu,\la)=(1,1).$ From (\ref{m5.2}) we see
that $O(g,h)$ contains 
$$a_i[-2]b_i+\sum_{l=2}^{\infty}(2l-1)E_{2l}(\t)\otimes
a_i[2l-2]b_i.$$ Since $\wt [a_i[2l-2]b_i]=k-2l$ then the sum
$$\sum_{l=2}^{\infty}(2l-1)E_{2l}(\t)\otimes
a_i[2l-2]b_i$$ lies in $M+O(g,h)$ by the inductive hypothesis, whence so too 
does $a_i[-2]b_i.$

On the other hand, it follows from (\ref{g2.8}) that we have
$$v(n)=v[n]+\sum_{j>n}\alpha_jv[j]$$
for $v\in V,$ $j\in \Z$ and scalars $\a_j.$ In particular we get
$$a_i(-2)b_i= a_i[-2]b_i+\sum_{j>-2}\alpha_ja_i[j]b_i.$$
Having already shown that each of the summands $a_i[j]b_i$
lies in $M+O(g,h),$ $j\geq -2,$ we get $a_i(-2)b_i\in M+O(g,h)$ as desired.

Now suppose that $(\mu,\la)\ne (1,1).$ In this case (\ref{m5.4}) tells
us that $O(g,h)$ contains the element 
$$-a_i[-1]b_i+
\sum_{l=1}^{\infty}Q_l(\mu,\la,\tau)\otimes a_i[l-1]b_i$$ 
(cf. (\ref{m4.24})). More to the point,
$O(g,h)$ also contains the same expression with $a_i$ replaced by $L[-1]a_i.$
Since $(L[-1]a_i)[t]=-ta_i[t-1]$ by (\ref{g2.7}),
we see that $O(g,h)$ contains the element
$$a_i[-2]b_i+ \sum_{l=1}^{\infty}(l-1)Q_l(\mu,\la,\tau)\otimes a_i[l-2]b_i.$$

Now we proceed as before to conclude that $a_i(-2)b_i\in M+O(g,h).$
\qed

There is a natural grading on $V(T,T_1).$ Namely, the subspace of elements
of degree $n$ is defined to be
\begin{equation}\label{m5.5}
V(T,T_1)_n=\oplus_{k+l=n}M_k(T,T_1)\otimes V_{[l]}.
\end{equation}
Observe that $O(g,h)$ is a graded subspace of $V(T,T_1).$

\begin{lem}\label{l5.3} Suppose $V$ satisfies condition $C_2.$ If $v\in V$ 
then there is $m\in\N$ and elements $r_i(\tau)\in M(T,T_1)$, $0\leq i\leq m-1,$ 
such that
\begin{equation}\label{m5.6}
L[-2]^mv+\sum_{i=0}^{m-1}r_i(\tau)\otimes L[-2]^iv\in O(g,h).
\end{equation}
\end{lem}

\pf Let $I$ be the $M(T,T_1)$-submodule of $V(T,T_1)/O(g,h)$ generated by
$\{L[-2]^iv,i\geq 0\}.$ Since $M(T,T_1)$ is 
 a Noetherian ring, Lemma \ref{l5.2} tells us that $I$ is finitely
generated and so some relation of the form (\ref{m5.6}) must hold.
\qed

\bigskip

We now define the {\em space of $(g,h)$ 1-point functions} ${\cal C}_1(g,h)$
to be the $\C$-linear space consisting of functions
$$S: V(T,T_1)\otimes {\frak h}\to \C$$
which satisfy

(C1) $S(v,\tau)$ is holomorphic in $\t$ for $v\in V(T,T_1).$

(C2) $S(v,\t)$ is $M(T,T_1)$-linear in the sense that $S$ is $\C$-linear in
$v$ and satisfies
\begin{equation}\label{m5.8}
S(f\otimes v,\tau)=f(\tau)S(v,\tau)
\end{equation}
for $f\in M(T,T_1)$ and $v\in V.$

(C3) $S(v,\tau)=0$ if $v\in O(g,h).$

(C4) If $v\in V$ satisfies $gv=hv=v$ then
\begin{equation}\label{m5.10}
S(L[-2]v,\tau)=\partial S(v,\tau)+\sum_{l=2}^{\infty}E_{2l}(\tau)S(L[2l-2]v,\t).
\end{equation}
In (\ref{m5.10}), $\partial S$ is the operator which is linear in $v$ and
satisfies
\begin{equation}\label{m5.11}
\partial S(v,\t)=\partial_kS(v,\t)=\frac{1}{2\pi i}\frac{d}{d\t}S(v,\tau)
+kE_2(\tau)S(v,\t)
\end{equation}
for $v\in V_{[k]}$ (cf. (\ref{m4.27})).

\begin{th}\label{t5.4} (Modular-Invariance) For $S\in {\cal C}_1(g,h)$ and
$\gamma=\left(\begin{array}{cc}
a & b\\
c & d
\end{array}\right)\in SL(2,\Z)$ define 
\begin{equation}\label{m5.12}
S|\gamma(v,\t)=S|_k\gamma(v,\t)=(c\t+d)^{-k}S(v,\gamma\t)
\end{equation}
for $v\in V_{[k]},$ and extend linearly. Then $S|\gamma\in
{\cal C}_1((g,h)\gamma).$
\end{th}

\pf We need to verify that $S|\gamma$ satisfies (C3)-(C4) with 
$(g,h)\gamma=(g^ah^c,g^bh^d)$ in place of $(g,h).$

Step 1: $S|\gamma$ vanishes on $O((g,h)\gamma).$ Pick $v,w\in V$ homogeneous
in $(V,Y[\ ])$ and with $g^ah^{c}v=\mu^{-1}v,$ $g^{b}h^dv=\la^{-1}v.$ Suppose
to begin with that $(\mu,\la)=(1,1).$ We must show that $S|\gamma(u,\tau)=0$
when $u$ is one of the elements (\ref{m5.1}) and (\ref{m5.2}). This follows
easily from (\ref{m5.12}), the equality $S(u,\tau)=0,$ and the fact that
$E_{2k}$ is modular of weight $2k.$

Now assume that $(\mu,\la)\ne (1,1).$ If $gv=\alpha^{-1}v$ and $hv=\beta^{-1}v$
then $(\a,\b)\ne (1,1),$ so that certainly $S|\gamma(v,\t)=(c\t+d)^{-\wt[v]}
S(v,\gamma\t)=0.$ So it remains to show that $S|\gamma(u,\tau)=0$ for
$$u=\sum_{k=0}^{\infty}Q_k(\mu,\la,\tau)\otimes v[k-1]w.$$
First note that we have $(\a,\b)=(\mu,\la)\gamma^{-1}.$ Then
with Lemma \ref{l5.1} and
 Theorem \ref{t4.5} we calculate that
\begin{eqnarray*}
& &\ \ S|\gamma(u,\t)=\sum_{k=0}^{\infty}Q_k(\mu,\l,\t)S|\gamma(v[k-1]w,\t)\\
& &=\sum_{k=0}^{\infty}Q_k(\mu,\l,\t)(c\tau+d)^{-\wt[v]-\wt[w]+k}S(v[k-1]w,\gamma\t)\\
& &=(c\tau+d)^{-\wt[v]-\wt[w]}\sum_{k=0}^{\infty}Q_k(\a,\b,\gamma\t)S(v[k-1]w,\gamma\t)\\
& &=(c\tau+d)^{-\wt[v]-\wt[w]}\sum_{k=0}^{\infty}S(Q_k(\a,\b,\gamma\t)\otimes v[k-1]w,\gamma\t)
\end{eqnarray*}
which is indeed 0 since $S\in {\cal C}_1(g,h).$

Step 2: $S|\gamma$ satisfies (\ref{m5.10}). First note that
if $g^ah^cv=g^bh^dv=v$ then
also $gv=hv=v.$ Then we calculate using (\ref{m4.28}) that
\begin{eqnarray*}
& &S|\gamma(L[-2]v,\t)=(c\tau+d)^{-\wt[v]-2}S(L[-2]v,\gamma\t)\\
& &=(c\tau+d)^{-\wt[v]-2}(\partial S(v,\gamma\tau)+\sum_{k=2}^{\infty}E_{2k}(\gamma\t)S(L[2k-2]v,\gamma\t))\\
& &=(\partial_{\wt[v]}S)|_{\wt[v]+2}\gamma(v,\tau)+
\sum_{k=2}^{\infty}(c\tau+d)^{2k-\wt[v]-2}E_{2k}(\t)S(L[2k-2]v,\gamma\t)\\
& &=\partial_{\wt[v]}(S|_{\wt[v]}\gamma)(v,\tau)+
\sum_{k=2}^{\infty}E_{2k}(\t)S|\gamma(L[2k-2]v,\t).
\end{eqnarray*}
This completes the proof of Step 2, and with it that of the theorem. \qed
\section{The differential equations}
\setcounter{equation}{0}

In this section we study certain differential equations which are satisfied by functions $S(v,\tau)$ in the space of $(g,h)$ 1-point functions. The idea is to exploit Lemma \ref{l5.3} together with (\ref{m5.10}).

We fix an element $S\in{\cal C}_1(g,h).$ 
\begin{lem}\label{l6.1} Let $v\in V$ and suppose that $V$ satisfies condition
$C_2.$ There are $m\in\N$ and $r_i(\tau)\in M(T,T_1)$, $0\leq i\leq m-1,$ such
that 
\begin{equation}\label{m6.1}
S(L[-2]^mv,\t)+\sum_{i=0}^{m-1}r_i(\t)S(L[-2]^iv,\tau)=0.
\end{equation}
\end{lem}

\pf Combine Lemma \ref{l5.3} together with (C2) and (C3). \qed
\bigbreak
In the following we extend (\ref{m5.11}) by setting for $f\in M_l(T,T_1),$
$v\in V_{[k]},$
\begin{equation}\label{6.2}
\partial S(f\otimes v,\t)=\partial_{k+l}(S(f\otimes v,\t))=\partial_{k+l}
(f(\tau)S(v,\tau))
\end{equation}
(cf. (\ref{m4.27}). Then define
\begin{equation}\label{6.3}
\partial^iS(f\otimes v,\t)=\partial_{k+l+2(i-1)}(\partial^{i-1}S(f\otimes v,\tau))
\end{equation}
for $i\geq 1.$ Note that 
\begin{equation}\label{6.4}
\partial S(f\otimes v,\t)=(\partial_{l}f(\t))S(v,\tau)+f(\t)\partial S(v,\t).
\end{equation}
Moreover $\partial_lf(\t)\in M_{l+2}(T,T_1)$ as we see from (\ref{m4.28}).

The simplest case to study is that corresponding to a {\em primary field,}
i.e. a vector $v$ which is a highest weight vector for the Virasoro algebra
in $(V,Y[\ ]).$ Thus $v$ satisfies $L[n]v=0$ for $n>0.$ We assume that this
holds until further notice.

First note that we have
\begin{equation}\label{6.5}
S(L[-2]v,\t)=\partial S(v,\t).
\end{equation}
This follows from (\ref{m5.10}) if $gv=hv=v.$ In general, it is a consequence of this special case, the linearity of $S(v,\t)$ in $v,$ and the identity
$S(w,\tau)=0$ if $gw=\mu^{-1}w,$ $hw=\la^{-1}w$ and $(\mu,\la)\ne (1,1).$
This latter equality follows from (\ref{m5.3}) and (C3). In the
same way, we find from (\ref{m5.10}) that
\begin{equation}\label{6.6}
S(L[-2]^{i+1}v,\t)=\partial S(L[-2]^iv,\tau)
+\sum_{k=2}^{\infty}E_{2k}(\t)S(L[2k-2]L[-2]^iv,\t).
\end{equation}
Using the Virasoro algebra relation we easily find that for $i\in\N$ and 
$k\geq 2$ there are scalars $\a_{ijk},0\leq j\leq i-1$ such that
\begin{equation}\label{6.7}
L[2k-2]L[-2]^iv=\sum_{j=0}^{i-1}\alpha_{ijk}L[-2]^jv,
\end{equation}
so that (\ref{6.6}) becomes
\begin{equation}\label{6.8}
S(L[-2]^{i+1}v,\t)=\partial S(L[-2]^iv,\tau)
+\sum_{j=0}^{i-1}\sum_{k=2}^{\infty}\alpha_{ijk}E_{2k}(\t)S(L[-2]^jv,\t).
\end{equation}

Now proceeding by induction on $i,$ the case $i=1$ being (\ref{6.5}), one proves
\begin{lem}\label{l6.2} Suppose that $L[n]v=0$ for $n>0.$ Then for $i\geq 1$ there are elements $f_j(\t)\in M(1), 0\leq j\leq i-1,$ such that 
\begin{equation}\label{6.9}
S(L[-2]^{i}v,\t)=\partial^iS(v,\tau)
+\sum_{j=0}^{i-1}f_j(\t)\partial^jS(v,\t).
\end{equation}
\end{lem}

Combine Lemmas \ref{l6.2} and \ref{l6.1} to obtain
\begin{lem}\label{l6.3}
Suppose that $V$ satisfies condition $C_2,$ and that $v\in V$ satisfies
$L[n]v=0$ for $n>0.$ Then there are $m\in\N$ and $g_i(\t)\in M(T,T_1),$
$0\leq i\leq m-1,$ such that 
\begin{equation}\label{6.10}
\partial^m S(v,\tau)+\sum_{i=0}^{m-1}g_i(\t)\partial^iS(v,\t)=0.
\end{equation} 
\end{lem}

Bearing in mind the definition of $\partial$ (cf. (\ref{m5.11}),
(\ref{6.3})), (\ref{6.10}) may be reformulated as follows:
\begin{prop}\label{p6.4} Let $R=R(T,T_1)$ be the ring
of holomorphic functions generated by $E_2(\t)$ and $M(T,T_1).$ Suppose that
 $V$ satisfies condition $C_2,$ and that $v\in V$ satisfies
$L[n]v=0$ for $n>0.$ 
Then there are $m\in\N$ and $r_i(\t)\in R(T,T_1),$
$0\leq i\leq m-1,$ such that 
\begin{equation}\label{6.11}
(q_{1\over T}\frac{d}{dq_{1\over T} })^mS(v,\tau)+\sum_{i=0}^{m-1}r_i(\t)(q_{1\over T} \frac{d}{dq_{1\over T} })^iS(v,\t)=0
\end{equation} 
where $q_{1\over T} =e^{2\pi i\tau/T}.$
\end{prop}

We observe here only that $\displaystyle{q\frac{d}{dq}=\frac{1}{T}q_{1\over T} \frac{d}{dq_{1\over T} }
=\frac{1}{2\pi i}\frac{d}{d\tau}.}$

Now (\ref{6.11}) is a homogeneous linear differential equation with
holomorphic coefficients $r_i(\t)\in R,$ and such that $0$ is a
regular singular point. The forms in $R(T,T_1)$ have Fourier expansions
at $\infty$ which are power series in $q_{1\over T} $ because they are
invariant under $\left(\begin{array}{cc} 1 & T\\ 0 &
1\end{array}\right).$ We are therefore in a position to apply the
theory of Frobenius-Fuchs concerning the nature of the solutions to
such equations. A good reference for the elementary aspects of this theory
is [I], but the reader may also consult [AM] where they arise in a context 
related to that of the present paper. We will also need some results
{}from (loc. cit.) in Section 11.
Frobenius-Fuchs theory tells us that $S(v,\tau)$ may be expressed
in the following form: for some $p\geq 0,$
\begin{equation}\label{6.12}
S(v,\tau)=\sum_{i=0}^p(\log q_{1\over T} )^iS_i(v,\tau)
\end{equation}
where
\begin{equation}\label{6.13}
S_i(v,\t)=\sum_{j=1}^{b(i)} q^{\lambda_{ij}}S_{i,j}(v,\tau)
\end{equation}
\begin{equation}\label{6.14}
S_{i,j}(v,\tau)=\sum_{n=0}^{\infty}a_{i,j,n}(v)q^{n/T}
\end{equation}
are holomorphic on the upper half-plane, and
\begin{equation}\label{6.15}
\la_{i,j_1}\not\equiv \la_{i,j_2}\ (\mod \frac{1}{T}\Z)
\end{equation}
for $j_1\ne j_2.$

We are going to prove
\begin{th}\label{t6.5}
Suppose that $V$ satisfies condition $C_2.$ For {\em every} $v\in V,$ 
the function 
$S(v,\tau)\in {\cal C}_1(g,h)$ can be expressed in the form 
(\ref{6.12})-(\ref{6.15}). Moreover, $p$ is bounded independently of $v.$
\end{th}

We begin by proving by induction on $k$ that if $v\in V_{[k]}$ then $S(v,\tau)$
has an expression of the type
(\ref{6.12}).
We have already shown this if $v$ is a highest weight vector for the
Virasoro algebra and in particular if $v$ is in the top 
level $V_{[t]}$ of $(V, Y[\ ]),$ i.e., if $V_{[t]}\ne 0$ and
$V_{[s]}=0$ for $s<t.$ This begins the induction. The proof is an elaboration
of the previous case. We may assume that $gv=hv=v.$
\begin{lem}\label{l6.6}
Suppose that $l\geq 2$ and $i\geq 0.$ Then there are scalars 
$\a_{ijl}$ and $w_{ijl}\in V_{[2i+2-2l-2j+k]},$ $0\leq j\leq i-1,$ such that
\begin{equation}\label{6.16}
L[2l-2]L[-2]^{i}v=
L[-2]^{i}L[2l-2]v+\sum_{j=0}^{i-1}\alpha_{ijl}L[-2]^jw_{ijl}.
\end{equation}
Moreover $\wt [w_{ijl}]\leq \wt [v]$ with equality only if $w_{ijl}=v.$
\end{lem}

\pf By induction on $i+l,$ the case $i=0$ being trivial. Now we calculate
\begin{eqnarray*}
& &\ \ \ \ \ L[2l-2]L[-2]^{i+1}v=(L[-2]L[2l-2]+2lL[2l-4]+\delta_{l,2}\frac{c}{2})L[-2]^iv\\
& &=L[-2]^{i+1}L[2l-2]v+\sum_{j=0}^{i-1}\alpha_{ijl}L[-2]^{j+1}w_{ijl}+
  2lL[2l-4]L[-2]^iv+\delta_{l,2}\frac{c}{2}L[-2]^iv.
\end{eqnarray*}
Either $l=2$ or the inductive hypothesis applies to $L[2l-4]L[-2]^iv,$
and in either case the lemma follows. \qed

Now use (\ref{m5.10}) and Lemma \ref{l6.6} to see that
\begin{eqnarray}
& &S(L[-2]^{i+1}v,\tau)=\partial S(L[-2]^iv,\tau)\nonumber\\
& &+\sum_{l=2}^{\infty} E_{2l}(\t)(S(L[-2]^iL[2l-2]v,\tau)+\sum_{j=0}^{i}
\alpha_{ijl}S(L[-2]^jw_{ijl},\t)).\label{6.17}
\end{eqnarray}
Note that (\ref{6.17}) is the appropriate generalization of (\ref{6.8}).
By induction based on (\ref{6.17}) we find
\begin{lem}\label{l6.7}
For $i\geq 1$ we have 
\begin{equation}\label{6.18}
S(L[-2]^{i}v,\tau)=\partial^i S(v,\tau)+\sum_{j=0}^{i-1} f_{ij}(\t)\partial^jS(v,\t)+\sum_{j=0}^{i-1}\sum_lg_{ijl}(\t)\partial^jS(w_{ijl},\tau),
\end{equation}
where $f_{ij}(\t), g_{ijl}(\t)\in M(1)$ and $\wt [w_{ijl}]<\wt [v].$ \qed
\end{lem}

The analogue of Lemma \ref{l6.3} is now
\begin{lem}\label{l6.8} There is $m\in\N$ such that 
\begin{equation}\label{6.19}
\partial^m S(v,\tau)+\sum_{i=0}^{m-1}g_i(\t)\partial^iS(v,\t)+
\sum_{j=0}^{m}\sum_lh_{jl}(\t)\partial^jS(w_{jl},\t)=0
\end{equation} 
for $g_i(\tau), h_{jl}(\t)\in M(T,T_1),$ and $\wt [w_{jl}]<\wt [v].$ \qed
\end{lem}

We are now in a position to complete the proof that $S(v,\tau)$ has an
expression  of the form (\ref{6.12})-(\ref{6.15}). By induction this is true 
of the terms $S(w_{jl},\t)$ in (\ref{6.19}), and hence the third 
summand on the r.h.s of (\ref{6.19}) also has such an expression. Thus as 
before we may view (\ref{6.19}) as a differential equation of regular
singular type, this time inhomogeneous, namely,
\begin{equation}\label{6.20}
(q_{1\over T} \frac{d}{dq_{1\over T} })^mS(v,\tau)+\sum_{i=0}^{m-1}r_i(\t)(q_{1\over T} \frac{d}{dq_{1\over T} })^iS(v,\t)+\sum_{i=0}^p(\log q_{1\over T} )^iu_i(v,\tau)=0
\end{equation} 
where $r_i(\tau)\in R(T,T_1)$ (cf. Proposition \ref{p6.4}) and $u_i(v,\tau)$
satisfies (\ref{6.13})-(\ref{6.15}).

One easily sees (cf. [I], [AM]) that the functions $(\log q_{1\over T} )^iu_i(v,\tau),$ $0\leq i\leq p,$ are themselves solutions of a differential 
equation of regular singular type (\ref{6.11}) with coefficients
analytic in the upper half plane. Let us formally state this by saying that
they are solutions of the differential equation $L_1f=0$ where $L_1$ is a
suitable linear differential operator with $0$ as regular singular point and
coefficients analytic in the upper half plane. Now
(\ref{6.20})  takes the form $L_2S+f=0$ for the corresponding linear
differential operator $L_2,$ so that we get $L_1L_2S=0.$ But $L_1L_2$ is
once again a linear differential operator of the appropriate type, so again
the Frobenius-Fuchs theory allows us to conclude that $S(v,\tau)$ indeed
satisfies (\ref{6.12})-(\ref{6.15}).    

It remains to prove that the integer $p$ in (\ref{6.12}) can be
bounded independently of $v.$ Indeed, we showed in Lemma \ref{l5.2}
that if $W=\oplus_{n\leq m}V_n$ and $V_n\subset C_2(V)$ for $n>m$ then
$V(T,T_1)/O(g,h)$ is generated as $M(T,T_1)$-module by $W.$ So for
$v\in V$ we have $v\equiv \sum_{i}f_i(\tau)\otimes w_i$
($\mod\,O(g,h)$) where $\{w_i\}$ is a basis for $W,$ whence
$S(v,\tau)=\sum_{i}f_i(\tau)S(w_i,\tau)$ since $S$ vanishes on $O(g,h).$
Clearly then, we may take $p$ to be the maximum of the corresponding
integers determined by $S(w_i,\tau).$
This completes the proof of Theorem \ref{t6.5}.

\section{Formal 1-point functions}
\setcounter{equation}{0}

Although we dealt with holomorphic functions in Section 6, the
arguments were all formal in nature. In this short section we record 
a consequence of this observation.

We identify elements of $M(T,T_1)$ with
their Fourier expansions at $\infty,$ which lie in the ring
of formal power series $\C[[q_{1\over T} ]].$ Similarly, the functions $E_{2k}(\t),$
$k\geq 1,$ are considered to lie in $\C[[q]].$ 
The operator $\p$ (cf. (\ref{m5.11}), (\ref{6.2})) operates on these
and other power series via the identification $\frac{1}{2\pi i}\frac{d}{dq}=
\frac{q_{1\over T} }{T}\frac{d}{dq_{1\over T} }.$

A {\em formal} $(g,h)$ 1-point function is a map 
$$S: V(T,T_1)\otimes {\frak h}\to P$$
where $P$ is the space of formal power series of the form
\begin{equation}\label{7.1}
q^{\l}\sum_{n=0}^{\infty}a_nq^{n/T}
\end{equation}
for some $\l\in\C,$ and which satisfies the formal analogues
of (C2)-(C4) in Section 5. We will establish
\begin{th}\label{t7.1} Suppose that $S$ is a formal $(g,h)$ 1-point
function. Then $S$ defines an element of ${\cal C}_1(g,h),$ also
denoted by $S,$ via the identification
\begin{equation}\label{7.2}
S(v,\tau)=S(v,q), q=q_{\tau}=e^{2\pi i\tau}.
\end{equation}
\end{th}

The main point is to show that if $S$ is a formal $(g,h)$ 1-point function,
and if $v\in V$ is such that 
$$S(v,q)=q^{\lambda}\sum_{n=0}^{\infty}a_nq^{n/T}$$
then $q_{\t}^{\lambda}\sum_{n=0}^{\infty}a_nq_{\tau}^{n/T}$ is holomorphic
in $\t.$
We prove this as in Section 6. Namely, by first showing that if $v$ is
a highest weight vector for the Virasoro algebra then $S(v,q)$ satisfies
a differential equation of the type (\ref{6.11}). Since the coefficients
are holomorphic in $\frak h,$ the Frobenius-Fuchs theory tells us that $S(v,q)$
has the desired convergence.

Proceeding by induction on $\wt [v],$ in the general case we arrive at an 
inhomogeneous differential equation of the type (\ref{6.20}). Again
convergence of $S(v,q)$ follows from the Frobenius-Fuchs theory. Since
the proofs of these assertions are {\em precisely} the same as those of Section
6, we omit further discussion.

\section{Correlation functions}
\setcounter{equation}{0}

In this section we start to relate the theory of 1-point functions to that
of twisted $V$-modules. We keep the notation
(a)-(h) introduced at
the beginning of Section 5, and introduce now a simple $g$-twisted
$V$-module $M=M(g)=\oplus_{n=0}^{\infty}M_{\la+n/T}$ 
(cf. (\ref{3.10'})).
 We further assume that $h$ 
leaves $M$ stable, that is $h\circ M\simeq M.$ As remarked in Section 3, there
is a projective representation on $M$ of the stablizer (in $\Aut V$) of $M,$
and we let $\phi(h)$ be a linearized action on $M$ of the
element corresponding to $h.$ This all means (cf. (\ref{g3.19}), (\ref{g3.20}))
that if $v\in V$ operates on $M$ via the vertex operator $Y_M(v,z)$ then we
have (as operators on $M$) 
\begin{equation}\label{8.1}
 \phi(h)Y_M(v,z)\phi(h)^{-1}=Y_M(hv,z).
\end{equation}

We define $M'=\oplus_{n=0}^{\infty}M_{\l+n/T}'$ to be the
{\em restricted dual} of $M,$ so that $M_n'\!=\!\Hom_{\C}(M_n,\C)$ and there
is a pairing $\<,\>:M'\times M\to \C$ such that $\<M_n',M_m\>=0$ if $m\ne n.$

With this notation, a $(g,h)$ 1-{\em point correlation function}
is essentially a trace function, namely
\begin{equation}\label{8.2}
\tr Y_M(v,z)=\sum_{w}\<w',Y_M(v,z)w\>
\end{equation}
where $w$ ranges over a homogeneous basis of $M$ and $w'$ ranges over the 
dual basis of $M'.$ As power series in $z$ we have 
\begin{equation}\label{8.3}
\tr Y_M(v,z)=\sum_{w}\sum_{n\in\frac{1}{T}\Z}\<w',v(n)w\>z^{-n-1}.
\end{equation}
It is easy to see that the trace function is independent of the choice
of basis.

Now we introduce the function $T$ which is linear in $v\in V,$ and defined for
homogeneous $v\in V$ as follows:
\begin{equation}\label{8.4}
T(v)=T_M(v,(g,h),q)=z^{\wt v}{\tr}Y_M(v,z)\phi(h)q^{L(0)-c/24}.
\end{equation}
Here, $q$ is indeterminate and $c$ is the central charge of $V.$ 

Next observe that for $m\in \frac{1}{T}\Z,$ $v(m)$ maps $M_{n}$ to
$M_{n+\wt v-m-1}.$ So unless $m=\wt v-1,$ we have $\sum\<w',v(m)\phi(h)w\>=0.$
So only the zero mode $o(v)=v(\wt v-1)$ contributes to the sum in
(\ref{8.3}). Thus $T(v)$ is independent of $z,$ and
\begin{equation}\label{8.6}
T(v)=q^{\l-c/24}\sum_{n=0}^{\infty}\tr_{M_{\l+n/T}}o(v)\phi(h)q^{n/T}.
\end{equation}
We could equally write
\begin{equation}\label{8.7}
T_M(v)=\tr_{M}o(v)\phi(h)q^{L(0)-c/24}.
\end{equation}

We are going to prove
\begin{th}\label{t8.1} $T(v)\in {\cal C}_1(g,h).$
\end{th}

The strategy is to prove that $T$ is a formal $(g,h)$ 1-point
function, then invoke Theorem \ref{t7.1}. Certainly $T(v)$ has the correct
shape as a power series in $q$ (cf. (\ref{7.1}). So we must establish
that $T(v)$ satisfies the formal analogues of (C2)-(C4).

We can impose $M(T,T_1)$-linearity (C2) by extension of
scalars. As we shall explain, the proof of (C4) is contained in
Zhu's paper [Z]. So it remains to discuss (C3), that is we must
show that $T(v)$ vanishes on $O(g,h),$ i.e., on the elements of type
(\ref{m5.1})-(\ref{m5.4}). Again we shall later explain that (\ref{m5.1})
and (\ref{m5.2}) may be deduced from results in [Z], so we concentrate 
on (\ref{m5.3}) and (\ref{m5.4}).

To this end, let us now fix a homogeneous $v\in V$ such that $gv=\mu^{-1}v,$
$hv=\l^{-1}v$ and $(\mu,\l)\ne (1,1).$ We need to establish
\begin{lem}\label{l8.2}
$T(v)=0.$
\end{lem}
\begin{th}\label{t8.3} $\sum_{k=0}^{\infty}Q_k(\mu,\la,\tau)T(v[k-1]w)=0$
for any $w\in V.$ 
\end{th}

The proof of Lemma \ref{l8.2} is easy. We have already seen that only
the zero mode $o(v)$ of $v$ contributes a possible non-zero term in the
calculation of $T(v).$ On the other hand, if $\mu\ne 1$ then from
(\ref{g3.12}) we see that $o(v)=0.$ So Lemma \ref{8.2} certainly holds
if $\mu\ne 1.$ 

Suppose that $\la\ne 1.$ We have 
$$T(v)=\tr_Mo(v)\phi(h)q^{L(0)-c/24}=\tr_M\phi(h)o(v)q^{L(0)-c/24},$$
i.e.,
\begin{equation}\label{8.7'}
\tr Y_M(v,z)\phi(h)q^{L(0)}=\tr \phi(h)Y_M(v,z)q^{L(0)}.
\end{equation}
But (\ref{8.1}) yields
\begin{equation}\label{8.8}
\phi(h)Y_M(v,z)=\la^{-1}Y_M(v,z)\phi(h).
\end{equation}
As $\la\ne 1, $ (\ref{8.7'}) and (\ref{8.8}) yield $\tr Y_M(v,z)\phi(h)q^{L(0)}=0.$ This completes the proof of Lemma \ref{l8.2}.

The proof of Theorem \ref{t8.3} is harder. We first need to define
$n$-point correlation functions. These are multi-linear functions $T(v_1,...,v_n),$ $v_i\in V,$ defined for $v_i$ homogeneous via
 \begin{eqnarray}\label{8.9}
& &T(v_1,...,v_n)=T((v_1,z_1),...,(v_n,z_n),(g,h),q)\nonumber\\
& &\ \ =z_1^{\wt v_1}\cdots z_n^{\wt v_n}{\tr}Y_M(v_1,z_1)\cdots Y_M(v_n,z_n)\phi(h)q^{L(0)-c/24}.
\end{eqnarray}

We only need
the case $n=2.$ We will prove 
\begin{th}\label{t8.4} Let $v,v_1\in V$ be homogeneous with $gv=\mu^{-1}v,$
$hv=\la^{-1}v$ and $(\mu,\la)\ne (1,1).$ Then
\begin{equation}\label{8.9a}
T(v,v_1)=\sum_{k=1}^{\infty}\bar P_{k}(\mu,\la,\frac{z_1}{z},q)T(v[k-1]v_1)
\end{equation}
\begin{equation}\label{8.9b}
T(v_1,v)=\la\sum_{k=1}^{\infty}\bar P_{k}(\mu,\la,\frac{z_1}{z}q,q)T(v[k-1]v_1)
\end{equation}
where $\bar P_k$ is as in (\ref{m4.29}).
\end{th}

We start with

\begin{lem}\label{l8.5} Let $k\in\frac{1}{T}\Z.$ Then  
\begin{eqnarray}\label{8.9c}
& &\ \ \  \ (1-\la q^k)\tr v(\wt v-1+k)Y_M(v_1,z_1)\phi(h)q^{L(0)}\nonumber\\
& &=\sum_{i=0}^{\infty}{\wt v-1+k\choose i}z_1^{\wt v-1+k-i}
\tr Y_M(v(i)v_1,z_1)\phi(h)q^{L(0)}
\end{eqnarray}
\begin{eqnarray}\label{8.9d}
& &\ \ \ \ (1-\la q^k)\tr Y_M(v_1,z_1)v(\wt v-1+k)\phi(h)q^{L(0)}\nonumber\\
& &=\la q^k\sum_{i=0}^{\infty}{\wt v-1+k\choose i}z_1^{\wt v-1+k-i}
\tr Y_M(v(i)v_1,z_1)\phi(h)q^{L(0)}.
\end{eqnarray}
\end{lem}

\pf We have 
\begin{eqnarray}\label{8.10}
& &\ \ \ \ \tr v(\wt v-1+k)Y_M(v_1,z_1)\phi(h)q^{L(0)}\nonumber\\
& &=tr [v(\wt v-1+k),Y_M(v_1,z_1)]\phi(h)q^{L_0}\nonumber\\
& &\ \ \ \ \ \ +\tr Y_M(v_1,z_1)v(\wt v-1+k)\phi(h)q^{L(0)}.
\end{eqnarray}
{}From (\ref{8.8}) we get 
$$v(\wt v-1+k)\phi(h)=\la\phi(h)v(\wt v-1+k),$$
 moreover,
$$v(\wt v-1+k)q^{L(0)}=q^kq^{L(0)}v(\wt v-1+k).$$
Hence
\begin{equation}\label{8.12}
\tr Y_M(v_1,z_1)v(\wt v-1+k)\phi(h)q^{L(0)}
=\la q^k\tr v(\wt v-1+k)Y_M(v_1,z_1)\phi(h)q^{L(0)}.
\end{equation}
Using the relation 
$$[v(m),Y_M(v_1,z_1)]=\sum_{i=0}^{\infty}{m\choose i}z_1^{m-i}Y_M(v(i)v_1,z_1)$$
which is a consequence of the Jacobi identity (\ref{g2.3}) we get
\begin{eqnarray}\label{8.11}
\tr[v(\wt v-1+k),Y_M(v_1,z_1)]=\sum_{i=0}^{\infty}{\wt v-1+k\choose i}
z_1^{\wt v-1+k-i}Y_M(v(i)v_1,z_1).
\end{eqnarray}

Both parts of the Lemma follow from (\ref{8.12})-(\ref{8.11}). \qed

Now we turn to the proof of (\ref{8.9a}) of Theorem \ref{t8.4}. Using
(\ref{8.9c}) in the last lemma and setting $\mu=e^{2\pi ir/T}$ we have
\begin{eqnarray*}
& &T(v,v_1)=T((v,z),(v_1,z_1),(g,h),q)\\ 
& &=z^{\wt v}z_1^{\wt v_1}{\tr}Y_M(v,z)Y_M(v_1,z_1)\phi(h)q^{L(0)-c/24}\\
& &=z^{\wt v}z_1^{\wt v_1}\sum_{k\in\Z+\frac{r}{T}}z^{-\wt v-k}{\tr}
v(\wt v-1+k)Y_M(v_1,z_1)\phi(h)q^{L(0)-c/24}\\
& &=z_1^{\wt v_1}\sum_{k\in \Z+\frac{r}{T}}z^{-k}(1-\la q^k)^{-1}\sum_{i=0}^{\infty}{\wt v-1+k\choose i}z_1^{\wt
     v-1+k-i}Y_M(v(i)v_1,z_1)\phi(h)q^{L(0)-c/24} \\
& &=\sum_{k\in\Z+{r\over T}}(\frac{z_1}{z})^k(1-\la q^k)^{-1}
\sum_{i=0}^{\infty}{\wt v-1+k\choose i}T(v(i)v_1)\\
& &=\sum_{k\in\Z+{r\over T}}(\frac{z_1}{z})^k(1-\la q^k)^{-1}
\sum_{i=0}^{\infty}\sum_{m=0}^ic(\wt v,i,m)k^mT(v(i)v_1)\\
& &=\sum_{i=0}^{\infty}\sum_{m=0}^i m!\bar P_{m+1}(\mu,\la,z_1/z,q)c(\wt v,i,m)
T(v(i)v_1)\\
& &=\sum_{m=0}^{\infty}\bar P_{m+1}(\mu,\la,z_1/z,q)T(v[m]v_1),
\end{eqnarray*}
where we have used (\ref{g2.9}) and (\ref{m4.29}).
This is precisely (\ref{8.9a}) of Theorem \ref{t8.4}. (\ref{8.9b}) follows in the same way by using (\ref{8.9d}). \qed

Before proving Theorem \ref{t8.3} we still need
\begin{lem}\label{l8.6} We have
$$ \sum_{k=0}^{\infty}\frac{1}{k!}B_k(1-\wt v+r/T)v[k-1]=\sum_{i=0}^{\infty}
{r/T\choose i}v(i-1).$$
\end{lem}

\pf The l.h.s. of the equality is equal to 
\begin{eqnarray*}
& &\Res_{w}Y[v,w]\frac{e^{(1-\wt v+r/T)w}}{e^w-1}\\
& &=\Res_{w}Y(v,e^w-1)\frac{e^{(1+r/T)w}}{e^w-1}\\
& &=\Res_zY(v,z)\frac{(1+z)^{r/T}}{z}\\
& &=\sum_{i=0}^{\infty}{r/T\choose i}v(i-1)
\end{eqnarray*}
as required. \qed

{\em Proof of Theorem \ref{t8.3}:} Combine Lemma \ref{l8.6} and Proposition
\ref{p4.8} to get 
\begin{eqnarray*}
& &\ \ \sum_{k=0}^{\infty}Q_k(\mu,\la,q)T(v[k-1]w)\\
& &=\sum_{k=1}^{\infty}\Res_{z}\left(\iota_{z,z_1}((z-z_1)^{-1})z_1^{\wt v-r/T}z^{-\wt v+r/T}\bar P_{k}(\mu,\la,\frac{z_1}{z},q)\right)T(v[k-1]w)\\
& &\ \ -\la\sum_{k=1}^{\infty}\Res_{z}\left(\iota_{z_1,z}((z-z_1)^{-1})z_1^{\wt v-r/T}z^{-\wt v+r/M}\bar P_{k}(\mu,\la,\frac{z_1q}{z},q)\right)T(v[k-1]w)\\
& &\ \ -\sum_{i=0}^{\infty}{r/T\choose i}T(v(i-1)w).
\end{eqnarray*}

On the other hand, use (\ref{ad2}) to obtain 
\begin{eqnarray*}
& &\sum_{i=0}^{\infty}{r/T\choose i}T(v(i-1)w)=\sum_{i=0}^{\infty}{r/T\choose i}z_1^{\wt v+\wt w-i}\tr Y_M(v(i-1)w,z_1)\phi(h)q^{L(0)-c/24}\\
& & =\sum_{i=0}^{\infty}{r/T\choose i}z_1^{\wt v+\wt w-i}\Res_{z-z_1}(z-z_1)^{i-1}
\tr Y_M(Y(v,z-z_1)w,z_1)\phi(h)q^{L(0)-c/24}\\
& &=\Res_{z-z_1}\iota_{z_1,z-z_1}\left(\frac{z}{z_1}\right)^{r/T}(z-z_1)^{-1}z_1^{\wt v+\wt w}
\tr Y_M(Y(v,z-z_1)w,z_1)\phi(h)q^{L(0)-c/24}\\
& &=\Res_z\iota_{z,z_1}(z-z_1)^{-1}({z_1\over z})^{\wt v-r/T}T(v,w)
\\
& &\ \ \ \ \ -\Res_z\iota_{z_1,z}(z-z_1)^{-1}({z_1\over z})^{\wt v-r/T}T(w,v)
\end{eqnarray*}
which by Theorem \ref{t8.4} is equal to 
\begin{eqnarray*}
& &\sum_{k=1}^{\infty}\Res_z\iota_{z,z_1}(z-z_1)^{-1}({z_1\over z})^{\wt v-r/T}
\bar P_{k}(\mu,\la,\frac{z_1}{z},q)T(v[k-1]w)\\
& &-\la\sum_{k=1}^{\infty}\Res_z\iota_{z_1,z}(z-z_1)^{-1}({z_1\over z})^{\wt v-r/T}
\bar P_{k}(\mu,\la,\frac{z_1q}{z},q)T(v[k-1]w).
\end{eqnarray*}
This completes the proof of the theorem. \qed

In order to complete the proof of Theorem \ref{t8.1} we need to explain how 
(\ref{m5.10}), and the fact that $T(v)$ vanishes on (\ref{m5.1}) and (\ref{m5.2}),
follow from [Z]. These results concern the case in which the critical
vector $v$ satisfies $gv=hv=v.$ Thus $v$ lies in the invariant
sub \voa $V^A$. Now Zhu's proof of (\ref{m5.10}), for
example, is quite general in the sense that it does not depend on 
any special properties of $V.$ In particular, his argument applies to $V^A,$
which is what we need. (Note that $M$ is a {\em module} for $V^A.$) 
 Similarly, Zhu's argument establishes that
$T$ vanishes on (\ref{m5.1}) and (\ref{m5.2}) in the case that $g$ and $h$ 
both fix $v$ and $w.$ On the other hand we may certainly assume that
$w$ is an eigenvector for $g$ and $h.$ If $gw=\alpha w,$ $hw=\beta w$ and 
$(\alpha,\beta )\ne (1,1),$ we have already see (cf. the proof of 
Lemma \ref{l8.2}) that $T(v[2k-2]w)=0,$ so it is clear
in this case that $T$ vanishes on (\ref{m5.2}). This completes
the proof of Theorem \ref{t8.1}.

\begin{th}\label{t8.7} Let $M^1,M^2,...$ be inequivalent simple $g$-twisted
$V$-modules, each of which is $h$-stable. Let $T_1,T_2,...$ be the
corresponding trace functions (\ref{8.7}). Then $T_1,T_2,...$ are linearly independent
elements of ${\cal C}_1(g,h).$
\end{th}

\pf Suppose false. Then we may choose notation so that for some
$m\in\N$ there are non-zero scalars $c_1,c_2,...,c_m$ such that
\begin{equation}\label{8.13}
c_1T_1+\cdots +c_mT_m=0.
\end{equation}

Let $\Omega_i$ be the top level of $M^i, 1\leq i\leq m,$ and let $\la_i$ be the
conformal weight of $M^i$ (cf. Section 3). Thus $M^i$ is graded by $\frac{1}{T}\Z,$
$$M^i=\oplus_{n=0}^{\infty}M^i_{\la_i+n/T}$$
and $\Omega_i=M^i_{\la_i}.$

Define a partial order $<<$ on the $\l_i$ by declaring that 
\begin{equation}\label{8.14}
\la_i<<\l_j, \ {\rm if\ and \ only\ if }\ \l_j-\l_i\in\frac{1}{T}\Z_+.
\end{equation}
We may, and shall, assume that $\l_1$ is minimal with respect to $<<.$
 
By Theorem \ref{t3.5} the $\Omega_i$ realize inequivalent irreducible
representations of the algebra $A_g(V).$ Moreover by
(\ref{3.25}) the $A_g(V)$-modules $\Omega_i$ are $h$-invariant, and
the corresponding trace functions $\tr|_{\Omega_i}o(v)\phi(h)$ are linearly
independent. Thus we may choose $v$ so that  $\tr|_{\Omega_1}o(v)\phi(h)=1,$
$\tr|_{\Omega_i}o(v)\phi(h)=0$ for $i>1.$ 

Because of our assumption on $\lambda_1,$ applying this to (\ref{8.13}) yields $c_1=0,$ contradiction. \qed

\section{An existence theorem for $g$-twisted modules}
\setcounter{equation}{0}

We will prove 
\begin{th}\label{t9.1} Suppose that $V$ is a simple \voa which satisfies 
condition $C_2,$ and that $g\in \Aut V$ has finite order. Then $V$ has 
at least one simple $g$-twisted module.
\end{th}

The idea is to prove that $A_g(V)\ne 0.$ Then the theorem is a
consequence of Proposition \ref{p3.8}.  We start with more general 
considerations that we shall need in Section 10. Let $(g,h)$ be a pair of
commuting elements in $\Aut (V).$ 

\begin{lem}\label{l9.3} Let $v\in V$ satisfy $gv=\mu^{-1}v, hv=\la^{-1}v.$
Then the following hold:

(a) The constant term of $\sum_{k=0}^{\infty}Q_{k}(\mu,\la,q)v[k-1]w$ is equal
to $-v\circ_g w$ if $\mu\ne 1.$

(b) The constant term of $\sum_{k=0}^{\infty}Q_{k}(\mu,\la,q)(L[-1]v)[k-1]w$ is equal
to $-v\circ_g w$ if $\mu=1,$ $\la\ne 1.$

(c) The constant term of $v[-2]w+\sum_{k=2}^{\infty}(2k-1)E_{2k}(q)v[2k-2]w$ is 
$v\circ_g w$ if $\mu=\la=1.$
\end{lem}

\pf As usual, (c) follows from the corresponding result in [Z]. The proof
of (a) is similar to that of Lemma \ref{l8.6}. For from (\ref{m4.23})
we see that if we take $\mu=e^{2\pi ir/T}$ with $1\leq r\leq T$ then 
the constant term of the expression in (a) is equal to the following 
(take $v$ homogeneous):
\begin{eqnarray*}
& &\ \ -\sum_{k=0}^{\infty}\frac{B_k(r/T)}{k!}v[k-1]w\\
& &=-\Res_{z}Y[v,z]w\frac{e^{rz/T}}{e^z-1}\\
& &=-\Res_{z}Y(v,e^z-1)w\frac{e^{(\wt v+r/T)z}}{e^z-1}\\
& &=-\Res_zY(v,z)w\frac{(1+z)^{\wt v-1+r/T}}{z}
\end{eqnarray*}
and by (\ref{g3.15}) this is exactly $-v\circ_g w$ if $r\ne T.$

As for (b), we replace $v$ by $L[-1]v =L(-1)v+L(0)v$ (cf. (\ref{g2.17}))
and set 
$r=T$ in the foregoing. From (\ref{m4.23}) the constant term of 
$$\sum_{k=0}^{\infty}Q_k(1,\la,q)v[k-1]w$$ is
\begin{eqnarray*}
& &\ \ -\sum_{k=0}^{\infty}\frac{B_k(1)}{k!}v[k-1]w+\frac{1}{1-\l}v[0]w\\
& &=-\Res_zY(v,z)w\frac{(1+z)^{\wt v}}{z}+\frac{1}{1-\l}v[0]w.
\end{eqnarray*}
Note that $(L[-1]v)[0]w=0.$ 
Then the constant term of the expression in
(b) is equal to  
\begin{eqnarray*}
& & \ \ -\Res_zY(L(-1)v,z)w\frac{(1+z)^{\wt v+1}}{z}-\wt v\Res_zY(v,z)w\frac{(1+z)^{\wt v}}{z}\\
& &=-\Res_z\left(\frac{d}{dz}Y(v,z)\right)w\frac{(1+z)^{\wt v+1}}{z}
 -\wt v\Res_zY(v,z)w\frac{(1+z)^{\wt v}}{z}\\
& &=\Res_zY(v,z)w\frac{d}{dz}\frac{(1+z)^{\wt v+1}}{z}
-\wt v\Res_zY(v,z)w\frac{(1+z)^{\wt v}}{z}\\
& &=-\Res_zY(v,z)w\frac{(1+z)^{\wt v}}{z^2}\\
& &=-v\circ_g w.
\end{eqnarray*}
This completes the proof of the lemma. \qed

Now take $S\in{\cal C}_1(g,h)$ and assume that $S\ne 0.$ After Theorem
\ref{t6.5} we may choose $p$ so that (\ref{6.12})-(\ref{6.15}) hold for all
$v\in V,$ and such that $S_p\ne 0.$  We may further choose notation such that
$\l_{p,1}$ is minimal among all $\l_{p,j}$ with respect to
the partial order (\ref{8.14}) and $a_{p,1,0}(v)\ne 0$ for some
$v\in V.$ Setting 
\begin{equation}\label{9.1}
S_{p,1}(v,\tau)=\a(v)+\sum_{n=1}^{\infty}a_{p,1,n}(v)q^{n/T}
\end{equation}
defines a function $\a: V\to\C$ which is not identically
zero. Because $S(v,\tau)$ is linear in $v,$ $\alpha$ is a linear
functional on $V.$
\begin{lem}\label{la9.3} $\alpha$ vanishes on $O_g(V).$ 
\end{lem}

\pf We know that $S$ vanishes on $O(g,h),$ hence on the elements
(\ref{m5.1})-(\ref{m5.4}). Using Lemma \ref{l9.3} leads to the required vanishing conditions. For example, if $\mu\ne 1$ in the notation of
Lemma \ref{l9.3} then 
\begin{equation}\label{9.2}
 \sum_{k=0}^{\infty}Q_k(\mu,\la,\tau)S(v[k-1]w)=0.
\end{equation}
This identity holds if $S$ is replaced by $S_p,$ and 
then if $S_p$ is in turn replaced by $S_{p,1}.$ Hence
\begin{equation}\label{9.3}
\sum_{k=0}^{\infty}Q_k(\mu,\la,\tau)(\alpha(v[k-1]w)+\sum_{n=1}^{\infty}a_{p,1,n}(v[k-1]w)q^{n/T})=0.
\end{equation}
The constant term in (\ref{9.3}), necessarily zero, is equal to
 $\a(-v\circ_g w)=0$ by Lemma \ref{l9.3}, so that 
$\alpha$ vanishes on $v\circ_g w$ if $gv\ne v.$ The other
vanishing conditions follow similarly. \qed

\bigskip

To complete the proof of Theorem \ref{t9.1} we must show that there is some
non-zero element $S$ of ${\cal C}_1(g,h)$ (for suitable $h$). For then we know 
that the function $\alpha$ is non-vanishing on $V$ but vanishes on $O_g(V),$
whence $V\ne O_g(V)$ and $A_g(V)=V/O_g(V)$ is non-zero, as required.

Consider ${\cal C}_1(1,g):$ $V$ is itself a $g$-stable simple $V$-module,
so that the corresponding trace function $T_V(v,(1,g),q)$ lies in 
${\cal C}_1(1,g)$ by Theorem \ref{t8.1}. So ${\cal C}_1(1,g)\ne 0,$ and
since $\left(\begin{array}{cc} 0 & -1\\ 1 & 0\end{array}\right)$ induces 
a linear isomorphism between ${\cal C}_1(1,g)$ and ${\cal C}_1(g,1)$ by Theorem
\ref{t5.4} then also ${\cal C}_1(g,1)\ne 0.$ This completes the proof
of Theorem \ref{t9.1}.

\section{The main theorems}
\setcounter{equation}{0}
We continue to use the notation introduced in Section 5. The next theorem is decisive.
\begin{th}\label{t10.1}
Suppose that $V$ is $g$-rational and satisfies condition $C_2,$ and
let $M^1,$ ..., $M^m$ be all of the inequivalent, simple, $h$-stable,
$g$-twisted $V$-modules. Let $T_1,...,T_m$ be the corresponding trace
functions (\ref{8.7}). Then $T_1,...,T_m$ form a basis of ${\cal
C}_1(g,h).$ 
\end{th}

There are several important corollaries.
\begin{th}\label{t10.2} 
Suppose that $V$ is rational and satisfies condition $C_2.$ Suppose
further that the group $\<g,h\>$ generated by $g$ and $h$ is cyclic
with generator $k.$ Then the dimension of ${\cal C}_1(g,h)$ is equal
to the number of inequivalent, $k$-stable, simple $V$-modules. In particular, the number of inequivalent, simple $g$-twisted
$V$-modules is at most equal to the number of $g$-stable, simple $V$-modules,
with equality if $V$ is $g$-rational.
\end{th}

Notice that parts (ii) and (iii) of Theorem \ref{thmm1} are
included in Theorem \ref{t10.2} together with Theorem \ref{t9.1}.
Recall next that a simple \voa $V$ is called {\em holomorphic} in case $V$ 
is rational and if $V$ is the unique simple $V$-module.

\begin{th}\label{t10.3} Suppose that $V$ is holomorphic and satisfies 
condition $C_2.$ For each automorphism $g$ of $V$ of finite order, there is
a unique simple $g$-twisted $V$-module $V(g).$ Moreover if $\<g,h\>$
is cyclic then ${\cal C }_1(g,h)$ is spanned by $T_{V(g)}(v,g,h,q).$
\end{th}

This establishes Theorem \ref{thmm2} (i). 
First we show how Theorems \ref{t10.2} and \ref{t10.3} follow from Theorem
\ref{t10.1}. In the situation of Theorem \ref{t10.2} we have
$\<g,h\>=\<k\>.$ Since $V$ is rational, Theorem \ref{10.1} tells us that 
${\cal C}_1(1,k)$ has a basis consisting of the trace functions $T_M(v,1,k,q)$
where $M$ ranges over the inequivalent, $k$-stable, simple $V$-modules. 

We can find $\gamma\in SL(2,\Z)$ such that $(g,h)\gamma=(1,k).$ By the theorem
of modular-invariance, $\gamma$ induces a linear isomorphism from
${\cal C}_1(g,h)$ to ${\cal C}_1(1,k).$ Theorem 
\ref{t10.2} follows from this together with Theorem \ref{t8.7}.

As for Theorem \ref{t10.3}, since $V$ is holomorphic it is certainly rational, so Theorem \ref{t10.2} applies. So if $\<g,h\>=\<k\>$ is cyclic then $\dim
{\cal C}_1(g,h)$ is equal to 1 since $V$ is the only simple $V$-module and it is certainly $k$-stable. By Theorem \ref{t9.1}, $V$ has at least one simple
$g$-twisted $V$-module, call it $V(g),$ and from Theorem
\ref{t8.7} there can be no more than one since $\dim {\cal C}_1(g,1)=1.$ So
$V(g)$ is unique, hence $h$-stable whenever $gh=hg.$ So $T_{V(g)}(v,g,h,q)$ spans ${\cal C}_1(g,h)$ by Theorem \ref{t8.1}. This establishes
Theorem \ref{t10.3}.

Turning to the proof of Theorem \ref{t10.1}, we consider first an
arbitrary function $S\in{\cal C}_1(g,h).$ We have seen in Section 9
that
$S$ can be represented as
\begin{equation}\label{10.1}
S(v,\tau)=\sum_{i=0}^p(log q_{1\over T} )^iS_i(v,\tau)
\end{equation}
for fixed $p$ and all $v\in V,$ with each $S_i$ satisfying (\ref{6.13})-(\ref{6.15}). We will prove
\begin{prop}\label{p10.4} Each $S_i$ is a linear combination of the functions $T_1,...,T_m.$
\end{prop}
\begin{prop}\label{p10.5} $S_i=0$ if $i>0.$
\end{prop}

Theorem \ref{t10.1} obviously follows from these two
propositions. First we show that Proposition \ref{p10.5} follows from
Proposition \ref{p10.4}. To this end, pick $v\in V$ such that
$gv=\mu^{-1}v, hv=\la^{-1}v.$ It suffices to show that
$pS_p(v,\tau)=0.$ If $(\mu,\la)\ne (1,1)$ this follows from (C3) and
(\ref{m5.3}), so we may assume $gv=hv=v.$

Set
\begin{equation}\label{10.2}
w=L[-2]v-\sum_{k=1}^{\infty}E_{2k}(\t)L[2k-2]v. 
\end{equation}
{}From (\ref{m5.10}) and (\ref{m5.11}) we get
\begin{equation}\label{10.3}
S(w,\t)=\frac{q_{1\over T} }{T}\frac{d}{dq_{1\over T} }S(v,\tau).
\end{equation}

Now Proposition \ref{p10.4} combined with Theorem \ref{t8.1} tells us that
(\ref{10.3}) is satisfied by each $S_i.$ Then we calculate that 
\begin{equation}\label{10.4}
S(w,\tau)=\sum_{i=0}^p(\frac{i}{T}(\log q_{1\over T} )^{i-1}S_i(v,\tau)
+(\log q_{1\over T} )^iS_i(w,\tau)).
\end{equation}

We may identify the parts of (\ref{10.1}) which involve a given power $(\log q_{1\over T} )^i.$ Taking $i=p-1,$ we see that 
$$S_{p-1}(w,\tau)=\frac{p}{T}S_p(v,\tau)+S_{p-1}(w,\tau),$$
so that $pS_p(v,\tau)=0,$ as desired.

We turn our attention to the proof of Proposition \ref{p10.4}. We assume
without loss that $S_p\ne 0$ and that each $S_{p,j}\ne 0$ (cf. (\ref{6.13}).
We are then in the situation that was in effect in Section 9. We adopt the notation (\ref{9.1}). It was shown that 
$\a: V\to\C$ vanishes on $O_g(V),$ and thus defines a linear functional
$$\alpha: A_g(V)\to \C.$$
We continue this line of reasoning, and now prove
\begin{lem}\label{l10.6}
Suppose that $u,v\in V$ and satisfy $hu=\rho u,hv=\sigma v,$ $\rho,\sigma\in \C.$ Then 
\begin{equation}\label{10.5}
\a(u*_g v)=\rho\delta_{\rho\sigma,1}\a(v*_gu).
\end{equation}
\end{lem}

\pf We may assume that $gu=\xi u$ and $gv=\nu v$ for scalars $\xi,\nu.$ If $\xi$ or $\nu$ is not equal to 1 then $u$ (resp. $v$) lies in $O_g(V)$
(cf. Lemma 2.1 of [DLM3]), whence so too do $u*_gv$ and $v*_gu$ by Theorem
\ref{t3.4}. So in this case both sides of (\ref{10.5}) are equal to
0. So we may assume $gu=u, gv=v.$

Similarly, $u*_gv$ is an eigenvector for $h$ with eigenvalue $\rho\sigma,$ so
if $\rho\sigma\ne 1$ then $u*_gv$ and $v*_gu$ lie in $O(g,h)$ by (\ref{m5.3}). Then $S(u*_gv)=S(v*_g u)=0$ by (C3), which again leads to both sides of
(\ref{10.5}) being 0. So we may assume that $\rho\sigma=1$ and try to prove that $\a(u*_gv)=\rho \a(v*_g u).$

Now we know from [Z] (also see Lemma 2.2 (iii) of [DLM3]) that if $V^g$ is the
space of $g$-invariants of $V$ then for $u$ homogeneous
$$ u*_gv-v*_gu\equiv \Res_{z}Y(u,z)v(1+z)^{\wt u-1}\ (\mod O(V^g))$$
Using (\ref{g2.14}) we get 
\begin{equation}\label{10.6}
   u*_gv-v*_gu\equiv \sum_{i=0}^{\infty}{\wt u-1\choose i}u(i)v\equiv
u[0]v\ (\mod O(V^g)).
\end{equation}

Now certainly $O(V^g)\subset O_g(V)$ (loc. cit.), and if $\rho=1$ then
$u[0]v\in O(g,h)$ by (\ref{m5.1}). So in this case (\ref{10.6}) leads to 
$\a(u*_gv-v*_gu)=0$ as desired.

So we may take $\rho \ne 1.$ In this case we follow the calculation of Lemma
\ref{l9.3}. Bearing in mind that $gu=u$ and $hv=\rho v$ with
$\rho\ne 1,$ we see from the proof of Lemma \ref{l9.3} (b)
that the constant term  of 
$$\sum_{k=0}^{\infty}Q_k(1,\rho^{-1},q)u[k-1]v$$
is 
$$-u*_gv+\frac{1}{1-\rho^{-1}}u[0]v.$$
Since $S$ vanishes on $\sum_{k=0}^{\infty}Q_k(1,\rho^{-1},q)u[k-1]v\in O(g,h),$
this shows that 
$$\a(u*_g v)=\frac{1}{1-\rho^{-1}}\a(u[0]v).$$
However (\ref{10.6}) still applies, so that
$$\a(u*_g v)=\frac{1}{1-\rho^{-1}}(\a(u*_gv)-\a(v*_gu)),$$
which is equivalent to the desired result. \qed

We will need 
\begin{lem}\label{l10.7} Let $A$ be a finite-dimensional semi-simple algebra 
over $\C$ with decomposition
$A=\oplus_{i\in I}A_i$ into simple components.

Let $h:A\to A$ be an automorphism of $A$ of finite order, and suppose
that $F: A\to\C$ is a linear map which satisfies
\begin{equation}\label{10.7}
F(ab)=\rho\delta_{\rho\sigma, 1}F(ba)
\end{equation}
whenever $ha=\rho a$ and $hb=\sigma b,$ $\rho,\sigma\in\C.$ Then $F$ 
can be written as a linear combination with scalars $\alpha_j:$
\begin{equation}\label{10.8}
F(a)=\sum_{j\in J}\alpha_j\tr_{W_j}(a \gamma_j)
\end{equation}
where in (\ref{10.8}), $\{A_j\}_{j\in J}$ ranges over the $h$-invariant
simple components of $A,$ $W_j$ is the simple $A$-module such that
$A_jW_j=W_j,$ and $\gamma_j\in A_j^*$ satisfies 
$ha=\gamma_j a\gamma_j^{-1}$ for $a\in A_j.$
\end{lem}

\noindent{\bf Remark:} The existence of $\gamma_j$ is the Skolem-Noether
Theorem.

\pf Proceed by induction on the order of $h$ and cardinality
of $I.$ The group $\<h\>$ permutes the $A_i$ among themselves, and the
conditions of the lemma apply to any $h$-invariant sum of $A_i.$ So we
may assume that $\<h\>$ is transitive in its action on the $A_i.$

First assume that 
there are at least two components. Then there are no $h$-invariant
components, so we must show that $F=0$ in this case. If $\sigma\ne 1$ and 
$hb=\sigma b,$ taking $a=1$ in (\ref{10.7}) shows that $F(b)=0.$ So we only
need show that $F$ is zero on the algebra $A^h$ of $h$-invariants of $A.$

If there is $1\ne k\in\<h\>$ such that $k$ fixes each $A_i$ then the algebra
of $k$-invariants $B$ is a semi-simple
algebra admitting $\<h\>/\<k\>.$ By induction we see that $F(B)=0,$ so we are
done as $A^h\subset B.$ So without loss there is no such $k.$ So
$h$ has order $|I|,$ the number 
of components. Thus if $A_0$ is the first component then we may set
$A_i=h^i A_0,$ $0\leq i\leq |I|-1.$ Then 
$$A^h=\{\sum_{i=0}^{|I|-1}h^ix|x\in A_0\}\simeq A_0.$$

By (\ref{10.7}) $F(ab)=F(ba)$ for $a,b\in A^h,$ so $F$ is a trace function
, i.e., $F(a)=\alpha \tr_W a$ for some $\a\in\C,$ $W$ the simple $A^h$-module.
So we must show $F(1_A)=0.$ Let $\l$ be an $|I|$-th root of unity, let $u,v\in A_0$ be units such that $uv=1_{A_0},$ and let 
$$a=\sum_{i=0}^{|I|-1}\l^ih^i(u), b=\sum_{i=0}^{|I|-1}\l^{-i}h^i(v).$$
Then $ba=ab=\sum_{i=0}^{|I|-1}h^i(u)h^i(v)=\sum_ih^i(uv)=1_A.$ On the
other hand $ha=\la^{-1}a, hb=\la b$ and $\la\ne 1.$ So (\ref{10.7})
yields $F(1_A)=F(ab)=F(ba)=\la F(ab).$ So $F(1_A)=0$ as desired.

This reduces us to the case that $A$ is itself a simple
algebra. Pick $\gamma\in A^*$ such that $h(a)=\gamma a\gamma^{-1}$ 
and consider $F_1: A\to\C$ defined by $F_1(a)=F(a\gamma^{-1}).$
If $ha=\rho a$ and $\rho\ne 1$ then $F(a\gamma^{-1})=0$ by (\ref{10.7}),
so $F_1(a)=0$ for such $a.$ On the other hand we get 
for $hb=\sigma b,$
\begin{eqnarray*}
& & \ \ \ F_1(ab)=F(ab\gamma^{-1})=\rho\delta_{\rho\sigma,1}F(b\gamma^{-1}a)\\
& &=\rho\delta_{\rho\sigma,1}F(b\gamma^{-1}a\gamma\gamma^{-1})
=\delta_{\rho\sigma,1}F(ba\gamma^{-1})=\delta_{\rho\sigma,1}F_1(ba).
\end{eqnarray*}
{}From this we conclude that $F_1(ab)=F_1(ba)$ for
all $a,b\in A.$ So $F_1$ is a trace function $F_1(a)=\a\tr_Wa,$ so
that $F(a)=\a\tr_Wa\gamma.$ This completes the proof of the
lemma. \qed

Now we return to the situation of Lemma \ref{l10.6}. From (\ref{g3.22})
$h$ induces an automorphism of $A_g(V)$ via $h:v\mapsto hv,$ and since
$V$ is $g$-rational, then $A_g(V)$ is semi-simple and Lemma
\ref{l10.7} applies.  From the discussion in Section 3, the
$h$-invariant components of $A_g(V)$ correspond precisely to the
$h$-invariant simple $A_g(V)$-modules, and these correspond to the
$h$-invariant simple $g$-twisted $V$-modules. For such a simple
$A_g(V)$-module $\Omega$ we have $\phi(h)o(v)\phi(h)^{-1}=o(hv)$ (cf. (\ref{8.1})), 
$o(v)$ being the corresponding zero mode (\ref{g3.18}). Also
$o(hv)=\gamma o(v)\gamma^{-1}$ if $\gamma$ represents $h$ in the sense of Lemma 
\ref{l10.7}. So $\gamma$ and $\phi(h)$ differ by a scalar when considered
as operators on $\Omega.$ By Lemmas \ref{l10.6} and \ref{l10.7} we get
\begin{lem}\label{l10.8} The linear function $\a: A_g(V)\to \C$ can be
represented in the form
\begin{equation}\label{10.9}
\a(v)=\sum_{j}\a_j\tr_{\O(M^j)}o(v)\phi(h)
\end{equation}
where $\a_j$ are scalars and the spaces $\O_{M^j}$ range over the top 
levels of the $h$-invariant simple $g$-twisted $V$-modules $M^j.$
\end{lem}

Recall that we have
\begin{equation}\label{10.10}
S_p(v,\tau)=\sum_{j=1}^bq^{\la_{p,j}}S_{p,j}(v,\tau)
\end{equation}
with $S_{p,1}$ as in (\ref{9.1}).

\begin{lem}\label{l10.9} Suppose that $\a_j\ne 0$ in (\ref{10.9}). Then
the conformal weight of the corresponding $g$-twisted module $M^j$ is
equal to $\la_{p,1}+c/24.$ 
\end{lem} 

\pf We use the method of
proof of Proposition \ref{p10.5} once more. For $v\in V,$ let $w=w(v)$ be as
in (\ref{10.2}). Thus (\ref{10.3}) holds whenever $S\in {\cal C}_1(g,h).$

Applying (\ref{10.3}) with $S=T_{M^j}$ (cf. (\ref{8.4})) and considering 
leading terms yields
\begin{equation}\label{10.11}
\tr_{\Omega(M^j)}o(w)\phi(h)=(\l_j-c/24)\tr_{\O(M^j)}o(v)\phi(h)
\end{equation}
where $\la_j$ is the conformal weight of $M^j.$ Similarly applying 
(\ref{10.3}) to $S$ itself and considering the leading term of $S_p$ yields
\begin{equation}\label{10.12}
\a(w)=\l_{p,1}\a(v).
\end{equation}
Using (\ref{10.5}), we find that for $v\in V,$
\begin{eqnarray*}
& &\ \ \ \la_{p,1}\sum_j\a_j\tr_{\O(M^j)}o(v)\phi(h)
=\l_{p,1}\a(v)=\a(w)\\
& &=\sum_j\a_j\tr_{\Omega(M^j)}o(w)\phi(h)
=\sum_j\a_j(\l_j-c/24)\tr_{\Omega(M^j)}o(v)\phi(h).
\end{eqnarray*}
The linear independence of characters of $A_g(V)$ implies that
$\l_{p,1}\a_j=\a_j(\l_j-c/24).$ The lemma follows. \qed

We are ready for the final argument. We have in the previous notation
\begin{eqnarray*}
& &q^{\l_{p,1}}S_{p,1}(v,\tau)=q^{\l_{p,1}}(\sum_j\a_j\tr_{\O(M^j)}o(v)\phi(h)
+\sum_{n=1}^{\infty}a_{p,1,n}(v)q^{n/T})\\
& &=\sum_jq^{\l_j-c/24}\a_j\tr_{\O(M^j)}o(v)\phi(h)+q^{\l_{p,1}}\sum_{n=1}^{\infty}a_{p,1,n}(v)q^{n/T}).
\end{eqnarray*}
Now also
$$ T_{M^j}(v,\tau)=q^{\l_j-c/24}(\tr_{\O(M^j)}o(v)\phi(h)+\sum_{n=1}^{\infty}
\tr_{M^j_{\l_j+n/T}}o(v)\phi(h)q^{n/T}).$$

So we see that the function $$S'(v,\tau)=S(v,\tau)-(\log
q_{1\over T})^p\sum_j\a_jT_{M^j}(v,\tau)$$ again has the form
(\ref{6.12})-(\ref{6.15}), but the leading term of the piece
corresponding to $S_p$ now has a higher degree than $S_p$ itself. We
now continue the argument, replacing $S$ with $S'$ and $S_p$ with
$S_p'.$ We find, since each $T_{M^j}$ already lies in ${\cal C}_1(g,h),$ and
since there are only finitely many $M^j,$ that indeed $S_p$ is a
linear combination of $T_{M^j}.$ But our argument applies equally well
to each $S_i,$ so each $S_i$ is a linear combination of $T_{M^j}.$
This completes the proof of Proposition \ref{p10.4}.
\section{Rationality of central charge  and conformal\newline  weights}
\setcounter{equation}{0}

Recall from (\ref{3.10'}) that a simple $g$-twisted $V$-module $M$ has
grading of the form $M=\oplus_{n=0}^{\infty}M_{\la+n/T}$ for some
$\la\in\C$ called the conformal weight of $M.$ We will show that,
under suitable rationality conditions on $V,$ the conformal weight $\la$
of $M$ is a rational number. We prove even more, namely
\begin{th}\label{t11.1}
Suppose that $V$ is a holomorphic \voa which satisfies condition
$C_2,$ and let $g\in\Aut V$ have finite order. Let $V(g)$ be the unique simple $g$-twisted
$V$-module whose existence is guaranteed by Theorem \ref{t10.3}. Then
the conformal weight of $V(g)$ is rational, and the central charge $c$ of
$V$ is also rational.
\end{th}

\begin{th}\label{t11.2} Suppose that $V$ is a \voa which satisfies condition 
$C_2,$ and let $g\in\Aut V$ have finite order. Suppose that $V$ is
$g^i$-rational for all integers $i.$ Then each simple $g^i$-twisted
$V$-module has rational conformal weight, and the central charge $c$ of $V$
is rational. 
\end{th}

\begin{th}\label{t11.3} Suppose that $V$ is a rational \voa which satisfies 
condition $C_2.$ Then each simple $V$-module has rational conformal weight, and the central charge of $V$ is rational.
\end{th}

These theorems complete the proofs of Theorems \ref{thmm1} and
\ref{thmm2}. Note that Theorem \ref{t11.3} is simply a restatement of Theorem \ref{t11.2}
in the special case that $g=1.$ We will prove Theorems \ref{t11.1} and
\ref{t11.2} simultaneously. Indeed, at this point
in the paper the proof follows from ideas in
a paper of Anderson and Moore [AM]: we have only to assemble the
relevant facts.

First observe that to prove Theorem \ref{t11.2} it  suffices to show that 
each simple $g$-twisted $V$-module has rational conformal weight, and that
$c$ is rational. With this in mind, let $f(q)$ be one of the following
$q$-expansion: $q^{-c/24}\sum_{n\geq n_0}(\dim V_n)q^n$ where 
$V=\oplus_{n\geq n_0}V_n;$   $q^{\la-c/24}\sum_{n\geq n_0}(\dim V(g)_{\la+n/T})q^n$ with $\la$ the conformal weight of $V(g),$ where $V(g)$ is either 
the unique simple $g$-twisted $V$-module in the situation of Theorem
\ref{t11.1}, or any simple $g$-twisted $V$-module in the situation of Theorem 
\ref{t11.2}.

Let $U$ be the $SL(2,\Z)$-module of holomorphic functions on $\frak h$
generated by $f(q).$ In each case $U$ is a finite-dimensional $\C$-linear 
space, and the elements of $U$ have $q$-expansions in (not necessarily
rational) powers of $q.$ This assertion follows from Theorems \ref{t10.1}
and \ref{t10.3}. This puts us in the position of being able to apply
methods and results of Anderson and Moore (loc.cit.) . The argument proceeds as follows.

Define by $\l=\l(\tau)$ the usual Picard function which generates the field
of rational functions on the compactification of $\frak h/\Gamma(2).$ With
$E=\frac{d}{d\la},$ there are unique meromorphic functions $k_i$ such that
$U$ is precisely the space of solutions of the differential equation
\begin{equation}\label{11.1}
E^ny+\sum_{i=0}^{n-1}k_iE^iy=0.
\end{equation}
The $k_i$ are then in $\C(\l)$ (Proposition 1 of (loc.cit.)). 

For a given $\phi\in\Aut(\C),$ and for $r(q)\in U,$ let $r^{\phi}$ be
as defined in (loc.cit.). By the Frobenius-Fuchs theory, 
the $r^{\phi}$ are then
$q$-expansions of the solutions of the $\phi$-transform of
(\ref{11.1}), namely
\begin{equation}\label{11.2}
E^ny+\sum_{i=0}^{n-1}k_i^{\phi}E^iy=0.
\end{equation}

We claim that the solutions of (\ref{11.2}) also afford
a representation of $SL(2,\Z).$ First note that since each $k_i$ lies 
$\C(\l)$ then the actions of $\phi$ and $\gamma\in SL(2,\Z)$ on
$\C(\l)$ commute: this follows from the well-known formulae for the 
action of the modular group on $\la.$ Then if $y|\gamma$ is the 
$\gamma$-image of a solution $y$ of (\ref{11.1}) we find that $(y|\gamma)^{\phi}|\gamma^{-1}$ is a solution of (\ref{11.2}). The claim follows
{}from this observation.

Now $f(q)\in U$ has the form  
$$f=q^{\la-c/24}\sum_{n\geq N}a_nq^{n/T}$$
where $\l=0$ and $T=1$ in the first case, and where $a_n\in\Z$ in all cases.
Then 
$$f^{\phi}=q^{\phi(\la-c/24)}\sum_{n\geq N}a_nq^{n/T}$$
i.e.,
\begin{equation}\label{11.3}
f^{\phi}=q^{\phi(\la-c/24)-(\la-c/24)}f.
\end{equation}

One now applies $S$ to both sides of (\ref{11.3}) to obtain
\begin{equation}\label{11.4}
f^{\phi}|S=e^{-\alpha/\tau}f|S
\end{equation}
where $\alpha=2\pi i(\phi(\la-c/24)-(\la-c/24)).$ On the other hand, we showed
above that both $f^{\phi}|S$ and $f|S$ have $q$-expansions. This leads to a 
contradiction by using the limit argument of (loc.cit.) unless $\alpha=0.$ As this
holds for all $\phi,$ we conclude that $c/24\in\Q$ and $\la-c/24\in \Q,$ which
completes the proofs of the theorems.
\bigbreak
Let us formalize the situation which prevails in case $V$ is a holomorphic \voa
which satisfies condition $C_2$ and is  equipped with a finite 
group $G$ of automorphisms. Let $V(g)$ be the unique simple $g$-twisted
$V$-module (Theorem \ref{t10.3}) for $g\in G,$ and let $C(g)=\{h\in G|gh=hg\}$ 
be the {\em centralizer} of $g$ in $G.$  

According to Theorem \ref{t11.1} the grading of $V(g)$ has the form 
$V(g)=\oplus_{n=0}^{\infty}V(g)_{\l+n/T}$ where $g$ has order $T$ and $\l=\l(g)\in\Q.$ As $V(g)$ is unique, it admits a (projective) representation
of $C(g)$ (cf. (\ref{3.25})), so for any $h\in C(g)$ we may consider the trace
function $T_{V(g)}(v,g,h,q).$ By Theorem \ref{t10.3} this trace function spans
the 1-dimensional space ${\cal C}_1(g,h)$ if $\<g,h\>$ is cyclic. The shape
of the trace function is as in (\ref{8.6}), with $\l-c/24\in\Q.$ Thus it
has a $q$-expansion with {\em rational} powers of $q$ of bounded denominator,
and is holomorphic as a function on $\frak h.$

Note also that by [Z] (or Theorem \ref{11.1} and a short argument) we know
that $c$ is in fact an integer divisible by 8.

We now assume that $\<g,h\>$ is cyclic and fix $v\in V_{[k]}.$ It follows from Theorem \ref{t5.4} that
$T_{V(g)}|_k\gamma (v,g,h,\t)=(c\t+d)^{-k} T_{V(g)}(v,g,h,\gamma\tau)$ lies in
${\cal C}_1((g,h)\gamma,\t)$ and hence is a scalar multiple of 
$T_{V(g^ah^c)}(v,(g,h)\gamma,\t).$ Here $\gamma=\left(\begin{array}{cc}
a &b\\c &d\end{array}\right)$ lies in $SL(2,\Z).$ Thus there are scalars
$\sigma(g,h,\gamma)$ such that the following holds:
\begin{equation}\label{11.5}
(c\t+d)^{-k} T_{V(g)}(v,g,h,\gamma\tau)=\sigma(g,h,\gamma)T_{V(g^ah^c)}(v,(g,h)\gamma,\t).
\end{equation}
Equation (\ref{11.5}) together with the rationality of the corresponding 
$q$-expansions says precisely that each $T_{V(g)}(v,g,h,\t)$ is a
{\em generalized modular form of weight} $k$ in the language of [KM]. Note
that Theorem \ref{thmm4} follows from these results. Theorem \ref{thmm3}
follows in the same way. 

A case of particular interest is when we take $v$ to be the vacuum element
${\bf 1},$ in which case $k=0.$ One
sets 
\begin{equation}\label{11.6}
Z(g,h)=T_{V(g)}({\bf 1},g,h,\tau)
\end{equation}
in this case. This is essentially the graded trace of $\phi(h)$ on the 
$g$-twisted module $V(g),$ sometimes called a {\em partition function}
or {\em McKay-Thompson series.} In this case, we have proved
\begin{th}\label{t11.4} Let $V$ be a holomorphic \voa which satisfies condition
$C_2,$ and let $G$ be a finite group of automorphisms of $V.$ For
each pair of commuting elements $(g,h)$ which generates a cyclic group,
$Z(g,h)$ is a generalized modular function (i.e., of weight zero) 
which is holomorphic on
$\frak h$ and satisfies 
\begin{equation}\label{11.7}
\gamma: Z(g,h)\mapsto \sigma(g,h,\gamma)Z((g,h)\gamma)
\end{equation}\
for $\gamma\in SL(2,\Z).$
\end{th}

\section{Condition $C_2$}
\setcounter{equation}{0}

In order to be able to apply the preceding results to known vertex operator 
algebras, we
need verify that condition $C_2$ is satisfied. We do this for
some of the best known rational \voas in this section. Refer to Section 3
for the definition of condition $C_2.$

\begin{lem}\label{l12.1}
If $V$ is a \voa and $M$ is $V$-module, then $C_2(M)$ contains 
$v(-n)M$ for all $v\in V$ and $n\geq 2.$
\end{lem}

\pf This follows from the definition (\ref{3.26}) together with the equality
$$(L(-1)^mv)(-2)=(m+1)!v(-m-2).\ \ \ \qed $$
 
\begin{lem}\label{l12.2} Let $V_1,...,V_k$ be \voas such that for each $i,$
all simple $V_i$-modules satisfy condition $C_2.$ Then the same is
true for the tensor product \voa $V_1\otimes \cdots \otimes V_k.$
\end{lem} 

\pf See [FHL] for tensor product \voas and their modules. We
may assume that $k=2.$ One knows (loc.cit.) that the simple 
$V_1\otimes V_2$-modules are precisely those of the form $M_1\otimes M_2$ with
$M_i$ a simple $V_i$-module.

If $v\in V_1$ then $(v\otimes {\bf 1})(-2)=v(-2)\otimes id,$ from 
which it follows that $C_2(M_1\otimes M_2)$ contains 
$C_2(M_1)\otimes M_2.$ Similarly it contains $M_1\otimes C_2(M_2).$ 
The lemma follows immediately. \qed

Now we discuss condition $C_2$ for the most well-known rational vertex operator
algebras, namely,

(i) The \voa $L(c_{p,q},0)$ associated with the (discrete series) simple
Virasoro algebra $Vir$-module of highest weight 0 and central charge 
$c=c_{p,q}=1-\frac{6(p-q)^{2}}{pq}$ ([DMZ], [FZ], [Wa]).

(ii) The moonshine module $V^{\natural}$ ([B1], [FLM2]).

(iii) The \voa $V_L$ associated with a positive definite even lattice $L$
([B1], [D1], [FLM2]).

(iv) The \voa $L(k,0)$ associated to a $\hat \frak g$-module of highest 
weight $0$ and positive integral level $k,$ $\frak g$ a simple Lie algebra
([DL],[FZ],[Li]).

\begin{lem}\label{l12.3}
$L(c_{p,q},0)$ satisfies condition $C_2.$
\end{lem}

\pf Set $L=L(c_{p,q},0).$ It is a quotient of the corresponding Verma module
$M=M(c_{p,q},0)$ and we have $M\simeq U(Vir_{-})\cdot 1$ (cf. [FZ]) where 
$Vir_{-}=
\oplus_{n=1}^{\infty}\C L(-n)$ and the $L(n)$ are the usual generators
of the Virasoro algebra $Vir.$

Now $Y(\omega,z)=\sum_{n\in\Z}L(n)z^{-n-2},$ so that $C_2(V)$ contains
$L(-n)L$ for all $n\geq 3$ by Lemma \ref{l12.1}. 

We have $L=M/J$ and $J$ contains two singular vectors [FF]. The first
is $L(-1)\cdot 1,$ which shows that $C_2(L)$ contains $U(Vir_{-})L(-1)1.$
{}From this we see that $L=C_2(L)+\sum_{k=0}^{\infty}\C L(-2)^k1.$

The second singular vector has the form 
\begin{eqnarray}\label{12.1}
v=L(-2)^{pq}{\bf 1}+\sum a_{n_{1},\cdots, n_{r}}L(-n_1-2)\cdots L_(-n_r-2)1
\end{eqnarray}
where the sum ranges over certain $(n_{1},\cdots, n_{r})\in {\Z}_{+}^{r}$
with $n_{1}+\cdots +n_{r}\ne0,$ $a_{n_1,,..,n_r}\in\C$ (cf. equation
(3.11) of [DLM2]). From the previous paragraph we see that the terms
under the summation sign in (\ref{12.1}) each lie in $C_2(L),$ whence also
$L(-2)^{pq}1$ lies in $C_2(L).$ By Lemma \ref{l3.7a}, $C_2(L)$ is
invariant under $L(-2).$  We conclude that $L=C_2(L)+\sum_{k=0}^{pq-1}\C L(-2)^k1,$ and the proposition follows. \qed

The following result was stated without proof in [Z]

\begin{prop}\label{p12.4} The moonshine module $V^{\natural}$ satisfies
condition $C_2.$
\end{prop}

\pf Let $U$ be the tensor product $L(\frac{1}{2},0)^{\otimes 48}.$ It is 
shown in [DMZ] that $V^{\natural}$ contains a sub \voa isomorphic to $U.$ 
Moreover when considered as a $U$-module, $V^{\natural}$ is a direct
sum of finitely many simple $U$-modules.

Suppose that each simple module for $L(\frac{1}{2},0)$ satisfies condition
$C_2.$ Then this is true also for $U$ by Lemma \ref{l12.2}, so that
the space spanned by $u(-2)v$ for $u\in U$ and $v\in V^{\natural}$ already has finite codimension in $V^{\natural}.$ So it suffices to show that 
the simple $L(\frac{1}{2},0)$-modules indeed satisfy condition $C_2.$

The proof of this later assertion is similar to that of the last proposition.
Apart from $L=L(\frac{1}{2},0)$ itself there are 
just two other simple modules for $L,$ namely $L(\frac{1}{2},\frac{1}{2})$
and $L(\frac{1}{2},\frac{1}{16})$ [DMZ]. Let $M(\frac{1}{2},h)$ $(h=\frac{1}{2},
\frac{1}{16})$
be the corresponding Verma module with
$L(\frac{1}{2},h)=M(\frac{1}{2},h)/J_h.$ As in the previous
proposition we have $L(-n)L(\frac{1}{2},h)\subset
C_2(L(\frac{1}{2},h))$ for $n\geq 3,$ so that $L(\frac{1}{2},h)=
C_2(L(\frac{1}{2},h))+\sum_{a,b\geq 0}\C L(-2)^aL(-1)^b v_h$ where
$v_h$ is a highest weight vector. Moreover $J_h$ contains two singular
vectors (cf. [FF]). One of them is
$(L(-2)-\frac{3}{4}L(-1)^2)v_{\frac{1}{2}}$ or
$(L(-2)-\frac{4}{3}L(-1)^2)v_{\frac{1}{16}}$ (cf. [DMZ]), from which
we conclude that
$L(\frac{1}{2},h)=C_2(L(\frac{1}{2},h))+\sum_{b\geq 0}\C L(-1)^bv_h.$ Because
of the existence of a second singular vector we see that $C_2(L(\frac{1}{2},h))$
necessarily contains $L(-1)^bv_h$ for big enough $b,$ whence our claim follows.
\qed

The reader is referred to [FLM2] for the construction
of $V_L$ and associated notation, which we use in the next result. 

\begin{prop}\label{p12.5} $V_L$ satisfies condition $C_2.$
\end{prop}

\pf Set $H=L\otimes_{\Z}\C.$  Then it is easy to see that
\begin{equation}\label{12.2}
V_L=C_2(V_L)+S(H\otimes t^{-1})\otimes \C\{L\}.
\end{equation}
Let $0\ne \alpha,\gamma\in L$ and set $\beta=\gamma-\alpha.$ Then
$C_2(V_L)$ contains $u(-2)v$ where we take $u=L(-1)^k\iota(a)$
$(k\geq 0)$ and $v=\iota(b)$ where $a,b,c\in\hat L$ such that $\bar a=\alpha,$ $\bar b=\beta$ and $c=ab.$ Then $\bar c=\gamma.$ So $C_2(V_L)$ contains
$$\Res_zz^{-k-2}Y(\iota(a),z)\iota(b)
=\Res_zz^{-k-2}
\exp\left(\sum_{n<0}\frac{-\a(n)}{n}z^{-n}\right)
\exp\left(\sum_{n>0}\frac{-\a(n)}{n}z^{-n}\right)az^{\alpha}\iota(b).$$
Now $az^{\alpha}\iota(b)=z^{\<\a,\b\>}\iota(c).$ Using
(\ref{12.2}) we now see that $C_2(V_L)$ contains 
$$\Res_zz^{-k-2}Y(\iota(a),z)\iota(b)e^{\alpha(-1)z}z^{\<\a,\b\>}\iota(v)$$
and we conclude that $C_2(V_L)$ contains
\begin{equation}\label{12.3}
\a(-1)^{1+k-\<\a,\beta\>}\iota(c)
\end{equation}
whenever $k\geq 0$ and $k\geq 1-\<\a,\b\>.$

So if $\<\a,\b\>\geq 1$ we may choose $k$ appropriately to see that
$C_2(V_L)$ contains $\iota(c).$ So we have shown that $C_2(V_L)$
contains $\iota(c),$ and hence $S(H\otimes t^{-1})\otimes\iota(c)$ by
Lemma \ref{l3.7a}, unless $\<\a,\gamma\>\leq \<\a,\a\>$ for all
$\alpha\in L.$ Let $\Gamma$ denote the set of $\gamma\in L$ with this
latter property.

Now $\Gamma$ is a finite set. Fix a $\Z$-basis $B$ of $L$ and let
$M=\max_{\gamma\in\Gamma,\b\in B}(1-\<\beta,\gamma-\beta\>).$
{}From (\ref{12.3}) we see that $\beta(-1)^M\otimes\iota(c)\in C_2(V_L)$
for all $\beta\in B,\gamma\in\Gamma.$ Now from the above calculations,
we see that $C_2(V_L)$ contains $S^r(H\otimes t^{-1})\otimes \C\{L\}$ for
all big enough integers $r$ and also $C_2(V_L)$ contains 
$S(H\otimes t^{-1})\otimes \C\iota(c)$ for all $c\in\hat L$ such that
$\bar c\in L\backslash \Gamma.$ It follows from (\ref{12.2}) that indeed
$C_2(V_L)$ has finite codimension in $V_L.$ \qed

\begin{prop}\label{p12.6} Let $k$ be a positive integer and $\frak g$ a 
complex simple
Lie algebra. Then the \voa $L(k,0)$ associated to $\frak g$ and
$k$ satisfies condition $C_2.$
\end{prop}

\pf See [FZ] for the \voa structure of $L(k,0)$ and also the corresponding
Verma module $M(k,0).$ By definition $$M(k,0)=U(\hat \frak
g)\otimes_{U(\sum_{n=0}^{\infty}t^n\otimes \frak g+\C c)}\C\simeq
U(\sum_{n=1}^{\infty}t^{-n}\otimes\frak g)$$ (linearly). Then
$L=L(k,0)$ is the quotient of $M(k,0)$ by the maximal $\hat\frak
g$-submodule. For $a\in \frak g,$ $Y(a,z)=\sum_{n\in\Z}a(n)z^{-n-1},$
so $C_2(L)$ contains $a(-n)L$ for all $a\in\frak g$ and all $n\geq 2$
by Lemma \ref{l12.1}.
Thus $L=C_2(L)+U(t^{-1}\otimes g)1.$ It is enough to show that
$C_2(L)$ contains $a_1(-1)^{m_1}\cdots a_d(-1)^{m_d}1$ whenever
$m_i\geq 0$ and $m_1+\cdots +m_d$ is large enough; here $a_1,...,a_d$
is a basis of $\frak g.$

By Lemma 3.6 of [DLM2] we may choose that $a_i$ so that 
$$[Y(a_i,z_1),Y(a_i,z_2)]=0,\ \ Y(a_i,z)^{3k+1}=0$$
for each $i.$ Now the constant term in $Y(a_i,z)^{3k+1}1$ is equal to 
$a_i(-1)^{3k+1}1+r$ where $r$ is a sum of products of the form
$a_i(n_1)^{e_1}\cdots a_i(n_{3k+1})^{e_{3k+1}}1$ with some $n_j\leq -2.$
Since the operators $a_i(n_j)$ commute, $r\in C_2(L).$ Hence
$a_i(-1)^{3k+1}1\in C_2(L).$  We can now conclude that $C_2(L)$ contains
$a_1(-1)^{m_1}\cdots a_d(-1)^{m_d}1$ whenever $m_i\geq 3k+1$ for some
$i.$ The proposition follows immediately. \qed

\section{Applications to the Moonshine Module}
\setcounter{equation}{0}

We now apply our results to the study of the conjectures of Conway-Norton-Queen
as discussed in the Introduction. Recall that the moonshine module $V^{\natural}$ is a \voa whose automorphism group is precisely the Monster ${\Bbb M}.$
See [B1], [FLM2], [G] for details. The first author proved that $V^{\natural}$
has a unique simple module, namely $V^{\natural}$ itself, in [D2], and
in [DLM3] we showed that in fact $V^{\natural}$ is rational, that is 
every admissible $V^{\natural}$-module is completely reducible. Thus $\V$
is holomorphic. It also satisfies condition $C_2$ (Proposition \ref{p12.4}).
By Theorem \ref{t10.3} we conclude that there is a unique simple
$g$-twisted $\V$-module $\V(g)$ for each ${\Bbb M}.$ For each pair
of commuting elements $(g,h)$ in ${\Bbb M},$ recall from Section 11 that
$Z(g,h,\tau)=Z(g,h)$ is the corresponding partition function. The function
$Z(1,h)$ is precisely the graded character of $h\in{\Bbb M}$ on
$\V.$ By the results of Borcherds [B2], which confirm the original 
Conway-Norton conjecture [CN], each $Z(1,h)$ is a {\em Hauptmodul} -- in fact
the Hauptmodul conjectured in [CN]. We use this to prove the next result,
which completes the proof of Theorem \ref{thmm5}. 
\begin{th}\label{t13.1} The following hold:

(i) There is a scalar $\sigma=\sigma(g)$ such that the graded dimension
$Z(g,1,\tau)$ of $\V(g)$ is equal to $\sigma Z(1,g,S\t).$ In particular,
$Z(g,1,\t)$ is a Hauptmodul.

(ii) More generally, if a commuting pair $g,h\in{\Bbb M}$ generates a cyclic
group then $Z(g,h,\tau)$ is Hauptmodul.
\end{th}

\pf Suppose that $\<g,h\>=\<k\>.$ Then there is $\gamma\in SL(2,\Z)$ such that
$(g,h)\gamma^{-1}=(1,k).$ By Theorem \ref{t11.4} we have
$$Z(g,h,\tau)=\sigma Z(1,k,\gamma\t)$$
for some scalar $\sigma.$ Since $Z(1,k,\tau)$ is a Hauptmodul, so too
is $Z(g,h,\tau).$ If $h=1$ then we may take $k=g$ and $\gamma=S.$ 
Both parts of the Theorem now follow. \qed
\bigbreak
More is known in special cases. Huang has shown [H] that if $g$ is of
type $2B$ (in ATLAS notation) then in fact the constant $\sigma(g)$
in part (i) is equal to 1. This also follows from our results
and [FLM2]. Similarly, if $g$ is of type $2A$ it is shown in [DLM1],
on the basis of Theorem \ref{t13.1}, that again $\sigma(g)=1.$ In this case,
the precise description of $Z(g,h,\tau)$ for
$gh=hg$ and $h$ of odd order is given in [DLM1].

As discussed in [DM1] and [DM2], the uniqueness of $\V(g)$ leads to a
projective representation of the centralizer $C_{\Bbb M}(g)$ on
$\V(g).$ These are discussed in some detail, in the case that $g$ has
order $2,$ in (loc.cit.). The conjectures of Conway-Norton-Queen
state that there are (projective) representations of $C_{\Bbb M}(g)$
on suitably graded spaces such that all graded traces are either
Hauptmoduln or zero. There is no doubt that the twisted modules $\V(g)$
are the desired spaces. One would still like 
to show that $\sigma(g)=1$ in part
(i) of the theorem for all $g\in {\Bbb M},$ and to compute the
Mckay-Thompson series $Z(g,h,\tau)$ for $\<g,h\>$ {\em not} cyclic.

Finally, we discuss some aspects of correlation functions. We fix a holomorphic
vertex operator algebra $V$ satisfying condition $C_2.$ Let $Z=Z(v,(g,h),\tau)$
be the $(g,h)$ 1-point correlation function associated with $v$ and 
a pair of commuting elements $g,h\in \Aut(V)$ as defined in 
(\ref{8.4}). If $g=h=1$ we write $Z(v)$ or $Z(v,\tau)$ for this function. If
$\wt[v]=k$ then we have seen that $Z(v)$ is a generalized modular
form, holomorphic in the upper half-plane $\frak h.$ In fact, 
equation (\ref{11.5}) tells us that $Z(v)$ spans a 1-dimensional
$SL(2,Z)$-module under the action (\ref{e1.13}), that is we have
\begin{equation}\label{13.1}
Z|\gamma(v,\tau)=\sigma(\gamma)Z(v,\tau)
\end{equation}
for some character $\sigma: SL(2,Z)\to \C^*.$ We use this to prove
\begin{lem}
Let $V$ be a holomorphic vertex operator algebra which satisfies condition
$C_2$ and let $v\in V$ satisfy $\wt[v]=k.$ If $k$ is odd, then the
correlation function $Z(v,\tau)$ is identically zero.
\end{lem}

\pf We observe that the $q$-expansion of $Z(v,\tau)$ lies in 
$\C[[q^{1/3},q^{-1/3}]].$ This is because $V$ is holomorphic and so 
$8|c$ (cf. Theorem \ref{thmm2} (ii)). It follows that $T$ acts 
on $Z(v)$ as multipication be a cube root of unity, and since 
$T$ covers the abelianization of $SL(2,Z)$ then the kernel of
the character $\sigma$ has index  dividing 3. Since $S^4=id$ it follows
that $S,$ and in particular $S^2,$ lies in the kernel of $\sigma.$
Now setting $\gamma=S^2$ in (\ref{13.1}) yields 
$$Z(v,\tau)=Z|S^2(v,\tau)=(-1)^kZ(v,\tau).$$
The lemma follows. \qed

\bigskip

Now let $G$ be a finite group of automorphisms of $V.$ Then $Z(v,1,g,\tau)$
is essentially the graded trace of $o(v)g$ on $V$ for
$g\in G.$ If we choose $v$ to be the conformal vector $\tilde\omega$ of
$(V,Y[])$ then $\wt[\tilde \o]=2,$ so $Z(\tilde\o,1,g,\tau)$ is a form 
of weight 2. It is easy to describe:
\begin{lem} $Z(\tilde\o,1,g,\tau)=\frac{1}{2\pi i}
\frac{d}{d\tau}Z(1,g,\tau).$
\end{lem}

\pf One could proceed by setting $\tilde\o=\o-c/24$ and using
$Y(\omega,z)=\sum_nL(n)z^{-n-2},$ but it is 
simpler to use (\ref{m5.10}) and (\ref{m5.11}), as we may because $g\tilde\o
=\tilde\o.$  As $\tilde\o=L[-2]{\bf 1},$ the lemma follows. \qed

\bigskip

If we write $V=\oplus_nV_n,$ then of course we have 
$$Z(1,g,\tau)=q^{-c/24}\sum_{n} (\tr|_{V_n}g)q^n$$
$$Z(\tilde\o,1,g,\tau)=q^{-c/24}\sum_{n}(n-c/24) (\tr|_{V_n}g)q^n$$
and we may think of $Z(\tilde\o,1,g,\tau)$ as arising from a sequence 
of virtual characters of $G.$ That is, instead of 
``Moonshine of weight 0,'' one now has ``Moonshine of weight 2.''
This is relevant because of the work of Devoto [De] in which such 
things are interpreted as being elements of degree 2 in the
elliptic cohomology of $BG.$ 

Similarly suppose that $v\in V$ satisfies $\wt [v]=k$ with $gv=v$ for
all $g\in G.$ Then $G$ commutes with $o(v)$ in its action on each $V_n,$
so each eigenspace of the semisimple part $o(v)_s$ of $o(v)$ on $V_n$ is
a $G$-module and gives rise to a ``generalized module''
for $G$, i.e., of the form $\sum_{i}\l_n^iV_n^i$ with 
$\l_n^i\in\C$ the distinct eigenvalues of $o(v)_s$ on $V_n$ and 
$V_n^i$ the corresponding eigenspaces of $o(v)_s.$
In this way, the pair $(V,v)$ gives rise to a sequence of generalized
modules $\sum_{n,i}\l_n^iV_n^i$ for $G$ such that the corresponding
trace functions $Z(v,1,g,\tau)$ are modular forms of weight $k.$ This
is ``Moonshine of weight $k$,'' and together with the analogues 
for the twisted sectors gives rise to elements of
$Ell^kBG$ as in [De]. Actually, this is not quite what we have proved, because in [De] there are additional arithmetic requirements. It seems likely
that the appropriate conditions do hold, but that remains to be investigated.

We conclude with an illustration of some weight 4 Monstrous Moonshine 
modular forms. That is,
we take $V=V^{\natural}$ to be the Moonshine Module, we choose 
$v=L[-2]\tilde\o,$ so that $\wt [v]=4,$ and take $G$ to be the Monster.
Let $T_g(q)=Z(1,g,\tau),$ with $E_4(q)$ as in (\ref{m4.26}). Then 
\begin{eqnarray*}
& &Z(v,\tau)=12\cdot 71 E_4(q)(J(q)-240)\\
& &\ \ \ \ \ \ \ =\frac{71}{60}\left\{q^{-1}+(21\chi_1+\frac{51}{71}\chi_2)q
+(91\chi_1+\frac{701}{71}\chi_2+ \frac{221}{71}\chi_3)q^2+\cdots\right\}
\end{eqnarray*}
with $1=\chi_1,\chi_2,\chi_3,...$ the irreducible characters of the
Monster in ascending degree. Also,
$$Z(v,1,2B,\tau)=12\cdot 71 \left\{E_4(2\tau)(T_{2B}(\tau)+\frac{88}{71})
- -\frac{88}{71}E_4(\tau)\right\}$$
$$ Z(v,1,3B,\tau)=12\cdot 71 \left\{E_4(3\tau)(T_{3B}(\tau)+\frac{360}{71})
+\frac{360}{71}E_4(\tau)\right\}.$$
One can use (\ref{m5.10}) to calculate these and many other examples, of weights
4, 6, 8,....


\begin{thebibliography}{abcdefgh} 
\bibitem[AM]{AM} G. Anderson and G. Moore, Rationality in conformal
field theory, {\em Commu. Math. Phys.} {\bf 117} (1988), 441-450.
\bibitem[B1]{B1} R. E. Borcherds, Vertex algebras, Kac-Moody algebras,
and the Monster, {\em Proc. Natl. Acad. Sci. USA} {\bf 83} (1986), 3068-3071. 
\bibitem[B2]{B2} R. E. Borcherds, Monstrous moonshine and monstrous 
Lie superalgebras, {\em Invent. Math.} {\bf 109} (1992), 405-444.
\bibitem[B3]{B3} R. E. Borcherds, Modular Moonshine III, preprint.
\bibitem[BR]{BR} R. E. Borcherds and A. Ryba, Modular Moonshine II,
{\em Duke. Math. J.} {\bf 83} (1996), 435-459.
\bibitem[CN]{CN} J. H. Conway and S. P.  Norton, Monstrous Moonshine,
{\em Bull. London. Math. Soc.} {\bf 12} (1979), 308-339.
\bibitem[De]{De} J. Devoto, Equivariant cohomology and finite groups,
{\em Michigan Math. J.} {\bf 43} (1996), 3-32.
\bibitem[DGH]{DGH} L. Dixon, P. Ginsparg, J. A. Harvey, Beauty and the
beast: Super-conformal symmetry in a Monster module, {\em Commu. Math.
Phys.} {\bf 119} (1986), 221-241.
\bibitem[DHVW]{DHVW} L. Dixon, J. A. Harvey, C. Vafa and E. Witten,
String on orbifolds, {\em Nucl. Phys.} {\bf B261} (1985), 620-678;
String on orbifolds II {\em Nucl. Phys.} {\bf B274} (1986), 285-314.
\bibitem[DGM]{DGM} L. Dolan, P. Goddard and P. Montague, 
Conformal field theory of twisted vertex operators,
{\em Nucl. Phys.}  {\bf B338} (1990), 529.
\bibitem[D1]{D1} C. Dong, Vertex algebras associated with 
even lattices, {\it J. Algebra} {\bf 160} (1993), 245-265.
\bibitem[D2]{D2} C. Dong,  Representations of the moonshine module 
vertex operator algebra, {\em Contemporary Math.} {\bf 175} (1994).
\bibitem[DL]{DL} C. Dong and J. Lepowsky, Generalized Vertex
Algebras and Relative Vertex Operators, {\em Progress in Math.} Vol. 112,
Birkh\"{a}user, Boston, 1993.
\bibitem[DLM1]{DLM1} C. Dong, H. Li and G. Mason, Some twisted modules for the moonshine vertex operator algebras,
{\em Contemp. Math.} {\bf 193} (1996), 25-43.
\bibitem[DLM2]{DLM2} C. Dong, H. Li and G. Mason, Regularity of rational vertex operator algebras, {\em Advances. in Math.}
to appear, q-alg/9508018.
\bibitem[DLM3]{DLM3} C. Dong, H. Li and G. Mason,  Twisted representations of
vertex operator algebras, {\em Math. Annalen,} to appear, q-alg/9509005.
\bibitem[DM1]{DM1} C. Dong and G. Mason, 
Nonabelian orbifolds and the boson-fermion correspondence,
{\em Commu. Math. Phys.} {\bf 163} (1994), 523-559.
\bibitem[DM2]{DM2} C. Dong and G. Mason, 
Vertex operator algebras and moonshine: A survey, {\it Advanced Studies
in Pure Math.} {\bf 24} (1996), 101-136.
\bibitem[DMZ]{DMZ} C. Dong, G. Mason and Y. Zhu, 
Discrete series of the 
Virasoro algebra and the moonshine module, {\em Proc. Symp. Pure. Math., American Math. Soc.} {\bf 56} II (1994), 295-316.
\bibitem[EZ]{EZ} M. Eichler and D. Zagier, On the Theory of Jacobi Forms I,
{\em Progress in Math.} Vol. 55, Birkh\"{a}user, Boston, 1985.
\bibitem[FF]{FF}
B. L. Feigin and D. B. Fuchs,Verma modules over the Virasoro algebra, Lecture 
Notes in Math., Vol. 1060, Springer-Verlag, Berlin and New York, 1984.
\bibitem[FHL]{FHL} I. B. Frenkel, Y. Huang and J. Lepowsky, On
axiomatic approaches to vertex operator algebras and modules,
{\it Memoirs American Math. Soc.} {\bf 104}, 1993.
\bibitem[FLM1]{FLM1} I. B. Frenkel, J. Lepowsky and A. Meurman,
Vertex operator calculus, in: {\em Mathematical Aspects of String 
Theory, Proc. 1986 Conference, San Diego.} ed. by S.-T. Yau, World
Scientific, Singapore, 1987, 150-188.
\bibitem[FLM2]{FLM2} I. B. Frenkel, J. Lepowsky and A. Meurman, 
Vertex Operator Algebras and the Monster, {\em Pure and Applied
Math.,} Vol. {\bf 134}, Academic Press, 1988.
\bibitem[FZ]{FZ} I. B. Frenkel and Y. Zhu, Vertex operator algebras 
associated to  representations of affine and Virasoro algebras, {\em Duke
 Math. J.} {\bf 66} (1992), 123-168.
\bibitem[G]{G} R. Griess, The Friendly Giant, {\em Invent. Math.}
{\bf 69} (1982), 1-102.
\bibitem[HMV]{HMV} J. Harvey, G. Moore and C. Vafa, Quasicrystalline 
compactification, {\em Nucl. Phys.}  {\bf B304} (1988), 269-290.
\bibitem[H]{H} Y. Huang, A non-meromorphic extension of the
moonshine module vertex operator algebra, {\em Contemporary Math.}
{\bf 193} (1996), 123-148.
\bibitem[I]{I} E. Ince, Ordinary Differential Equations, Dover Publications,
Inc., 1956.
\bibitem[KP]{KP} V. Kac and D. Peterson, Infinite-dimensional Lie
algebras, theta functions and modular forms, {\em Advances in Math.} {\bf
53} (1984), 125-264.
\bibitem[KM]{KM} M. Knopp and G. Mason, Generalized modular forms, 
in preparation.
\bibitem[K]{K} N. Koblitz, Introduction to Ellipitic Curves and
Modular Forms, Springer-Verlay, New York, 1984.
\bibitem[La]{La} S. Lang, Elliptic Functions, G.T.M. Vol. {\bf 112},
2nd. ed., Springer-Verlag, New York, 1987.
\bibitem[Le]{Le} J. Lepowsky, Calculus of twisted vertex operators, 
{\em Proc. Natl. Acad Sci. USA} {\bf 82} (1985), 8295-8299.
\bibitem[Li]{Li} H. Li, Local systems of vertex operators, vertex superalgebras
and modules, {\em J. Pure Appl. Alg.} {\bf 109} (1996) 143-195.
\bibitem[M]{M} P. Montague, Orbifold constructions and the classification
of self-dual $c=24$ conformal field theory, {\em Nucl. Phys.}  {\bf B428}
(1994), 233-258.
\bibitem[MS]{MS} G. Moore and N. Seiberg, Classical and quantum 
conformal field theory, {\em Comm. Math. Phys.} {\bf 123} (1989), 177-254.
\bibitem[N]{N}  S. Norton, Generalized moonshine, 
{\em Proc. Symp. Pure. Math., American Math. Soc.} {\bf 47} (1987), 208-209.
\bibitem[Q]{Q} L. Queen, Some relations between finite groups, Lie groups
and modular functions, Ph.D. Thesis, University of Cambridge, 1980.
\bibitem[Ra]{Ra} H. Rademacher, Topics in Analytic Number Theory, Springer-Verlag, New York, 1973.
\bibitem[Ry]{Ry} A. Ryba, Modular Moonshine? {\em Contemp. Math.} {\bf 193} (1996), 307-336.
\bibitem[S]{S} A. N. Schellekens, Meromorphic $c=24$ Conformal Field 
Theories, {\it Commu. Math. Phys.} {\bf 153} (1993), 159.
\bibitem[T1]{T1} M. Tuite, Monstrous moonshine from orbifolds, {\em
Commu. Math. Phys.} {\bf 146} (1992), 277-309.
 \bibitem[T2]{T2} M. Tuite, On the relationship between monstrous moonshine 
and the uniqueness of the moonshine module, {\em
Commu. Math. Phys.} {\bf 166} (1995), 495-532.
\bibitem[T3]{T3} M. Tuite, Generalized moonshine and abelian orbifold 
constructions, {\em Contemp. Math.} {\bf 193} (1996), 353-368.
\bibitem[V]{V} C. Vafa, Modular invariance and discrete torsion
on orbifolds,  {\em Nucl. Phys.}  {\bf B273} (1986), 592.
\bibitem[Wa]{Wa} W. Wang, Rationality of Virasoro vertex operator algebras, {\it Duke
Math. J. IMRN}, {\bf Vol. 71}, No. 1 (1993), 197-211.
\bibitem[Wo]{Wo} K. Wohlfahrt, An extension of F. Klein's level concept,
{\em Illinois J. Math.} {\bf 8} (1964), 529-535.
\bibitem[Z]{Z} Y. Zhu, Modular invariance of characters of vertex operator algebras,
{\em J. Amer, Math. Soc.} {\bf 9} (1996), 237-302.
\end{thebibliography}
\end{document}